\documentclass[sn-chicago,iicol]{sn-jnl}

\jyear{2022}%

\theoremstyle{thmstyleone}%
\theoremstyle{thmstyletwo}%

\theoremstyle{thmstylethree}%
\newtheorem{remark}{Remark}%

\usepackage{amsmath}
\usepackage{mathtools}
\usepackage{amssymb}
\usepackage{amsfonts}
\usepackage{accents}
\usepackage{hyperref}
\usepackage{makecell}
\usepackage{bbm}
\hypersetup{
	colorlinks=true,
	linkcolor=blue,
	filecolor=magenta,      
	urlcolor=cyan,
}
\usepackage{subcaption}
\usepackage{titlesec}
\usepackage[linesnumbered,ruled,vlined,boxed,commentsnumbered]{algorithm2e}

\usepackage{etoolbox}

\newcommand\munderbar[1]{%
	\underaccent{\bar}{#1}}

\newcommand{\blds}[1]{\mbox{\Large $#1$}}

\begin{document}

\title[Article Title]{An isogeometric finite element formulation for frictionless contact of Cosserat rods with unconstrained directors}


\author*[1]{\fnm{Myung-Jin} \sur{Choi}}\email{choi@lbb.rwth-aachen.de}
\author[1]{\fnm{Sven} \sur{Klinkel}}\email{klinkel@lbb.rwth-aachen.de}

\author[2,3,4]{\fnm{Roger} \sur{A. Sauer}}\email{sauer@aices.rwth-aachen.de}

\affil[1]{\orgdiv{Chair of Structural Analysis and Dynamics}, \orgname{RWTH Aachen University}, \orgaddress{\street{Mies-van-der-Rohe Str. 1}, \city{Aachen}, \postcode{52074}, \country{Germany}}}

\affil[2]{\orgdiv{Aachen Institute for Advanced Study in Computational Engineering Science (AICES)}, \orgname{RWTH Aachen University}, \orgaddress{\street{Templergraben 55}, \city{Aachen}, \postcode{52062}, \country{Germany}}}

\affil[3]{\orgdiv{Faculty of Civil and Environmental Engineering}, \orgname{Gdańsk University of Technology}, \orgaddress{\street{ul. Narutowicza 11/12}, \city{Gdańsk}, \postcode{80-233}, \country{Poland}}}

\affil[4]{\orgdiv{Department of Mechanical Engineering}, \orgname{Indian Institute of Technology Guwahati}, \orgaddress{\city{Guwahati}, \postcode{781039}, \state{Assam}, \country{India}}}


\abstract{This paper presents an isogeometric finite element formulation for nonlinear beams with impenetrability constraints, based on the kinematics of Cosserat rods with unconstrained directors. The beam cross-sectional deformation is represented by director vectors of an arbitrary order. For the frictionless lateral beam-to-beam contact, a surface-to-surface contact algorithm combined with an active set strategy and a penalty method is employed. The lateral boundary surface of the beam is parameterized by its axis and cross-sectional boundary curves with NURBS basis functions having at least $C^2$-continuity, which yields a continuous surface metric and curvature for the closest point projection. Three-dimensional constitutive laws of hyperelastic materials are considered. Several numerical examples verify the accuracy and efficiency of the proposed beam contact formulation in comparison to brick element solutions. The lateral contact pressure distribution of the beam formulation is in excellent agreement with the contact pressure of the brick element formulation while requiring much less degrees-of-freedom.}



\maketitle
\section{Introduction}\label{sec1}
The simulation of interacting rods or rod-like bodies has been investigated across many applications, including wire strands \citep{menard2021solid}, cables \citep{bajas2010numerical}, biopolymer networks \citep{cyron2012numerical}, woven fabrics \citep{goyal2005nonlinear,durville2010simulation}, entangled fibrous materials \citep{rodney2016reversible}, DNA supercoiling \citep{lillian2011electrostatics}, and deformations in adhesive microstructures \citep{sauer2009multiscale}. In such examples the large number of bodies and their contact interactions typically cause significant computational costs, which calls for the development of efficient and accurate beam and beam-to-beam contact formulations.

A \textit{beam} in solid mechanics refers to a dimensionally reduced model of a three-dimensional slender body based on a suitable kinematic assumption. In Cosserat (or \textit{directed}) rod theory, a material point position in the current configuration, in three-dimensional space, is given by \citep{naghdi1981finite}
\begin{equation}
	\label{beam_kin_expansion_1st}
	{{\boldsymbol{x}}} = {\boldsymbol{\varphi }}(s) + {\zeta^1}{{\boldsymbol{d}}_1}(s) + {\zeta^2}{{\boldsymbol{d}}_2}(s),
\end{equation}
where $\boldsymbol{\varphi}(s)$ denotes the position vector of the current axis\footnote{The \textit{axis} is typically defined by a spatial curve connecting (mass) centroids of the cross-sections, which is often called \textit{line of centroids}.}, and $\boldsymbol{d}_\gamma(s)\,(\gamma\in\left\{1,2\right\})$ represent two director vectors spanning the planar cross-section. $s$ denotes the arc-length coordinate along the initial axis, and $\zeta^1$ and $\zeta^2$ denote two transverse coordinates of the cross-section. This kinematic expression is sufficiently general to encompass constrained theories like Euler-Bernoulli, and Timoshenko beams \citep{nordenholz1997steady}. One can refer to \citet{meier2019geometrically} and the references therein for various nonlinear beam formulations with orthonormality constraints on the two directors, which are typically satisfied by a parameterization using an orthogonal tensor with three rotational degrees-of-freedom (DOFs) in space, or constraints of shear-free deformations. There are several previous works employing Eq.\,(\ref{beam_kin_expansion_1st}) in its \textit{unconstrained} form, which allows to use nine DOFs per cross-section, for example, \citet{rhim1998vectorial}, \citet{durville2012contact} and \citet{choi2021isogeometric}. The same DOFs appear in the brick element formulation with six nodes in \citet{schweizerhof2014solid} and \citet{konyukhov2018consistent} that use polar coordinates for the elliptical cross-section, combined with a linear approximation along the longitudinal direction. 

The \textit{extensibility} of the two directors in Eq.\,(\ref{beam_kin_expansion_1st}) gives several advantages including an additive configuration update procedure, in-plane cross-sectional deformations, and a straightforward implementation of three-dimensional constitutive laws. Such formulations have been considered in several works. For example, \citet{frischkorn2013solid} presented a brick element formulation combined with an enhanced assumed strain (EAS) method, an assumed natural strain (ANS) method, and a reduced integration method in order to alleviate locking. \citet{wackerfuss2009mixed} developed a mixed variational formulation incorporating transverse normal strains, where an arbitrary three-dimensional constitutive laws can be easily implemented. The first order expression of Eq.\,(\ref{beam_kin_expansion_1st}) in terms of the transverse coordinates leads to constant in-plane cross-sectional strains, which suffers from artificial increase of bending stiffness for nonzero Poisson's ratio. In order to circumvent this Poisson locking, two quadratic order terms are additionally introduced in \cite{coda2009solid}; however, this formulation still suffers from locking due to the missing bilinear terms. An EAS method of enriching all linear in-plane Green-Lagrange strain components was verified to effectively alleviate Poisson locking in \citet{choi2021isogeometric}.

It was shown in \citet{naghdi1989significance} that the consideration of transverse normal strains of the cross-section significantly contributes to correctly predict the contact force distribution.
In \citet{olga2018contact}, a penalty parameter depending on the amount of penetration due to the Hertz theory of elastic contact is utilized within a small strain range in order to consider cross-sectional strains. The beam kinematics of Eq.\,(\ref{beam_kin_expansion_1st}) combined with a simplified constitutive equation to alleviate Poisson locking was employed in \citet{durville2012contact} for frictional beam-to-beam contact. 

In order to more accurately capture cross-sectional strains in a consistent way from three-dimensional elasticity, Eq.\,(\ref{beam_kin_expansion_1st}) can be generalized to have an arbitrary order of approximation in the transverse directions by using the series expansion \citep{antman1966dynamical}
\begin{equation}
	\label{exact_taylor_Nth}
	{{\boldsymbol{x}}} = \sum\limits_{p = 0}^N  {\sum\limits_{q = 0}^p {{{({\zeta ^1})}^{p - q}}\,{{({\zeta ^2})}^q}\,{{\boldsymbol{d}}^{(p - q,q)}}(s)} },
\end{equation}
with the directors
\begin{equation}
	\label{def_taylor_deriv_xt}
	{{\boldsymbol{d}}^{(m,n)}} \coloneqq \frac{1}{{m!{\,}n!}}\left(\left. {\frac{{{\partial ^{m+n}}{{\boldsymbol{x}}}}}{{\partial {{({\zeta ^1})}^{m}}\,\partial {{({\zeta ^2})}^n}}}}\right\rvert_{\zeta^1=\zeta^2=0}\right),
\end{equation}
where $m$ and $n$ are nonnegative integers, and $N$ is a positive integer representing the order of approximation in transverse direction. It is noted that Eq.\,(\ref{beam_kin_expansion_1st}) is a special case of Eq.\,(\ref{exact_taylor_Nth}) with $N=1$, ${\boldsymbol{\varphi }}\equiv{{\boldsymbol{d}}^{(0,0)}}$, ${\boldsymbol{d}}_1\equiv{\boldsymbol{d}}_{}^{(1,0)}$, and ${\boldsymbol{d}}_2\equiv{\boldsymbol{d}}_{}^{(0,1)}$. Here and hereafter, we often omit the argument $s$ in directors for brevity. The higher order kinematics of Eq.\,(\ref{exact_taylor_Nth}) with $N\ge2$ was employed in \cite{moustacas2021higher} to obtain a homogenized model of fiber bundles in contact with a rigid surface. Compared to existing works, our beam and contact formulations introduce the following novelties:
\begin{enumerate}
	\item Unconstrained directors enable an efficient and accurate description of cross-sectional strains. Further we verify the beam solutions by comparison with those using brick elements.
	\item We present a frictionless beam contact formulation considering cross-sectional strains, based on a Gauss point-to-surface contact algorithm. The impenetrability constraint is enforced using a penalty method combined with an active set method. We verify the contact pressure distribution obtained from the beam contact formulation by comparison with the results of the brick element formulation.
	\item A parameterization of the lateral boundary surface with NURBS basis functions with at least $C^2$-continuity in both longitudinal and circumferential directions enables to have continuous surface metric and curvature within the local Newton-Raphson iteration for the closest point projection.
\end{enumerate}

For very slender beams, a contact force on the lateral surface is often assumed to act directly on the center axis. Such \textit{curve-to-curve} contact formulation can be divided into the following two categories \citep{meier2016finite}. First, a \textit{point-to-point} contact formulation considers a discrete contact force between two interacting bodies. It is typically based on a \textit{bilateral} closest point projection, e.g, see \citet{wriggers1997contact}, and provides an efficient formulation, since the evaluation of the contact integral is not required. However, it suffers from the non-uniqueness of solution in the closet point projection if the intersection angle between the curves is small, e.g., in the case of two nearly parallel beams. Further discussions on the uniqueness of solution in the closest point projection can be found in \cite{konyukhov2008solvability} and \cite{meier2016finite}. One can also refer to \cite{weeger2017isogeometric} for a relevant discussion in the framework of isogeometric collocation method. Second, a \textit{line-to-line} (or Gauss point-to-axis) contact formulation considers a distributed contact force. It is also applicable to small intersection angles between curves, although it requires larger computational cost, compared to the point-to-point contact formulation, as the angle decreases. In order to consider an arbitrary intersection angle in a more efficient way, those two formulations are combined in \cite{meier2017unified} with a regularized transition. The \textit{curve-to-curve} contact formulation is efficient due to its one-dimensional contact search; however, for beams having low to moderate slenderness ratios, considering a surface load formulation rather than applying an equivalent load to the axis becomes crucial due to the additional effects like external moments and cross-sectional stretching \citep{choi2021isogeometric}. Further the closest point projection to the center axis might lead to an error in the calculation of contact forces since the cross-section is not always orthogonal to the axis due to transverse shear deformations \citep{sauer2014geometrically}. Therefore, in this paper, we present a \textit{surface-to-surface} contact formulation with a closest point projection to the lateral beam surface, based on a parameterization of the lateral surface with at least $C^2$-continuity. \citet{neto2016master} presented a finite element formulation of point-to-point frictionless contact on the lateral surface based on the bilateral closest point projection, which still suffers from the non-uniqueness of solution for small intersection angles. 

A seamless integration of geometry and analysis is developed, in \cite{hughes2005isogeometric}, by employing non-uniform rational B-splines (NURBS) basis functions for the spatial discretization of the solution field as well as the geometry, which is known as isogeometric analysis (IGA). The higher-order inter-element continuity in IGA yields smooth contact pressure distributions \citep{lu2011isogeometric, temizer2011contact}, and enables to reduce numerical instabilities and oscillations in classical node-to-segment (NTS) and Gauss-point-to-segment contact algorithms associated with kinks between elements \citep{matzen2013point,sauer2013local}. Thus, surface smoothing \citep{padmanabhan2001framework, wriggers2001smooth, stadler2003cn} and surface enrichment \citep{sauer2013local, corbett2014nurbs}, developed for classical finite element discretizations based on Lagrange polynomials are intrinsically captured by IGA. In this paper, we parameterize the initial lateral surface of the beam by the center axis curve and the cross-section's boundary curve defined in the plane spanned by two orthonormal directors. Especially, we employ an \textit{unclamped} knot vector in the cross-section's boundary curve in order to have higher order continuity, which gives us continuous metric and curvature components for the surface that are necessary for robustness of the local Newton-Raphson iteration in the closest point projection. 

The remainder of this paper is organized as follows. In Section \ref{nth_beam_crod}, we present the beam kinematics based on Cosserat rod theory, and the parameterization of the lateral boundary surface of the rod. In Section \ref{sec_frictionless_contact}, the frictionless beam-to-beam contact formulation is presented. In Section \ref{sec_spatial_disc_iga}, an isogeometric finite element discretization of the beam and contact formulations is presented. In Section \ref{num_ex_contact}, we present several numerical examples of beam-to-rigid body and beam-to-beam contact. In three appendices, we present detailed algorithms for the beam and its contact formulation, and supplementary information of the numerical examples.
\section{Cosserat rods with deformable cross-section}
\label{nth_beam_crod}
In this paper, in order to capture the cross-sectional strains more accurately, the isogeometric finite element formulation of \cite{choi2021isogeometric} is extended to incorporate an arbitrary order of approximation in the transverse directions, i.e., $N$, based on the kinematics of Eq.\,(\ref{exact_taylor_Nth}). For $N\ge{2}$, since the cross-section can represent at least linear in-plane strains properly, one does not require any special treatment to alleviate Poisson locking. For $N=1$, we use the EAS method in \cite{choi2021isogeometric}.
\subsection{Beam kinematics: Cosserat rod theory}
\label{nth_beam_kin}
\begin{figure} \centering	
	\includegraphics[width=0.9975\linewidth]{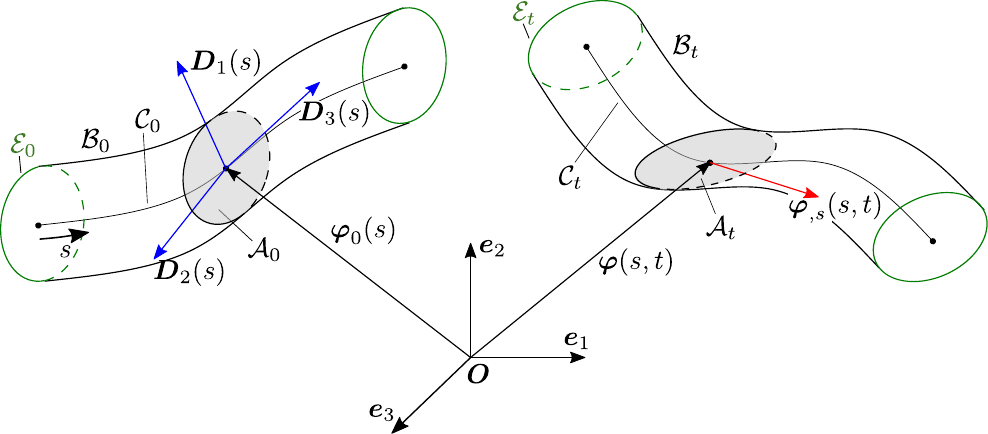}
	\caption{A schematic illustration of the beam kinematics in the initial and current configurations. Note that the initial cross-section ($\mathcal{A}_0$) is assumed planar; however, the current one ($\mathcal{A}_t$) is not always planar due to out-of-plane deformations for $N\geq2$ in Eq.\,(\ref{exact_taylor_Nth}).}
	\label{beam_kin}
\end{figure}

\begin{figure}
	\centering
	\includegraphics[width=0.85\linewidth]{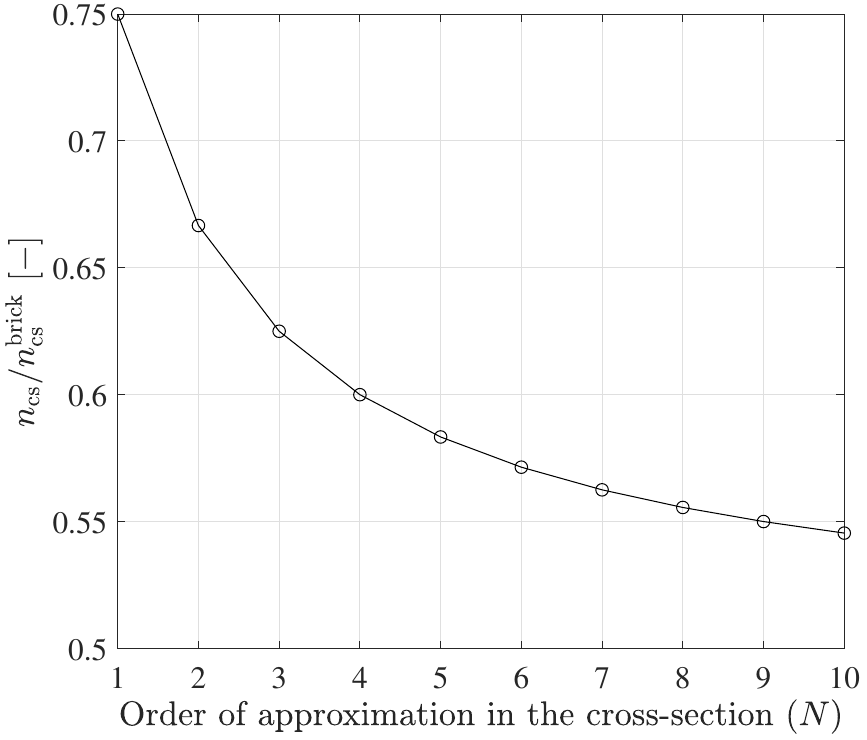}
	\caption{Comparison of the numbers of DOFs per cross-section for a beam element and a single brick element.}
	\label{beam_compare_ndof_cs}
\end{figure}

The initial configuration of a beam is typically described by a family of \textit{cross-sections} whose (mass) centroids\footnote{In this paper, we assume a constant mass density, so that the mass centroid coincides with the geometrical centroid.} are connected by a spatial curve called the \textit{line of centroids} or the \textit{initial axis}. For an initial (undeformed) axis, we consider a spatial curve $\mathcal{C}_0$ parameterized by a coordinate ${\xi}\in{\Bbb{R}^1}$, i.e., ${\mathcal{C}_0}:\,{\xi} \to {{\boldsymbol{\varphi }}_0}({\xi}) \in {{\Bbb{R}}^3}$. We reparameterize the curve as ${\mathcal{C}_0}:\,s \to {{\boldsymbol{\varphi }}_0}(s) \in {{\Bbb{R}}^3}$ by an arc-length parameter $s \in \left[ {0,L} \right] \subset {{\Bbb{R}}^1}$, where $L$ represents the length of the initial axis. The arc-length coordinate is defined by the mapping 
\begin{equation}
	\label{s_xi_relation}
	s({\xi}) \coloneqq \int_0^{\xi}  {{\left\| {{{\boldsymbol{\varphi}}_{0,{\tilde \xi}}}({\tilde \xi})} \right\|}\mathrm{d}\tilde \xi },
\end{equation}
and the Jacobian of the mapping is defined as ${\tilde j}\coloneqq{\mathrm{d}s/\mathrm{d}{\xi}}={\left\| {{{\boldsymbol{\varphi}}_{0,{\xi}}}({\xi})} \right\|}$. The arc-length coordinate enables to simplify the subsequent expressions by $\left\| {{{\boldsymbol{\varphi }}_{0,s}}} \right\| = 1$. Here are hereafter, we often use $s\equiv s(\xi)$ for brevity, and $(\bullet)_{,s}$ denotes the partial differentiation with respect to the arc-length parameter $s$. The initial cross-section domain $\mathcal{A}_0\subset {\Bbb{R}^2}$ is spanned by two \textit{initial directors} ${{\boldsymbol{D}}_{\gamma}}(s) \in {\Bbb R}^3$ $\left(\gamma\in\left\{1,2\right\}\right)$, which are orthonormal, and aligned with the principal directions of the second moment of inertia of the cross-section. Further, it is assumed that, in the initial configuration, the cross-section is orthogonal to the initial axis; thus, we simply obtain ${\boldsymbol{D}}_3(s)\coloneqq{\boldsymbol{\varphi}}_{0,s}(s)$. The current (deformed) configuration of the axis is defined by the spatial curve ${{\mathcal{C}}_t}:\,s \to {{\boldsymbol{\varphi }}}(s,t) \in {{\Bbb{R}}^3}$, where $t\in{\Bbb R}^{+}$ denotes time. It should be noted that the current axis does not always pass through the geometrical centroid of the cross-section if the order of approximation in the transverse direction, i.e., $N$ in Eq.\,(\ref{exact_taylor_Nth}) is greater than 1, see Remark \ref{remark_cent_axis_deform}. We define $\left\{ {{{\boldsymbol{e}}_1},{{\boldsymbol{e}}_2},{{\boldsymbol{e}}_3}} \right\}$ as a standard Cartesian basis in ${{\Bbb R}^3}$, and ${\boldsymbol{e}^i}\equiv{\boldsymbol{e}_i}$, $i\in\left\{1,2,3\right\}$. Fig.\,\ref{beam_kin} schematically illustrates the above kinematic description of the initial and current beam configurations.
\begin{figure}
	\centering	
	\includegraphics[width=0.75\linewidth]{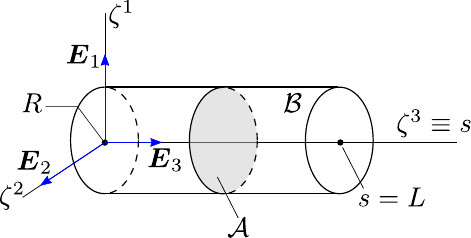}
	\caption{An example of the reference domain $\mathcal{B}$ in the case of circular cross-section with radius $R$.}
	\label{beam_ref_domain_cir}
\end{figure}
We define a \textit{reference domain} ${\mathcal{B}}\coloneqq (0,L) \times {\mathcal{A}}$, where ${\mathcal{A}}$ denotes the open domain of coordinates $\zeta^1$ and $\zeta^2$. For example, for a circular cross-section with radius $R$ we have $\mathcal{A} \coloneqq \left\{ {\left. {\left({\zeta^1},{\zeta ^2}\right)} \right\rvert{{\left( {{\zeta ^1}} \right)}^2} + {{\left( {{\zeta ^2}} \right)}^2} < R^2} \right\}$, see Fig.\,\ref{beam_ref_domain_cir} for an illustration. The location of each point in the reference domain is expressed in terms of the coordinates ${\zeta^1}$, ${\zeta^2}$, and ${\zeta^3}$ in the standard Cartesian basis in ${{\Bbb R}^3}$ denoted by ${{\boldsymbol{E}}_1}$, ${{\boldsymbol{E}}_2}$, and ${{\boldsymbol{E}}_3}$, where we use $\zeta^3\equiv s$. We then define two mappings from the reference domain to the initial configuration $\mathcal{B}_0$ and to the current configuration $\mathcal{B}_t$, respectively, by ${{\boldsymbol{X}}}:{\mathcal{B}} \to {{\mathcal{B}}_0}$ and ${{\boldsymbol{x}}}:{\mathcal{B}} \to {{\mathcal{B}}_t}$. The deformation from the initial to the current configuration is then expressed by the mapping 
\begin{equation}
	{{\boldsymbol{\Phi }}_t} \coloneqq {{\boldsymbol{x}}} \circ {{\boldsymbol{X}}}^{ - 1}:{{\mathcal{B}}_0} \to {{\mathcal{B}}_t}.
\end{equation}
We assume a smooth cross-section boundary such that the initial boundary surface ${\mathcal{S}}_0\equiv\partial{\mathcal{B}_0}$ consists of the lateral surface ${\mathcal{S}}^{\mathrm{L}}_0$, the cross-sections ${\mathcal{A}}_0$ at ends, and their interface edges ${\mathcal{E}_0}\equiv{\partial{{\mathcal{S}}^{\mathrm{L}}_0}}\equiv{\partial{\left.{{\mathcal{A}}_0}\right\rvert}_{s \in \{0,L\} }}$, i.e., ${{\mathcal{S}}_0} = {{\mathcal{S}}}_0^{\rm{L}} \cup {\left. {{{\mathcal{A}}_0}} \right\rvert_{s \in \{ 0,L\} }} \cup {\mathcal{E}_0}$. Further, in the current configuration, the boundary surface $\mathcal{S}_t\equiv\partial{\mathcal{B}_t}$ consists of the lateral surface ${\mathcal{S}}^{\mathrm{L}}_t$, the cross-sections ${\mathcal{A}}_t$ at ends, and their interface edges ${\mathcal{E}_t}\equiv{\partial{{\mathcal{S}}^{\mathrm{L}}_t}}\equiv{\partial{\left.{{\mathcal{A}}_t}\right\rvert}_{s \in \{0,L\} }}$, i.e., ${{\mathcal{S}}_t} = {{\mathcal{S}}}_t^{\rm{L}} \cup {\left. {{{\mathcal{A}}_t}} \right\rvert_{s \in \{ 0,L\} }} \cup {\mathcal{E}_t}$, see Fig.\,\ref{beam_kin} for an illustration. The initial position vector of any point of the beam is given as 
\begin{equation}
	\label{beam_init_config_kin}
	{\boldsymbol{X}}({\zeta^1},{\zeta^2},{\zeta^3}) = {\boldsymbol{\varphi }}_0({\zeta^3}) + {\zeta^\gamma }{{\boldsymbol{D}}_\gamma }({\zeta^3}).
\end{equation}
Here and hereafter, unless otherwise stated, repeated Latin indices like $i$ and $j$ imply summation over $1$ to $3$, and repeated Greek indices like $\alpha$, $\beta$\, and $\gamma$ imply summation over $1$ to $2$. It should be noted that the coordinates $\zeta^\gamma\,(\gamma\in\left\{1,2\right\})$ are chosen to have dimensions of length, such that the initial directors $\boldsymbol{D}_\gamma$ are dimensionless. Initial covariant base vectors are obtained by ${{\boldsymbol{G}}_{i}} \coloneqq \partial {{\boldsymbol{X}}}/\partial {\zeta ^i}$ ($i\in\left\{1,2,3\right\}$), so that we have
\begin{equation}\label{beam_th_str_init_cov_base}
	\left\{ \begin{array}{l}
		\begin{aligned}
			{{\boldsymbol{G}}_1}({{\zeta}^1},{{\zeta}^2},{{\zeta}^3}) &= {{\boldsymbol{D}}_1}(s),\\
			{{\boldsymbol{G}}_2}({{\zeta}^1},{{\zeta}^2},{{\zeta}^3}) &= {{\boldsymbol{D}}_2}(s),\\
			{{\boldsymbol{G}}_3}({{\zeta}^1},{{\zeta}^2},{{\zeta}^3}) &= {{\boldsymbol{D}}_3}(s) + {{\zeta ^\gamma }{{\boldsymbol{D}}_{\gamma ,s}}(s)}.\\
		\end{aligned}
	\end{array} \right.
\end{equation}
Further we define $j_0$ as the Jacobian of the mapping ${\boldsymbol{X}}\big(\zeta^1,\zeta^2,\zeta^3\big)$ such that the corresponding infinitesimal volume in the domain $\mathcal{B}_0$ can be expressed by \citep{choi2021isogeometric}
\begin{equation}
	\mathrm{d}\mathcal{B}_0={j_0}\,{\mathrm{d}\zeta^1}\,{\mathrm{d}\zeta^2}\,{\mathrm{d}s}
\end{equation}
with 
\begin{equation}
j_0=\big({\boldsymbol{G}_1}\times{\boldsymbol{G}_2}\big)\cdot{\boldsymbol{G}_3}.
\end{equation}
\begin{remark}
	We consider a pair of nonnegative integers $(p,q)$ in Eq.\,(\ref{exact_taylor_Nth}), which belong to
	\begin{equation}
		\label{def_cn}
		{C_n} = \left\{ {\left. {\left( {p,q} \right) \in {Z^*} \times {Z^*}} \right\rvert p + q = n} \right\},
	\end{equation}
	where ${\Bbb{Z}^+}$ denotes the set of positive integers, and ${Z^*}\coloneqq\{0\}\cup{\Bbb{Z}^+}$, and the number of elements in $C_n$ is $\left\lvert {{C_n}} \right\rvert = n + 1$. In this paper, $n$ of Eq.\,(\ref{def_cn}) is called the \textit{order} of directors, for examples, $\boldsymbol{\varphi}\equiv{\boldsymbol{d}^{(0,0)}}$ is called the 0th order director, and $\boldsymbol{d}_1\equiv{\boldsymbol{d}^{(1,0)}}$ and $\boldsymbol{d}_2\equiv{\boldsymbol{d}^{(0,1)}}$ are called the 1st order directors. Then, for an $N$-th order of approximation in the transverse directions, the number of directors in each cross-section is
	\begin{equation}
		n_\mathrm{dir}\coloneqq\sum\limits_{n = 0}^N {\left\lvert {{C_n}} \right\rvert}  = \frac{{(N + 1)(N + 2)}}{2},
	\end{equation}
	and the number of DOFs in each cross-section is simply obtained by $n_\mathrm{cs}=3{n_\mathrm{dir}}$. For a single brick element of order $N$, the number of nodes in each cross-section is typically $n_\mathrm{node}=(N+1)^2$, and then the number of DOFs is $n^\mathrm{brick}_\mathrm{cs}=3n_\mathrm{node}$. Fig.\,\ref{beam_compare_ndof_cs} compares the number of DOFs per cross-section in the beam and brick elements. It is shown that for the same order of approximation in the cross-section, the beam formulation uses less DOFs due to the symmetry with respect to the axis $\zeta^1=0$ or $\zeta^2= 0$ in the kinematic assumption of Eq.\,(\ref{exact_taylor_Nth}).
\end{remark}
\noindent Eq.\,(\ref{exact_taylor_Nth}) can be rewritten in the compact form \citep{Moustacas2019}
\begin{equation}
	\label{compact_xt_pi}
	{{\boldsymbol{x}}} = {\boldsymbol{\Pi }}({\zeta^1},{\zeta^2})^\mathrm{T}{\boldsymbol{q}}(\zeta^3),
\end{equation}
with
\begin{equation}
	\label{def_pi_operator_assemb}
	{\boldsymbol{\Pi }}({\zeta^1},{\zeta^2})\coloneqq{\left[ {\begin{array}{*{20}{c}}
				{{{\boldsymbol{\pi }}_0}({\zeta^1},{\zeta^2})}\\
				{{{\boldsymbol{\pi }}_1}({\zeta^1},{\zeta^2})}\\
				\vdots \\
				{{{\boldsymbol{\pi }}_N}({\zeta^1},{\zeta^2})}
		\end{array}} \right]_{{n_\mathrm{cs}}\times 3}},
\end{equation}
and
\begin{equation}
	\boldsymbol{q}(\zeta^3)\coloneqq{\left\{ {\begin{array}{*{20}{c}}
			{{{\boldsymbol{q}}_0}(\zeta^3)}\\
			{{{\boldsymbol{q}}_1}(\zeta^3)}\\
			\vdots \\
			{{{\boldsymbol{q}}_N}(\zeta^3)}
	\end{array}} \right\}_{{{n_\mathrm{cs}}\times 1}}},	
\end{equation}
where we define the operator
\begin{equation}
	\label{def_pi_operator}
	{{\boldsymbol{\pi }}_n}({\zeta ^1},{\zeta ^2}) \coloneqq {\left[ {\begin{array}{*{20}{c}}
				{{{({\zeta ^1})}^n}{{({\zeta ^2})}^0}{{\boldsymbol{1}}_{3}}}\\
				{{{({\zeta ^1})}^{n - 1}}{{({\zeta ^2})}^1}{{\boldsymbol{1}}_{3}}}\\
				\vdots \\
				{{{({\zeta ^1})}^0}{{({\zeta ^2})}^n}{{\boldsymbol{1}}_{3}}}
		\end{array}} \right]},
\end{equation}
$n\in\left\{0,1,...,N\right\}$, with
\begin{equation}
	{{\boldsymbol{q}}_n}(\zeta^3 ) \coloneqq {\left\{ {\begin{array}{*{20}{c}}
			{{{\boldsymbol{d}}^{(n,0)}}}\\
			{{{\boldsymbol{d}}^{(n - 1,1)}}}\\
			\vdots \\
			{{{\boldsymbol{d}}^{(0,n)}}}
	\end{array}} \right\}}.
\end{equation}
Here ${{\boldsymbol{q}}_n}(\zeta^3)$ denotes a column array of $n$th order directors, and ${\boldsymbol{q}}(\zeta^3)$ is a column array of ${{\boldsymbol{q}}_n}(\zeta^3 )\,(n\in\left\{1,...,N\right\})$, which is called a \textit{generalized director vector}. ${{\boldsymbol{1}}_{m}}$ denotes the identity matrix of dimension $m$. Note that ${\boldsymbol{\pi}}_n$ solely depends on the transverse coordinates, which is independent from the deformations. Thus, taking the first variation of Eq.\,(\ref{compact_xt_pi}) yields
\begin{equation}
	\label{del_x_t_general}
	\delta{\boldsymbol{x}}={\boldsymbol{\Pi }}({\zeta^1},{\zeta ^2})^\mathrm{T}\delta\boldsymbol{q}(\zeta^3).
\end{equation}
In the current configuration, the covariant base vectors are defined by $\boldsymbol{g}_i\coloneqq\partial{\boldsymbol{x}}/\partial{\zeta^i}\,(i=1,2,3)$, so that we have
\begin{equation}
	\left\{\setlength{\arraycolsep}{8pt}
	\renewcommand{\arraystretch}{1.5}\begin{array}{c}
		\begin{aligned}
			{{\boldsymbol{g}}_1}(\zeta^1,\zeta^2,\zeta^3) &= {{\boldsymbol{\Pi }}^\mathrm{T}_{,{\zeta ^1}}}({\zeta ^1},{\zeta ^2})\,{\boldsymbol{q}}(\zeta^3),\\
			{{\boldsymbol{g}}_2}(\zeta^1,\zeta^2,\zeta^3) &= {{\boldsymbol{\Pi }}^\mathrm{T}_{,{\zeta ^2}}}({\zeta ^1},{\zeta ^2})\,{\boldsymbol{q}}(\zeta^3),\\
			{{\boldsymbol{g}}_3}(\zeta^1,\zeta^2,\zeta^3) &= {\boldsymbol{\Pi }}^\mathrm{T}({\zeta ^1},{\zeta ^2})\,{{\boldsymbol{q}}_{,\zeta^3 }}(\zeta^3),
		\end{aligned}
	\end{array} \right.
\end{equation}
where $(\bullet)_{,{\zeta^i}}$ denotes the partial derivative with respect to $\zeta^i$. The deformation gradient can be expressed by \citep[p.\,478]{wriggers2006}
\begin{equation}
	{\boldsymbol{F}} = \boldsymbol{g}_i\otimes\boldsymbol{G}^i,
\end{equation}
such that the Green-Lagrange strain tensor is obtained by
\begin{equation}
	\label{gl_strn_def_init_cov}
	{\boldsymbol{E}} \coloneqq \frac{1}{2}({{\boldsymbol{F}}^\mathrm{T}}{\boldsymbol{F}} - {\boldsymbol{1}}) = \frac{1}{2}\left( {{g_{ij}} - {G_{ij}}} \right){{\boldsymbol{G}}^i} \otimes {{\boldsymbol{G}}^j},
\end{equation}
where $\boldsymbol{1}$ represents the identity tensor in $\Bbb{R}^3$, and $g_{ij}\coloneqq\boldsymbol{g}_i\cdot\boldsymbol{g}_j$, and $G_{ij}\coloneqq\boldsymbol{G}_i\cdot\boldsymbol{G}_j$.
\begin{remark} \label{remark_cent_axis_deform}\textit{Geometrical centroid of the current cross-section.} The position of the geometrical centroid in the current cross-section is defined by
\begin{subequations}
\begin{align}
	{\boldsymbol{C}} &\coloneqq \frac{1}{{A_t}}\int_{{\mathcal{A}_t}} {{{\boldsymbol{x}}}\,{\rm{d}}{\mathcal{A}_t}} \nonumber\\
	&= \frac{1}{{A}_t}\int_\mathcal{A} {{{\boldsymbol{x}}}\,\left\| {{{\boldsymbol{g}}_1} \times {{\boldsymbol{g}}_2}} \right\|\,{\rm{d}}{\zeta ^1}d{\zeta ^2}},
\end{align}
with the current cross-sectional area
\begin{equation}
	{A_t} \coloneqq \int_{{\mathcal{A}_t}} {{\rm{d}}{\mathcal{A}_t}} = \int_\mathcal{A} \left\| {{{\boldsymbol{g}}_1} \times {{\boldsymbol{g}}_2}} \right\|{{\rm{d}}{\zeta ^1}d{\zeta ^2}}.
\end{equation}
\end{subequations}
It can be easily verified that, for $N=1$ in Eq.\,(\ref{exact_taylor_Nth}), we have $\boldsymbol{C}=\boldsymbol{\varphi}$. That is, for $N=1$, the current axis of the beam always coincides with the geometrical centroid. For $N\gt1$, however, this is not always the case.
\end{remark}
\subsection{Parameterization of a lateral boundary surface using NURBS}
\label{nth_param_lateral_bd_surf}
\begin{figure}	
	\centering
	\begin{subfigure}[b] {0.45\textwidth} \centering
		\includegraphics[width=\linewidth]{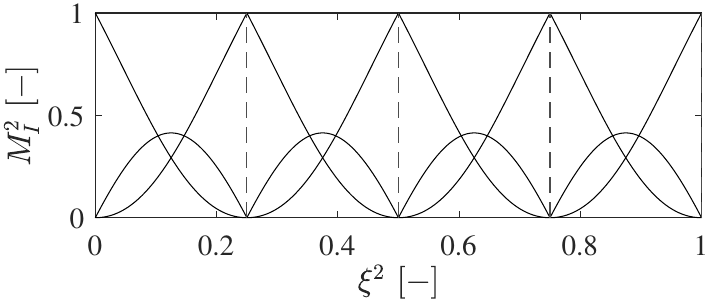}
		\caption{Original quadratic basis functions}
		\label{cs_bd_curve_basis_cl}
	\end{subfigure}
	\begin{subfigure}[b] {0.45\textwidth} \centering
		\includegraphics[width=\linewidth]{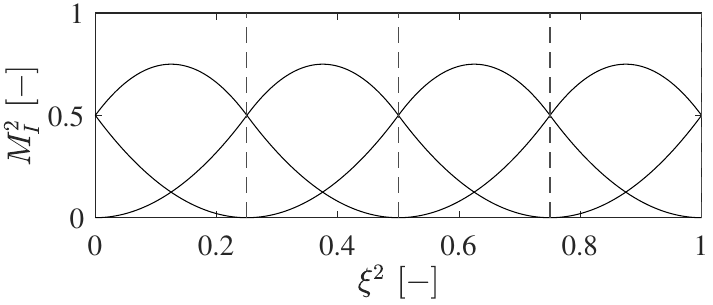}
		\caption{Periodic basis functions after knot removal}
		\label{cs_bd_curve_basis_uc}
	\end{subfigure}
	\caption{Comparison of non-periodic and periodic NURBS basis functions for the cross-section's boundary curve. (a) Non-periodic NURBS basis functions with knot vector $\tilde{\varXi}_\mathrm{cl}^2=\left\{0,0,0,1/4,1/4,1/2,1/2,3/4,3/4,1,1,1\right\}$, (b) Periodic NURBS basis functions with knot vector $\tilde{\varXi}_\mathrm{ucl}^2=\left\{0,1/8,1/4,3/8,1/2,5/8,3/4,7/8,1\right\}$. Note that the periodic basis functions are all translations of each other \citep{rogers2001introduction}. Vertical dashed lines divide the nonzero knot spans.}
	\label{cs_bd_curve_basis_plot}
\end{figure}

\begin{figure*}
	\centering
	\begin{subfigure}[b] {0.323\textwidth} \centering
		\includegraphics[width=\linewidth]{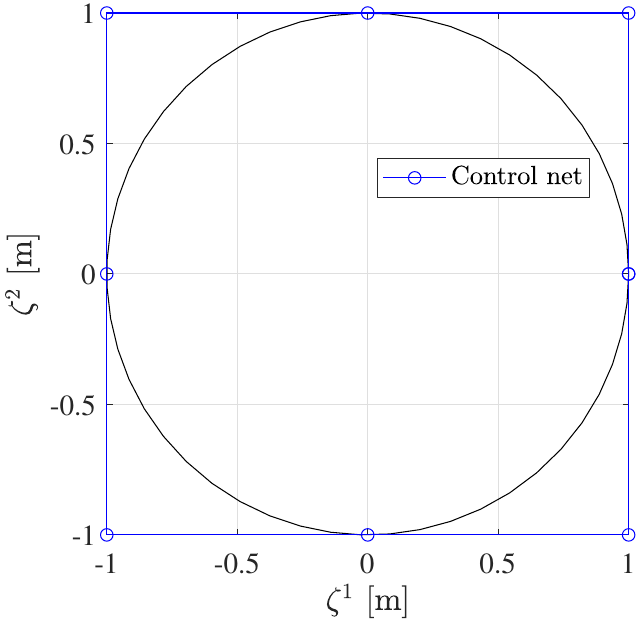}
		\caption{Exact, non-periodic, $m_\mathrm{cp}=8$}
		\label{cs_bd_curve_geom_cl}
	\end{subfigure}
	\begin{subfigure}[b] {0.33\textwidth} \centering
		\includegraphics[width=\linewidth]{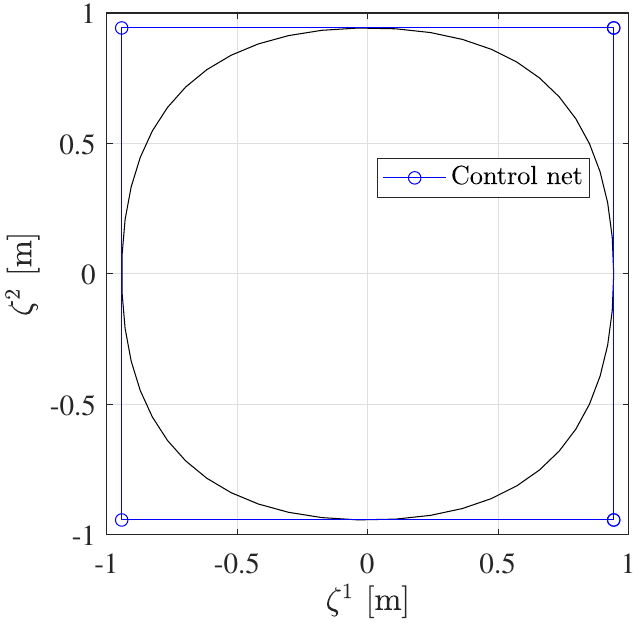}
		\caption{Approximate, periodic, $m_\mathrm{cp}=6$}
		\label{cs_bd_curve_geom_uc}
	\end{subfigure}
	\begin{subfigure}[b] {0.33\textwidth} \centering
		\includegraphics[width=\linewidth]{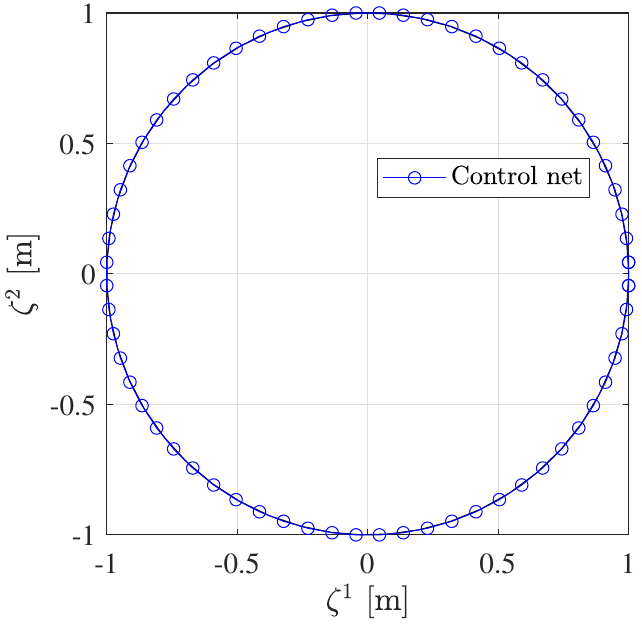}
		\caption{Approximate, periodic, $m_\mathrm{cp}=66$}
		\label{cs_bd_curve_geom_uc_ref}
	\end{subfigure}
	\caption{Representation of the boundary curve of a circular cross-section with radius $R=1\,\mathrm{m}$ by non-periodic or periodic NURBS. (a) Exact circle represented by non-periodic NURBS with $q=2$, $m_\mathrm{cp}=8$, whose corresponding basis functions are given in Fig.\,\ref{cs_bd_curve_basis_cl}. (b) Approximate circle represented by periodic NURBS with $q=2$, $m_\mathrm{cp}=6$, whose corresponding basis functions are shown in Fig.\,\ref{cs_bd_curve_basis_uc}. (c) Approximate circle represented by periodic NURBS with $q=2$ and $m_\mathrm{cp}=66$.}
	\label{cs_bd_curve_basis_plot}	
\end{figure*}

The initial geometry of the lateral boundary surface of the beam is decomposed into two parts: First, the \textit{directed} axis curve, described by a NURBS curve with two attached orthonormal directors, $\boldsymbol{D}_1$ and $\boldsymbol{D}_2$. Second, the cross-section defined by the NURBS curve in the plane spanned by the two initial directors
\begin{equation}
	\label{plane_zeta12_curve_nurbs}
	{\zeta ^\gamma }(\boldsymbol{\xi}) = \sum\limits_{I = 1}^{{{m}_{{{\mathrm{cp}}}}}} {{C^{\gamma}_{I}}({\xi ^1})\,{M^q_I}({\xi ^2})},\,\gamma\in\left\{1,2\right\},
\end{equation}
where ${\big({C^1_{I},C^2_{I}}\big)}\!\in\Bbb{R}\times\Bbb{R}$ represents the position of control point in the plane, and ${{m}}_{\mathrm{cp}}$ denotes the total number of control points along the curve, and $M^q_I(\xi^2)$ denotes the $I$-th NURBS basis function of order $q$. The definition and basic properties of NURBS basis function can be found in \cite{piegl1996nurbs}. $\xi^1$ and $\xi^2$ denote the convective coordinates along the axis and the cross-section's boundary curve, respectively. In this paper, we define those two convective coordinates by parametric coordinates of NURBS. Here and hereafter, we often use $\boldsymbol{\xi}\coloneqq\left[{\xi^1},{\xi^2}\right]^\mathrm{T}$ for brevity. It is noted that ${C}^\gamma_{I}$ depends on the coordinate $\xi^1$ only if a dimension of the initial cross-section is varying along the center axis. However, in this paper, we restrict our discussion to non-varying cross-sections along the axis, that is, we consider only those cases where ${C}^\gamma_{I}$ does not depend on the coordinate $\xi^1$. In the subsequent formulation of the closest point projection based on a local Newton-Raphson iteration, we need at least a $C^2$-continuous surface parameterization. However, if we use a \textit{clamped} knot vector for the \textit{closed} cross-section boundary curves, the displacement continuity is typically $C^0$ at the interface between the two end points where the end control points are matched, which may lead to difficult convergence of the local Newton-Raphson iteration. Further, contact tractions are singular at $C^0$ surface points, which can lead to spurious contact deformations. Thus, so as to have at least $C^2$-continuity of the whole lateral boundary surface, we employ the following approaches:
\begin{enumerate}
	\item In the representation of the cross-section's boundary curve, an unclamped knot vector is utilized, which allows higher order continuity at the interface between two end points of the closed boundary curves.
	\item We also reduce the multiplicity of internal knots to one, so that we have $C^{q-1}$ continuity in the entire curve of the cross-section's boundary, where $q$ denotes the order of basis functions.
	\item In the spatial discretization of kinematic variables as well as the initial geometry of the axis, we use NURBS basis functions, in the framework of \textit{isogeometric} analysis, in order to have higher order continuity in the axial direction of the beam as well.
\end{enumerate}
The \textit{removal} of knots may change the initial geometry of the cross-section. However, increasing the DOFs in the geometry, e.g., the number of control points, and the order of basis functions, before the knot removal enables to reduce the \textit{loss of geometry}. For example, we consider a circular cross-section of radius $R=1\,\mathrm{m}$. NURBS is capable of exactly representing conic sections like circle. Fig.\,\ref{cs_bd_curve_basis_cl} shows quadratic NURBS basis functions utilized to represent the circle in Fig.\,\ref{cs_bd_curve_geom_cl}, which have $C^0$-continuity at the interface between nonzero knot spans. To remove those non-smooth points we apply the operations, \textit{unclamping} and \textit{knot removal} for the end and internally repeated knots, respectively, and then we finally obtain the basis functions, shown in Fig.\,\ref{cs_bd_curve_basis_uc}. Further details on these \textit{unclamping}, and \textit{knot removal} procedures can be found in \citet{piegl1996nurbs} and \citet{rogers2001introduction}. We obtain those periodic NURBS curves by using the function \textit{MakePeriodic} in the commercial program \textit{Rhinoceros} 3D\footnote{Version 7, Robert McNeel \& Associates.}. It is clearly shown in Fig.\,\ref{cs_bd_curve_geom_uc} that the initial geometry of the cross-section is not an exact circle anymore. However, as Fig.\,\ref{cs_bd_curve_geom_uc_ref} shows, by increasing the number of control points before those operations, the loss of geometry approaches zero. Further, it should be noted that the DOFs in the geometry of initial cross-section is associated with the parameterization only, and does not have any influence on the computational costs, since the deformation of the current cross-section is solely described by director vectors.

\begin{figure} \centering	
	\includegraphics[width=0.9875\linewidth]{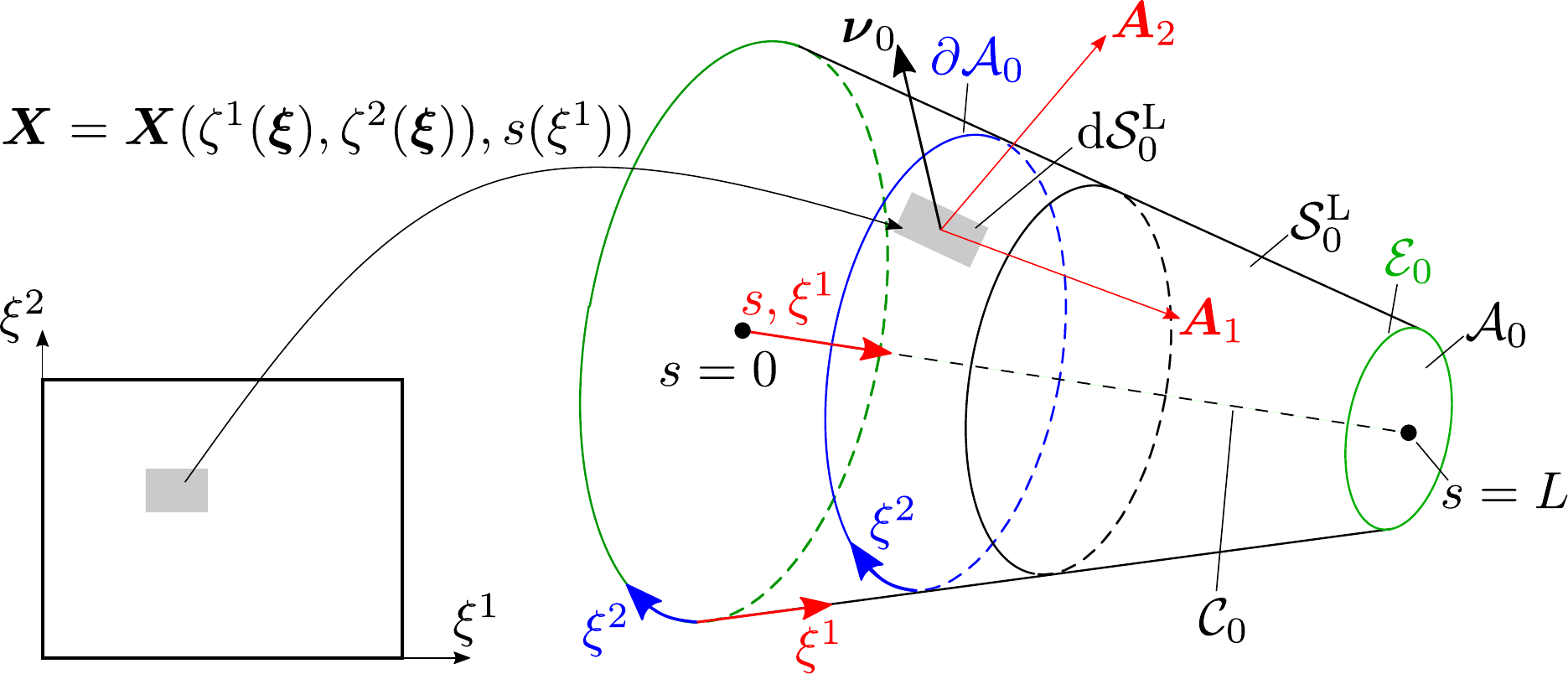}
	\caption{Parameterization of the beam's initial lateral boundary surface $\mathcal{S}^\mathrm{L}_0$. Note that the initial axis is parameterized by the same coordinate $\xi^1$ used for the longitudinal direction of the lateral boundary surface.}
	\label{beam_param_lat_bdc_sf}
\end{figure}
\subsubsection{Convective base vectors and metrics}
Using the parameterization of the coordinates $s=s(\xi^1)$ in Eq.\,(\ref{s_xi_relation}) and $\zeta^{\gamma}=\zeta^{\gamma}(\boldsymbol{\xi})$ of Eq.\,(\ref{plane_zeta12_curve_nurbs}), the initial position vector of the beam in Eq.\,(\ref{beam_init_config_kin}) can be reparameterized on the initial lateral boundary surface $\mathcal{S}_0^\mathrm{L}$, as
\begin{align}
	\label{init_pos_reparam_x}
	{\boldsymbol{X}} = {\boldsymbol{\varphi }}_0(s(\xi^1)) + {\zeta^\gamma }(\boldsymbol{\xi})\,{{\boldsymbol{D}}_\gamma }\left(s(\xi^1)\right).
\end{align}
Similarly, we can also reparameterize the current position vector of the beam in Eq.\,(\ref{compact_xt_pi}), on the current lateral boundary surface $\mathcal{S}_t^\mathrm{L}$, as
\begin{equation}
	\label{reparam_x_lat}
	{\boldsymbol{x}} = {\boldsymbol{\Pi }}({\zeta ^1}({\boldsymbol{\xi }}),{\zeta ^2}({\boldsymbol{\xi }}))^\mathrm{T}\,{\boldsymbol{q}}\left(s({\xi ^1})\right).
\end{equation}
For the given parameterization of Eq.\,(\ref{init_pos_reparam_x}), we define the covariant base vectors on the initial lateral boundary surface $\mathcal{S}_0^\mathrm{L}$, as 
\begin{equation}
	\label{app_1st_deriv_init_pt}
	\left. \begin{array}{c}
		\begin{aligned}
			\boldsymbol{A}_1\coloneqq{\boldsymbol{X}}_{,1}& = {\boldsymbol{\varphi}}_{0,1}+{\zeta _{,1}^\gamma}\,{\boldsymbol{D}}_{\gamma}+{\zeta}^{\gamma}{\boldsymbol{D}}_{\gamma,1},\\
			\boldsymbol{A}_2\coloneqq{\boldsymbol{X}}_{,2}& = {\zeta^{\gamma}_{,2}}\,{\boldsymbol{D}}_\gamma,
		\end{aligned}
	\end{array} \right\}
\end{equation}
and let $A_{\alpha\beta}\coloneqq{\boldsymbol{A}_\alpha}\cdot{\boldsymbol{A}_\beta}$ denote the covariant components of the surface metric, and here and hereafter, we define ${(\bullet)_{,\alpha }} \coloneqq \partial (\bullet)/\partial {\xi ^\alpha }$, and ${(\bullet)_{,\alpha\beta}} \coloneqq \partial^2 (\bullet)/\partial {\xi ^\alpha }\partial{\xi^\beta}$ $(\alpha,\beta\in\left\{1,2\right\})$. An initial convective surface basis $\left\{ {{{\boldsymbol{A}}_1},{{\boldsymbol{A}}_2},{\boldsymbol{\nu}}_0} \right\}$ on $\mathcal{S}_0^\mathrm{L}$ can be obtained by defining an outward unit normal vector on the surface as ${{\boldsymbol{\nu}}_0}\coloneqq{{\boldsymbol{A}}_1} \times {{\boldsymbol{A}}_2}/\left\| {{{\boldsymbol{A}}_1} \times {{\boldsymbol{A}}_2}} \right\|$, see Fig.\,\ref{beam_param_lat_bdc_sf} for an illustration.
\begin{remark} The first order derivatives of $\zeta^{\alpha}$ with respect to $\xi^\beta$, i.e., $\zeta^\alpha_{,\beta}$ $(\alpha,\beta\in\left\{1,2\right\})$ can be simply obtained, from Eq.\,(\ref{plane_zeta12_curve_nurbs}), as
\begin{subequations}
	\label{app_1st_deriv_pcoor0}
\begin{align}
	\zeta _{,1}^\alpha  &= 0, \label{app_1st_deriv_pcoor}\\
	\zeta _{,2}^\alpha  &= \sum\limits_{I = 1}^{{{m}_{{\rm{cp}}}}} {C^{\alpha}_{I}{M^q_{I,2}}({\xi ^2})},\label{app_1st_deriv_pcoor_1}
\end{align}
\end{subequations}
and the second order derivatives are
\begin{subequations}
\begin{align}
	\zeta _{,11}^\alpha  &= \zeta _{,12}^\alpha  = 0,\label{app_2nd_deriv_pcoor}\\
	\zeta _{,22}^\alpha  &= \sum\limits_{I = 1}^{{{m}_{{\rm{cp}}}}} {C^{\alpha}_{I}{M^q_{I,22}}({\xi ^2})},\label{app_2nd_deriv_pcoor_1}
\end{align}
\end{subequations}
where the dependence of $C^\alpha_I$ on $\xi^1$ vanishes due to the assumption of non-varying initial cross-sections along the axis. 
\end{remark}
\noindent Further we define the covariant base vectors on the current lateral surface $\mathcal{S}_t^\mathrm{L}$, as
\begin{equation}
	\label{app_1st_deriv_cur_pt}
	\left. \begin{array}{c}
		\begin{aligned}
			\boldsymbol{a}_1\coloneqq{\boldsymbol{x}}_{,1}& = {\zeta _{,1}^\alpha}\,{\boldsymbol{\Pi }}_{,{\zeta ^\alpha }}^{\rm{T}}\,{\boldsymbol{q}} + {{\boldsymbol{\Pi }}^{\rm{T}}}{{\boldsymbol{q}}_{,1}},\\
			\boldsymbol{a}_2\coloneqq{\boldsymbol{x}}_{,2}& = {\zeta _{,2}^\alpha}\,{\boldsymbol{\Pi }}_{,{\zeta ^\alpha }}^{\rm{T}}\,{\boldsymbol{q}},
		\end{aligned}
	\end{array} \right\}
\end{equation}
and let ${{a}_{\alpha \beta }} \coloneqq {{\boldsymbol{a}}_\alpha } \cdot {{\boldsymbol{a}}_\beta }$ denote the covariant components of the surface metric. The contravariant base vectors ${{\boldsymbol{a}}^\alpha}$ are defined by ${{\boldsymbol{a}}_\alpha }\cdot{{\boldsymbol{a}}^\beta}=\delta_{\alpha}^{\beta}$ $\left(\alpha,\beta\in\left\{1,2\right\}\right)$, where $\delta_{\alpha}^{\beta}$ denotes the Kronecker-delta symbol, which leads to ${{\boldsymbol{a}}^\alpha } = {a^{\alpha \beta }}{{\boldsymbol{a}}_\beta }$, where $a^{\alpha\beta}\coloneqq{{\boldsymbol{a}}^\alpha}\cdot{{\boldsymbol{a}}^\beta}$ denotes the contravariant components of the surface metric, calculated from
\begin{equation}
	\left[ {\begin{array}{*{20}{c}}
			{{a^{11}}}&{{a^{12}}}\\
			{{a^{21}}}&{{a^{22}}}
	\end{array}} \right] = {\left[ {\begin{array}{*{20}{c}}
				{{a_{11}}}&{{a_{12}}}\\
				{{a_{21}}}&{{a_{22}}}
		\end{array}} \right]^{ - 1}}.
\end{equation}
A current convective surface basis $\left\{ {{{\boldsymbol{a}}_1},{{\boldsymbol{a}}_2},{\boldsymbol{\nu}}_t} \right\}$ on $\mathcal{S}_t^\mathrm{L}$ can be obtained by defining an outward unit normal vector on the surface as ${{\boldsymbol{\nu}}_t}\coloneqq{{\boldsymbol{a}}_1} \times {{\boldsymbol{a}}_2}/\left\| {{{\boldsymbol{a}}_1} \times {{\boldsymbol{a}}_2}} \right\|$. The covariant components of curvature are defined by $b_{\alpha\beta}\coloneqq {{\boldsymbol{a}}_{\alpha ,\beta }} \cdot {{\boldsymbol{\nu }}_t}$, where $\boldsymbol{a}_{\alpha,\beta}=\boldsymbol{x}_{,\alpha\beta}$, and
\begin{equation}
	\label{app_2nd_deriv_cur_pt}
	\left. \begin{array}{c}
		\begin{aligned}
			{{\boldsymbol{x}}_{,11}} &= \left( {{\zeta _{,1}^\beta}\,\zeta _{,1}^\alpha\,{\boldsymbol{\Pi }}_{,{\zeta ^\alpha }{\zeta ^\beta }}^{\rm{T}} + \zeta _{,11}^\alpha\, {\boldsymbol{\Pi }}_{,{\zeta ^\alpha }}^{\rm{T}}} \right){\boldsymbol{q}} \\
			&+ 2{\zeta _{,1}^\alpha}\,{\boldsymbol{\Pi }}_{,{\zeta ^\alpha }}^{\rm{T}}\,{{\boldsymbol{q}}_{,1}} + {{\boldsymbol{\Pi }}^{\rm{T}}}{{\boldsymbol{q}}_{,11}},\\
			{{\boldsymbol{x}}_{,22}} &= \left( {\zeta _{,2}^\alpha\, \zeta _{,2}^\beta\, {\boldsymbol{\Pi }}_{,{\zeta ^\alpha }{\zeta ^\beta }}^{\rm{T}} + \zeta _{,22}^\alpha\,{\boldsymbol{\Pi }}_{,{\zeta ^\alpha }}^{\rm{T}}} \right){\boldsymbol{q}}, \\
			{{\boldsymbol{x}}_{,12}} &= \left( {\zeta _{,1}^\alpha \,\zeta _{,2}^\beta \,{\boldsymbol{\Pi }}_{,{\zeta ^\alpha }{\zeta ^\beta }}^{\rm{T}} + \zeta _{,12}^\alpha\, {\boldsymbol{\Pi }}_{,{\zeta ^\alpha }}^{\rm{T}}} \right){\boldsymbol{q}} \\
			&+ \zeta _{,2}^\alpha \,{\boldsymbol{\Pi }}_{,{\zeta ^\alpha }}^{\rm{T}}\,{{\boldsymbol{q}}_{,1}}.
		\end{aligned}
	\end{array} \right\}
\end{equation}
\subsubsection{Initial infinitesimal area element}
In the initial configuration, from the parametrization of the lateral boundary surface in Eq.\,(\ref{init_pos_reparam_x}), the infinitesimal area element of the initial lateral boundary surface can be expressed, using $\mathrm{d}s={\tilde j}\,\mathrm{d}\xi^1$ and the surface Jacobian ${\tilde J}\coloneqq\left\| {{{\boldsymbol{A}}_{1}}} \times {{\boldsymbol{A}}_{2}} \right\|$, by
\begin{align}
	\label{inf_area_init_config}	
	\mathrm{d}{{\mathcal{S}}^\mathrm{L}_0} = {\tilde J}\,\mathrm{d}{\xi^1}\mathrm{d}{\xi^2}=\frac{\tilde J}{{\tilde j}}\,\mathrm{d}\xi^2\,\mathrm{d}s.
\end{align} 
See Fig.\,\ref{beam_param_lat_bdc_sf} for an illustration.
\subsection{Variational formulation}
\subsubsection{Boundary value problem in strong form}
The (static) local linear momentum balance equations in combination with displacement and traction boundary conditions is stated in the initial configuration by the boundary value problem \citep{bonet2010nonlinear}
\begin{subequations}
	\label{bvp_str_form}
	\begin{alignat}{3}
		{\rm{Div}}\,\left({{\boldsymbol{F}}{\boldsymbol{S}}}\right) + {{\boldsymbol{b}}_0} &= \boldsymbol{0}\,\,\,\,&&\mathrm{in}\,\,&&\mathcal{B}_0,\\
		{{\boldsymbol{u}}} &= {{\boldsymbol{\bar u}}}\,\,\,&&\mathrm{on}\,\,&&\mathcal{S}^\mathrm{D}_0,\label{str_form_disp_bdc}\\
		{{\boldsymbol{F}}{\boldsymbol{S}}}{{\boldsymbol{\nu }}_0} &= {{\boldsymbol{\bar T}}_0}\,\,\,&&\mathrm{on}\,\,&&\mathcal{S}^\mathrm{N}_0,
	\end{alignat}
\end{subequations}
where ${\mathcal{S}^\mathrm{D}_0},\,{\mathcal{S}^\mathrm{N}_0}\!\subset\!\mathcal{S}_0$ denote the boundary surfaces where the displacement $\boldsymbol{u}\coloneqq{\boldsymbol{x}}-{\boldsymbol{X}}$ and traction are prescribed, respectively, and $\boldsymbol{S}$ denotes the second Piola-Kirchhoff stress tensor. Further, $\boldsymbol{\nu}_0$ denotes the unit outward normal vector on the surface of the undeformed configuration, $\boldsymbol{b}_0$ denotes the body force per unit undeformed volume, $\mathrm{Div}(\bullet)$ represents the divergence operator with respect to the initial configuration, and ${{\bar {\boldsymbol{u}}}}\in{{\Bbb{R}}^3}$ denotes the prescribed displacement vector. Here we consider hyperelastic materials, where the \textit{strain energy density}, defined by the strain energy per unit undeformed volume, is given in terms of the Green-Lagrange strain tensor $\boldsymbol{E}$, as $\Psi  = \Psi ({\boldsymbol{E}})$. Then, the constitutive equation is expressed by
\begin{equation}\label{2nd_pk_strs_comp}
	{\boldsymbol{S}} = {S^{ij}}{{\boldsymbol{G}}_i} \otimes {{\boldsymbol{G}}_j}\,\,\,\text{with}\,\,\,{S^{ij}} = \frac{{\partial \Psi }}{{\partial {E_{ij}}}}.
\end{equation}
\subsubsection{Weak form}
Exploiting the symmetries, ${\boldsymbol{E}}$ and ${\boldsymbol{S}}$ can be expressed in array form (Voigt notation), as ${\boldsymbol{\munderbar S}} \coloneqq {\left[ {{S^{11}},{S^{22}},{S^{33}},{S^{12}},{S^{13}},{S^{23}}} \right]^{\mathrm{T}}}$ and ${\boldsymbol{\munderbar E}} \coloneqq {\left[ {{E_{11}},{E_{22}},{E_{33}},2{E_{12}},2{E_{13}},2{E_{23}}} \right]^{\mathrm{T}}}$. From Eq.\,(\ref{gl_strn_def_init_cov}), the first variation of the covariant components of the Green-Lagrange strain tensor can then be expressed, in Voigt notation, by
\begin{align}
	\label{del_E_voigt_def_psi}
	\delta {\munderbar {\boldsymbol{E}}} = \left\{ {\begin{array}{*{20}{c}}
			{\delta {E_{11}}}\\
			{\delta {E_{22}}}\\
			{\delta {E_{33}}}\\
			{2\delta {E_{12}}}\\
			{2\delta {E_{13}}}\\
			{2\delta {E_{23}}}
	\end{array}} \right\} = {{\boldsymbol{\varXi }}}\delta {\boldsymbol{q}},
\end{align}
where
\begin{equation}
	{{\boldsymbol{\varXi }}} \coloneqq {\left[\setlength{\arraycolsep}{5pt}
		\renewcommand{\arraystretch}{1.25}{\begin{array}{*{20}{c}}
				{{{\boldsymbol{q}}^\mathrm{T}}{{\boldsymbol{\Pi }}_{,{\zeta ^1}}}\,{{\boldsymbol{\Pi }}_{,{\zeta ^1}}^\mathrm{T}}}\\
				{{{\boldsymbol{q}}^\mathrm{T}}{{\boldsymbol{\Pi }}_{,{\zeta ^2}}}\,{{\boldsymbol{\Pi }}_{,{\zeta ^2}}^\mathrm{T}}}\\
				{{{\boldsymbol{q}}_{,\zeta^3 }^\mathrm{T}}{{\bf{\Pi }}}\,{\bf{\Pi }}^\mathrm{T}{{(\bullet)}_{,\zeta^3}}}\\
				{{\boldsymbol{q}}^\mathrm{T}}\left({{\boldsymbol{\Pi }}_{,{\zeta ^2}}}\,{{\boldsymbol{\Pi }}_{,{\zeta ^1}}^\mathrm{T}}+{{\boldsymbol{\Pi }}_{,{\zeta ^1}}}\,{{\boldsymbol{\Pi }}_{,{\zeta ^2}}^\mathrm{T}}\right)\\
				{{{\boldsymbol{q}}_{,\zeta^3 }^\mathrm{T}}\,{{\boldsymbol{\Pi }}}\,{{\boldsymbol{\Pi }}_{,{\zeta^1}}^\mathrm{T}} + {{\boldsymbol{q}}^\mathrm{T}}{{\boldsymbol{\Pi }}_{,{\zeta^1}}}\,{\boldsymbol{\Pi }}^\mathrm{T}{{(\bullet)}_{,\zeta^3 }}}\\
				{{{\boldsymbol{q}}_{,\zeta^3 }^\mathrm{T}}\,{{\boldsymbol{\Pi }}}\,{{\boldsymbol{\Pi }}^\mathrm{T}_{,{\zeta^2}}} + {{\boldsymbol{q}}^\mathrm{T}}{{\boldsymbol{\Pi }}_{,{\zeta^2}}}\,{\boldsymbol{\Pi }}^\mathrm{T}{{(\bullet)}_{,\zeta^3 }}}
		\end{array}} \right]}.
\end{equation}
Then, the internal virtual work can be written as
\begin{equation}
	\label{intv_work}
	{G_{{\mathop{\rm int}} }}({\boldsymbol{q}},\delta {\boldsymbol{q}}) = \int_0^L {{{\boldsymbol{R}}^\mathrm{T}}\delta {\boldsymbol{q}}\,\mathrm{d}s},
\end{equation}
where the \textit{resultant of the stress} is obtained by
\begin{equation}
	\label{cs_integ_strs_res}
	\boldsymbol{R}\coloneqq{\int_\mathcal{A} {{{\boldsymbol{\varXi }}^\mathrm{T}}{\munderbar {\boldsymbol{S}} }\,{j_0}\,\mathrm{d}\mathcal{A}}},
\end{equation}
with $\mathrm{d}\mathcal{A}\coloneqq\mathrm{d}\zeta^1\mathrm{d}\zeta^2$. Further the external virtual work due to the body force, and surface tranctions on the lateral surface and the cross-section at the ends of rod, i.e., ${\left. {{\mathcal{A}_0}} \right\rvert_{s \in \Gamma_\mathrm{N} }}$ is given by
\begin{equation}
	\label{ext_vwork_q}
	{G_{{\rm{ext}}}}(\delta {\boldsymbol{q}}) = \int_0^L {\delta {{\boldsymbol{q}}^{\rm{T}}}{\boldsymbol{\bar R}}\,\mathrm{d}s}+{\left[ {\delta {{\boldsymbol{q}}^{\rm{T}}}{{{\boldsymbol{\bar R}}}_0}} \right]_{s \in {\Gamma _{\rm{N}}}}},
\end{equation}
where the resultant of the body force and traction on the lateral surface is expressed by
\begin{align}
	\label{ext_vwork_def_r_bar}
	{\boldsymbol{\bar R}} &\coloneqq \int_{\mathcal{A}} {{\boldsymbol{\Pi }}\,{{{\boldsymbol{b}}}_0}\,{j_0}\,\mathrm{d}\mathcal{A}}\nonumber\\
	&+\frac{1}{{\tilde j}}\int_{{\varXi^2}} {{\boldsymbol{\Pi }}\,{{{\boldsymbol{\bar T}}}_0}\left\| {{{\boldsymbol{A}}_1} \times {{\boldsymbol{A}}_2}} \right\|\mathrm{d}{\xi ^2}}.
\end{align}
Further, the resultant of the traction on the cross-section at the ends is expressed, using $\mathrm{d}\mathcal{A}_0=\left\| {{{\boldsymbol{G}}_1} \times {{\boldsymbol{G}}_2}} \right\|\mathrm{d}\mathcal{A}$, as
\begin{align}
	\label{res_ext_cs_trac_R0}
	{{\boldsymbol{\bar R}}_0} \coloneqq \int_{\mathcal{A}} {{\boldsymbol{\Pi }}\,{{{\boldsymbol{\bar T}}}_0}\left\| {{{\boldsymbol{G}}_1} \times {{\boldsymbol{G}}_2}} \right\|\mathrm{d}\mathcal{A}}.
\end{align}
Finally, the variational equation can be stated as: Find the generalized director vector $\boldsymbol{q}\in{{\mathcal{V}}}$ such that
\begin{align}
	\label{var_eq_q}
	&{G_{{\mathop{\rm int}} }}({\boldsymbol{q}},\delta {\boldsymbol{q}}) + {G_{{\mathop{\rm N}} }}({\boldsymbol{q}},\delta {\boldsymbol{q}}) \nonumber\\
	&= {G_{{\rm{ext}}}}(\delta {\boldsymbol{q}}) + {G_{{\rm{nc}}}}({\boldsymbol{q}},\delta {\boldsymbol{q}}),\,\,\forall \delta {\boldsymbol{q}} \in {\mathcal{V}_0},
\end{align}
where ${G_{{\mathop{\rm N}} }}({\boldsymbol{q}},\delta {\boldsymbol{q}})$ represents the internal virtual work due to normal contact whose detailed expression is given in Section \ref{sec_frictionless_contact}, and ${G_{{\mathop{\rm nc}} }}({\boldsymbol{q}},\delta {\boldsymbol{q}})$ denotes the external virtual work due to non-conservative loads, e.g., the distributed follower load in Remark \ref{remark_imp_dfol}. We also define
\begin{align}
	{\mathcal{V}} \coloneqq \left\{ {\left. {{\boldsymbol{q}} \in {{\left[ {{H^1}(0,L)} \right]}^{{n_{{\rm{cs}}}}}}} \right\rvert{\boldsymbol{q}} = {{{\boldsymbol{\bar q}}}_0}\,\,{\rm{on}}\;{\Gamma _{\rm{D}}}} \right\},
\end{align}
and
\begin{align}
	{\mathcal{V}_0} \coloneqq \left\{ {\left. {\delta{\boldsymbol{q}} \in {{\left[ {{H^1}(0,L)} \right]}^{{n_{{\rm{cs}}}}}}} \right\rvert\delta{\boldsymbol{q}} = {{\bf{0}}}\,\,{\rm{on}}\;{\Gamma _{\rm{D}}}} \right\}.
\end{align}
The director vectors are prescribed at the boundary ${\Gamma_\mathrm{D}}\ni{s}$. It is noted that ${\Gamma_\mathrm{D}} \cap {\Gamma_\mathrm{N}} = \emptyset$, and ${\Gamma_\mathrm{D}} \cup {\Gamma_\mathrm{N}} = \left\{ {0,L} \right\}$.
\subsubsection{Linearization}
For hyperelastic materials, in general, the constitutive relation between $\boldsymbol{S}$ and $\boldsymbol{E}$ is nonlinear. Taking the directional derivative of $\boldsymbol{S}$ gives
\begin{equation}
	\label{inc_2pk}
	D\boldsymbol{S}\cdot\Delta{\boldsymbol{x}}=\boldsymbol{\mathcal{C}}:D\boldsymbol{E}\cdot\Delta{\boldsymbol{x}},
\end{equation}
where $D(\bullet)\cdot(*)$ denotes the directional derivative of $(\bullet)$ in the direction of $(*)$, and $\Delta\boldsymbol{x}$ represents the increment of the current position of a material point. The material (Lagrangian) elasticity tensor $\boldsymbol{\mathcal{C}}$ is expressed by
\begin{equation}
	\label{elasticity_c_def}
	\boldsymbol{\mathcal{C}}= \frac{{\partial {\boldsymbol{S}}}}{{\partial {\boldsymbol{E}}}} = {C^{ijk\ell}}{{\boldsymbol{G}}_i} \otimes {{\boldsymbol{G}}_j} \otimes {{\boldsymbol{G}}_k} \otimes {{\boldsymbol{G}}_\ell},
\end{equation}
with
\begin{equation}
	{C^{ijk\ell}} = \frac{{{\partial ^2}\Psi }}{{\partial {E_{ij}}\,\partial {E_{k\ell }}}}.
\end{equation}
Eq.\,(\ref{inc_2pk}) can be rewritten, using Eq.\,(\ref{del_E_voigt_def_psi}), as
\begin{equation}
	\label{inc_2pk_voigt}
	D\munderbar{\boldsymbol{S}}\cdot{\Delta\boldsymbol{x}}=\munderbar{\munderbar{\boldsymbol{\mathcal{C}}}}\,\boldsymbol{\varXi}\,\Delta{\boldsymbol{q}},
\end{equation}
where
\begin{equation}
	{\boldsymbol{\munderbar{\munderbar{\mathcal{C}}}}} \coloneqq \left[ {\setlength{\arraycolsep}{1.225pt}
		\renewcommand{\arraystretch}{1}\begin{array}{*{20}{c}}{{{\mathcal{C}}^{1111}}}&{{{\mathcal{C}}^{1122}}}&{{{\mathcal{C}}^{1133}}}&{{{\mathcal{C}}^{1112}}}&{{{\mathcal{C}}^{1113}}}&{{{\mathcal{C}}^{1123}}}\\
			{}&{{{\mathcal{C}}^{2222}}}&{{{\mathcal{C}}^{2233}}}&{{{\mathcal{C}}^{2212}}}&{{{\mathcal{C}}^{2213}}}&{{{\mathcal{C}}^{2223}}}\\
			{}&{}&{{{\mathcal{C}}^{3333}}}&{{{\mathcal{C}}^{3312}}}&{{{\mathcal{C}}^{3313}}}&{{{\mathcal{C}}^{3323}}}\\
			{}&{}&{}&{{{\mathcal{C}}^{1212}}}&{{{\mathcal{C}}^{1213}}}&{{{\mathcal{C}}^{1223}}}\\
			{}&{{\rm{sym}}{\rm{.}}}&{}&{}&{{{\mathcal{C}}^{1313}}}&{{{\mathcal{C}}^{1323}}}\\
			{}&{}&{}&{}&{}&{{{\mathcal{C}}^{2323}}}
	\end{array}} \right].
\end{equation}
The internal virtual work of Eq.\,(\ref{intv_work}) is nonlinear in terms of the generalized director $\boldsymbol{q}$. Thus, in order to solve the variational equation using the Newton-Raphson iteration, we need to linearize Eq.\,(\ref{intv_work}). The directional derivative of $\boldsymbol{R}$ can be obtained by using Eq.\,(\ref{inc_2pk_voigt}), as
\begin{equation}
	D\boldsymbol{R}\cdot\Delta\boldsymbol{q}=\Bbb{C}\,\Delta\boldsymbol{q},\,\,\Bbb{C}\coloneqq\int_\mathcal{A} {{{\boldsymbol{\varXi }}^{\rm{T}}}\,{\boldsymbol{\munderbar{\munderbar{\mathcal{C}}}}}\,\boldsymbol{\varXi}\,{j_0}\,\mathrm{d}\mathcal{A}}.
\end{equation}
Thus, the directional derivative of the internal virtual work of Eq.\,(\ref{intv_work}) with the first variation of the Green-Lagrange strain tensor held constant, i.e., the material part of the tangent stiffness, is simply obtained by
\begin{equation}
	\label{tan_stiff_mat}
	{D_{\rm{M}}}{G_{{\mathop{\rm int}} }} \cdot \Delta {\boldsymbol{q}}= \int_0^L {\delta {{\boldsymbol{q}}^{\rm{T}}}\,{\Bbb{C}}\,\Delta {\boldsymbol{q}}}\,\mathrm{d}s.
\end{equation}
Further, the geometric part of the tangent stiffness is derived by taking the directional derivative of Eq.\,(\ref{intv_work}) with the second Piola-Kirchhoff stress part held constant, as
\begin{subequations}
	\label{tan_stiff_geom}
\begin{equation}
	\label{tan_stiff_geom_a}
	{D_{\rm{G}}}{G_{{\mathop{\rm int}} }} \cdot \Delta {\boldsymbol{q}} = \int_0^L {\delta {{\boldsymbol{q}}^{\rm{T}}}{{\boldsymbol{Y}}^{\rm{T}}}{{\boldsymbol{{k}}}_{\rm{G}}}{\boldsymbol{Y}}\Delta\boldsymbol{q}\,{\rm{d}}s},
\end{equation}
with
\begin{equation}
	\label{tan_stiff_geom_b}
	{{\boldsymbol{k}}_\mathrm{G}} \coloneqq \int_\mathcal{A} {{{{\boldsymbol{\bar k}}}_\mathrm{G}}}\,{j_0}\,\mathrm{d}\mathcal{A},
\end{equation}
\end{subequations}
where
\begin{subequations}
\begin{equation}
	{{\boldsymbol{\bar {{k}}}}_\mathrm{G}} = \left[ {\setlength{\arraycolsep}{3pt}
		\renewcommand{\arraystretch}{1.5}\begin{array}{*{20}{c}}
			\big({\bar {\boldsymbol{k}}_{\rm{G}}}\big)_{11} &\big(\bar {\boldsymbol{k}}_{\rm{G}}\big)_{12}\\
			{\mathrm{sym.}}&\big(\bar {\boldsymbol{k}}_{\rm{G}}\big)_{22}
	\end{array}} \right],	
\end{equation}
with
\begin{equation}
	\left. \renewcommand{\arraystretch}{1.5}\begin{array}{l}
		\big({\bar {\boldsymbol{k}}_{\rm{G}}}\big)_{11} \coloneqq {S^{\alpha\beta}}\,{{\boldsymbol{\Pi }}_{,{\zeta ^\alpha}}}\,{{\boldsymbol{\Pi }}_{,{\zeta ^\beta}}^\mathrm{T}},\\
		\big(\bar {\boldsymbol{k}}_{\rm{G}}\big)_{12} \coloneqq {{S^{13}}\,{{\boldsymbol{\Pi }}_{,{\zeta ^1}}}\,{{\boldsymbol{\Pi }}^\mathrm{T}} + {S^{23}}\,{{\boldsymbol{\Pi }}_{,{\zeta ^2}}}\,{{\boldsymbol{\Pi }}^\mathrm{T}}},\\
		\big(\bar {\boldsymbol{k}}_{\rm{G}}\big)_{22} \coloneqq {{S^{33}}\,{\boldsymbol{\Pi }}\,{{\boldsymbol{\Pi }}^\mathrm{T}}},
	\end{array} \right\}
\end{equation}
\end{subequations}
and
\begin{equation}
	{\boldsymbol{Y}} \coloneqq \left[\setlength{\arraycolsep}{3pt}
	\renewcommand{\arraystretch}{1.25}{\begin{array}{*{20}{c}}
			{{{\boldsymbol{1}}_{n_{{\rm{cs}}}}}}\\
			{{{\boldsymbol{1}}_{n_{{\rm{cs}}}}}{{(\bullet)}_{,s}}}
	\end{array}} \right].
\end{equation}
Combining Eqs.\,(\ref{tan_stiff_mat}) and (\ref{tan_stiff_geom}), we finally obtain the following increment of the internal virtual work
\begin{align}
	 \mathrm{D}{G_{{\mathop{\rm int}} }} \cdot \Delta \boldsymbol{q} &= \int_0^L {\delta {{\boldsymbol{q}}^{\rm{T}}}\left( {{\Bbb{C}} + {{\boldsymbol{Y}}^\mathrm{T}}{{\boldsymbol{k}}_{\rm{G}}}{\boldsymbol{Y}}} \right)} \Delta {\boldsymbol{q}}\,\mathrm{d}s\nonumber\\
	 &\eqqcolon\Delta {G_{{\mathop{\rm int}} }}({\boldsymbol{q}};\delta {\boldsymbol{q}},\Delta {\boldsymbol{q}}) .
\end{align}
\begin{remark} \label{remark_imp_dfol} \textit{Application of a moment load by a distributed follower load.} An end moment can be applied to a beam by employing the following linear distribution of the first Piola-Kirchhoff stress over the height, i.e., $-h/2\le{\zeta^1}\le{h/2}$ of the rectangular cross-section with dimension $h\times w$ on ${{{\left. {{\mathcal{A}_0}} \right\rvert}_{s \in {\Gamma _{\rm{N}}}}}}$\citep{betsch1995assumed, choi2021isogeometric}  
	\begin{equation}\label{end_moment_first_pk}
		{\boldsymbol{P}} = {\bar p}\,{{\boldsymbol{\nu }}_t} \otimes {{\boldsymbol{\nu }}_0}\,\,\text{with}\,\,{\bar p}\coloneqq{-\frac{M}{I}{\zeta^1}},\,\,{I=\frac{wh^3}{12}},
	\end{equation}
	and the outward unit normal vectors on the cross-sections at the ends of beam's axis in the initial and current configurations are obtained by
	\begin{equation}
		{{\boldsymbol{\nu }}_0} = \mathrm{sign}({{\boldsymbol{\nu }}_0})\frac{{{{\boldsymbol{G}}_1} \times {{\boldsymbol{G}}_2}}}{{\left\| {{{\boldsymbol{G}}_1} \times {{\boldsymbol{G}}_2}} \right\|}}\,\,\mathrm{on}\,\left.{{\mathcal{A}}_0}\right\rvert_{s\in\left\{0,L\right\}},
	\end{equation}
	and
	\begin{equation}
	{{\boldsymbol{\nu }}_t} = \mathrm{sign}({{\boldsymbol{\nu }}_t})\frac{{{{\boldsymbol{g}}_1} \times {{\boldsymbol{g}}_2}}}{{\left\| {{{\boldsymbol{g}}_1} \times {{\boldsymbol{g}}_2}} \right\|}}\,\,\mathrm{on}\,\left.{{\mathcal{A}}_t}\right\rvert_{s\in\left\{0,L\right\}},
	\end{equation}
	respectively, where the values of the signum functions can be simply determined by
	\begin{equation}
		\label{pnt_cnst_cpres_fix}
		\mathrm{sign}({{\boldsymbol{\nu }}_0}) = \mathrm{sign}({{\boldsymbol{\nu }}_t}) = \left\{ {\setlength{\arraycolsep}{1pt}
			\renewcommand{\arraystretch}{1}\begin{array}{*{20}{c}}
				\begin{array}{l}
					-1\\
					+1
				\end{array}&\begin{array}{l}
					\mathrm{at}\,\,{s=0},\\
					\mathrm{at}\,\,{s=L}.
				\end{array}
		\end{array}} \right.	
	\end{equation}
	Then, the prescribed surface traction vector can be expressed by ${\bar{\boldsymbol{T}}}_0\!=\!\boldsymbol{P}\,\boldsymbol{\nu}_0\!=\!{\bar p}\,{\boldsymbol{\nu}}_t$, and the external virtual work due to the follower load is obtained from Eq.\,(\ref{res_ext_cs_trac_R0}), as
	\begin{subequations}
	\begin{equation}
		\label{ext_vwork_mnt_follow}
		G_{{\rm{nc}}}(\boldsymbol{q},\delta{\boldsymbol{q}}) = {\left[ {\delta {{\boldsymbol{q}}^{\rm{T}}}{{{\boldsymbol{\bar R}}}_0}} \right]_{s \in {\Gamma _\mathrm{N}}}},
	\end{equation}
	with
	\begin{equation}
		{{\boldsymbol{\bar R}}_0} =  - \frac{M}{I}\int_{\mathcal{A}} {{\zeta ^1}\,{\boldsymbol{\Pi }}\,{{\boldsymbol{\nu }}_t}\left\| {{{\boldsymbol{G}}_1} \times {{\boldsymbol{G}}_2}} \right\|\mathrm{d}{\mathcal{A}}}.
	\end{equation}
	\end{subequations}
	Further, the increment of Eq.\,(\ref{ext_vwork_mnt_follow}) is derived as
	\begin{equation}
		\Delta G_{{\mathrm{nc}}}({\boldsymbol{q}};\delta {\boldsymbol{q}},\Delta {\boldsymbol{q}}) = {\left[ {\delta {{\boldsymbol{q}}^{\rm{T}}}{{\boldsymbol{S}}_{{\rm{nc}}}}\,\Delta {\boldsymbol{q}}} \right]_{s \in {\Gamma _\mathrm{N}}}},
	\end{equation}
	where
	\begin{subequations}
	\begin{equation}
		{{\boldsymbol{S}}_{{\rm{nc}}}} = \frac{M}{I}\int_\mathcal{A} {{\zeta ^1}{{\bar{\boldsymbol{S}}}_{{\rm{nc}}}}\,\mathrm{d}\mathcal{A}},
	\end{equation}
	with
	\begin{align}
	{{\bar{\boldsymbol{S}}}_{{\rm{nc}}}} &\coloneqq \mathrm{sign}(\boldsymbol{\nu}_t)\frac{{\left\| {{{\boldsymbol{G}}_1} \times {{\boldsymbol{G}}_2}} \right\|}}{{\left\| {{{\boldsymbol{g}}_1} \times {{\boldsymbol{g}}_2}} \right\|}}\times\nonumber\\
	&{\boldsymbol{\Pi }}\left( {{\boldsymbol{1}} - {{\boldsymbol{\nu }}_t} \otimes {{\boldsymbol{\nu }}_t}} \right)\left( {\widehat{{{\boldsymbol{g}}}_2}}\,{\boldsymbol{\Pi }}_{,{\zeta ^1}}^{\rm{T}} - {\widehat{{{\boldsymbol{g}}}_1}}\,{\boldsymbol{\Pi }}_{,{\zeta ^2}}^{\rm{T}} \right).
	\end{align}	
	\end{subequations}
\end{remark}
\noindent Here, $\widehat{{{\boldsymbol{g}}}_\alpha}$ denotes the skew-symmetric tensor associated with the dual vector ${{\boldsymbol{g}}}_\alpha$.
\section{A frictionless beam-to-beam contact formulation}\label{sec_frictionless_contact}
\label{sec_frictionless_contact}
\subsection{Impenetrability condition}
We present a continuum formulation for frictionless lateral contact between two slender bodies\footnote{This formulation can be easily extended to problems with more bodies by applying it to each pair of bodies.} based on the beam formulation presented in Sections \ref{nth_beam_kin}-\ref{nth_param_lateral_bd_surf}. We employ a Gauss point-to-surface contact formulation. The interacting bodies are denoted by ${{\mathcal{B}}}_0^{(\alpha)}$ ($\alpha\in\left\{1,2\right\}$) in their initial configurations. All the geometrical or physical quantities for the two bodies are indicated by the super- or subscript $(\alpha)$. Let $\mathcal{R}_t\coloneqq{\mathcal{R}_t^{(1)}}={\mathcal{R}_t^{(2)}}$ be the current contact interface between the two contacting bodies, which are subsets of $\mathcal{S}_t^\mathrm{L}$. The contact interface is pulled back to the boundary surface in the initial configuration of each body ${\mathcal{R}_0^{(\alpha)}}$, and the contact conditions are evaluated on $\mathcal{R}_0\coloneqq\mathcal{R}_0^{(1)}\ne\mathcal{R}_0^{(2)}$ \citep{laursen1993continuum}. We designate the bodies ${\mathcal{B}}_0^{(1)}$ and ${\mathcal{B}}_0^{(2)}$ as \textit{slave} and \textit{master}, which implies that all current points ${\boldsymbol{x}}\in{\mathcal{R}_t^{(1)}}$ are supposed to not penetrate the surface ${\mathcal{R}_t^{(2)}}$. The boundary value problem of Eq.\,(\ref{bvp_str_form}) can be rewritten for each of the interacting bodies, as
\begin{subequations}
	\label{nl_bvp_each_bd}
	\begin{alignat}{3}
		{\rm{Div}}\,{{\boldsymbol{P}}^{(\alpha )}} + {{\boldsymbol{b}}_0^{(\alpha )}} &= \boldsymbol{0}\,\,\,\,&&\mathrm{in}\,\,&&\mathcal{B}_0^{(\alpha )},\\
		{{\boldsymbol{u}}^{(\alpha )}} &= {{\boldsymbol{\bar u}}}^{(\alpha)}\,\,\,&&\mathrm{on}\,\,&&{\mathcal{S}^\mathrm{D}_0}^{(\alpha )},\\
		{{\boldsymbol{F}}^{(\alpha )}}{{\boldsymbol{S}}^{(\alpha )}}{\boldsymbol{\nu }}_0^{(\alpha )} &= {{\boldsymbol{\bar T}}^{(\alpha)}_0}\,\,\,&&\mathrm{on}\,\,&&{\mathcal{S}^\mathrm{N}_0}^{(\alpha )},
	\end{alignat}
\end{subequations}
$\alpha\in\left\{1,2\right\}$. Let ${\boldsymbol{\bar x}}\in\mathcal{S}_t^{\mathrm{L}(2)}$ be the solution of the following minimal distance problem
\begin{equation}
	\label{contact_min_disp_p_ul}
	{\boldsymbol{\bar x}} \coloneqq \arg \mathop {\min }\limits_{\boldsymbol{x}^{(2)}\in\mathcal{S}^\mathrm{L(2)}_t} \left\|{{\boldsymbol{x}}^{(1)}}-{{\boldsymbol{x}}^{(2)}}\right\|
\end{equation}
for a given point ${\boldsymbol{x}^{(1)}\in{{\mathcal{S}}_t^{\mathrm{L}(1)}}}$, and ${\bar {\boldsymbol{\nu}}}_t$ denote the outward unit normal vector at ${\boldsymbol{\bar x}}\in\mathcal{S}_t^{\mathrm{L}(2)}$. The impenetrability condition can then be evaluated, on the initial lateral surface $\mathcal{S}_0^{\mathrm{L}{(1)}}\ni{\boldsymbol{X}}^{(1)}$, as \citep{simo1992augmented}
\begin{subequations}
	\begin{align}
		{g_{\mathrm{N}}} \coloneqq \left( {{{\boldsymbol{x}}^{(1)}} - {{\boldsymbol{\bar x}}}} \right) \cdot {{\boldsymbol{\bar \nu }}_t}&\ge0,\label{cond_impenetrate_normal}\\
		{p_{\rm{N}}} \coloneqq  {\boldsymbol{\bar \nu }}_t \cdot {\boldsymbol{F}}^{(1)}{\boldsymbol{S}}^{(1)}{{\boldsymbol{\nu }}_0^{(1)}} &\ge 0,\label{cond_compf}\\
		{p_\mathrm{N}}\,{g_\mathrm{N}} &= 0,\label{cond_consist}\\
		{{\boldsymbol{t}}_\mathrm{T}} \coloneqq {\boldsymbol{F}}^{(1)}{\boldsymbol{S}}^{(1)}{{\boldsymbol{\nu }}_0^{(1)}} + {p_{\mathrm{N}}}\,{\boldsymbol{\bar \nu }}_t &= {\boldsymbol{0}}, \label{cond_nofriction}
	\end{align}
\end{subequations}
where ${\boldsymbol{\nu}}^{(1)}_0$ denotes the outward unit normal vector on the initial lateral surface of the slave body. Eq.\,(\ref{cond_impenetrate_normal}) denotes the impenetrability condition, and Eq.\,(\ref{cond_compf}) represents that contact should be compressive, where $p_\mathrm{N}$ denotes the contact pressure. Further, Eq.\,(\ref{cond_consist}) means the contact pressure applies only if the impenetrability condition is active, and Eq.\,(\ref{cond_nofriction}) implies no friction is considered, where $\boldsymbol{t}_\mathrm{T}$ denotes the tangential component of surface traction. The constraints of Eqs.\,(\ref{cond_impenetrate_normal})-(\ref{cond_consist}) can be replaced by the penalty regularization of the contact pressure \citep{simo1992augmented}
\begin{subequations}
	\label{pnt_cpres_reg}
	\begin{equation}
		p_\mathrm{N} = \epsilon_\mathrm{N}\left\langle {g_\mathrm{N}} \right\rangle\,\,\mathrm{on}\,\, \mathcal{S}_0^\mathrm{L}
	\end{equation}
	with
	\begin{equation}\label{gn_bracket}
		{\left\langle {{g_{\rm{N}}}} \right\rangle} \coloneqq \left\{ {\begin{array}{*{20}{c}}
				\begin{array}{l}
					- {g_\mathrm{N}}\\
					0
				\end{array}&\begin{array}{l}
					{\rm{if}}\,\,\,{g_\mathrm{N}} \le 0,\\
					\mathrm{otherwise},
				\end{array}
		\end{array}} \right.
	\end{equation}
\end{subequations}
where the constant $\epsilon_\mathrm{N}>0$ is the penalty parameter, and ${g_\mathrm{N}}\rightarrow0$ and $p_\mathrm{N}$ converges as ${\epsilon_\mathrm{N}}\rightarrow\infty$. That is, in order to reduce the error in the contact constraints or to avoid unphysical penetration, it is required to use a sufficiently large penalty parameter. However, it should be noted that, as we discuss in Section \ref{num_ex_slide_2beam}, a larger penalty parameter typically requires more surface Gauss integration points for the contact integral, and a smaller load increment, which makes the computation less efficient. An adaptive adjustment of penalty parameter remains future work. One can adjust the penalty parameter locally in order to control the penetration globally, such that the penetration does not exceed the maximum allowed value, see, e.g., the iterative adjustment approach of \cite{durville2012contact}. 
%
\subsubsection{Active set iteration}
If the contact surface $\mathcal{R}_0$ is known, the impenetrability condition becomes an equality constraint, and Eqs.\,(\ref{cond_impenetrate_normal})-(\ref{cond_consist}) can be rewritten as
\begin{subequations}
	\begin{align}
		{g_{\rm{N}}} &= 0\,\,\,{\rm{if}}\,\,{\boldsymbol{X}}^{(1)} \in \mathcal{R}_0,\label{equality_impen_cnst}\\
		{g_{\rm{N}}} &> 0\,\,\,{\rm{if}}\,\,{\boldsymbol{X}}^{(1)} \notin \mathcal{R}_0,\\
		{p_{\rm{N}}} &\ge 0\,\,\,{\rm{if}}\,\,{\boldsymbol{X}}^{(1)} \in \mathcal{R}_0,\\
		{p_{\rm{N}}} &= 0\,\,\,{\rm{if}}\,\,{\boldsymbol{X}}^{(1)} \notin \mathcal{R}_0,
	\end{align}
\end{subequations}
and the penalty regularization in Eq.\,(\ref{pnt_cpres_reg}) can be rewritten as
\begin{equation}
	\label{pnt_cnst_cpres_fix}
		{p_\mathrm{N}} = \left\{ {\begin{array}{*{20}{c}}
		\begin{array}{l}
			-\epsilon_\mathrm{N}\,{g_\mathrm{N}}\\
			0
		\end{array}&\begin{array}{l}
			{\rm{if}}\,\,\,{\boldsymbol{X}}^{(1)} \in \mathcal{R}_0,\\
			{\rm{if}}\,\,\,{\boldsymbol{X}}^{(1)} \in \mathcal{S}^\mathrm{L}_0\setminus\mathcal{R}_0.
		\end{array}
		\end{array}} \right.	
\end{equation}
However, since the contact surface is also unknown in general, we employ an iterative procedure, called the \textit{active set method}, in order to find the current contact area. The underlying idea of this method is to divide the inequality constraints into two groups -- the \textit{active set} and the \textit{inactive set} of constraints \citep{luenberger2016linear} -- and use two nested loops -- an \textit{outer loop} to find the correct set of active constraints, and an \textit{inner loop} to solve the nonlinear boundary value problem of Eq.\,(\ref{nl_bvp_each_bd}) using the Newton-Raphson iteration with a \textit{fixed} set of active constraints. Let $\mathcal{R}_0^k$ denote the contact area at the $k$th iteration ($k=1,2,...$). For a given $\mathcal{R}_0^{k-1}$ within the inner loop, we solve the variational equation of Eq.\,(\ref{var_eq_q}) using the Newton-Raphson iteration, and then we update the contact area by a contact search process, which is explained in Sections \ref{glob_csearch} and \ref{loc_csearch}. The outer loop continues until the \textit{active set} converges. The search procedure to find the active set has to be performed at every iteration step of the active set loop. It contains two cases:
\begin{itemize}
	\item Case 1: update the contact state of those contact pairs in the existing active set, and remove any inactive contact pairs,
	\item Case 2: add new contact pairs to the active set.
\end{itemize}
In the first case, we keep only the set elements (contact pairs) with positive contact pressure, i.e., $p_\mathrm{N}=-{\epsilon_\mathrm{N}}{g_\mathrm{N}}>0$, and remove the others from the active set. In the second case, new set elements are added to the active set if $g_\mathrm{N}<0$. For a more efficient contact search, we divide the search procedure into two steps: \textit{global} and \textit{local} searches.
\begin{remark}
	One can employ a regularized penalty law, which updates the active set simultaneously during the Newton-Raphson iteration, instead of employing an additional outer loop. For example, quadratically regularized penalty laws are used by \cite{durville2012contact} and \cite{meier2016finite}. Although in our formulation the active set method requires an additional outer loop, it has the following advantages:
	\begin{itemize}\small
		\item It fixes the active set within the inner Newton-Raphson iteration, so that it gives an improved convergence behavior for the same size of load increment.
		\item It does not require an additional user-defined parameter like a regularization threshold. 
		\item The quadratic regularization leads to zero contact stiffness initially, leading to inaccuracies. If a positive regularization threshold is used, unphysical contact force may initially occur.
	\end{itemize}
\end{remark}
\subsection{Global contact search}
\label{glob_csearch}
\begin{figure}	
	\centering
	\begin{subfigure}[b] {0.325\textwidth} \centering
		\includegraphics[width=\linewidth]{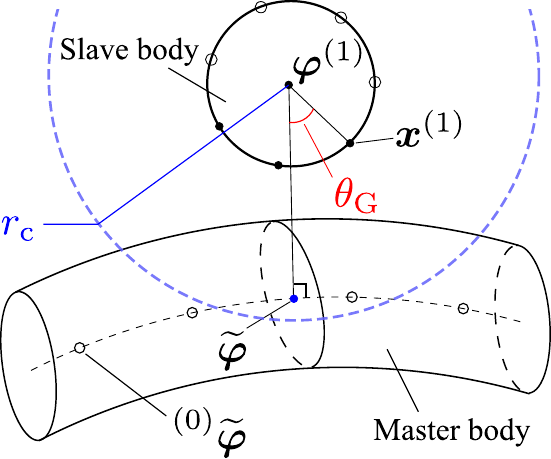}
	\end{subfigure}
	\caption{Global search for contact candidate surface Gauss points in the slave body using the normal distance between the beam axes. In the selection of initial guess, we first select several points (hollow dots) with uniform intervals of the parametric coordinate along the axis of the master body. From a selected initial guess ${}^{(0)}\widetilde{\boldsymbol{\varphi }}$, we find the closest point $\widetilde {\boldsymbol{\varphi}}$ on the axis, and we exclude the slave points ${\boldsymbol{x}^{(1)}}$, if the condition of Eq.\,(\ref{crad_1}) is not satisfied. The black hollow dots on the lateral surface of the slave body indicate those surface Gauss points excluded by the criterion of Eq.\,(\ref{gsearch_ang_tht}), and the other solid dots are those satisfying the criterion.}
	\label{contact_search_glob}
\end{figure}

The global search finds possible contact pairs using the \textit{normal distance} between the axes of interacting beams. We first find the convective axial coordinate $\widetilde \xi^1$ of the closest point ${\widetilde{\boldsymbol{\varphi}}}\coloneqq{\boldsymbol{\varphi}}^{(2)}\big(\widetilde \xi^{1}\big)$, for a given slave point $\boldsymbol{\varphi}^{(1)}$, where we use several initial guesses ${}^{(0)}{\widetilde {\boldsymbol{\varphi}}}$ along the axis of the master body, see Appendix \ref{app_glob_contact_search} for the detailed procedure, and Fig.\,\ref{contact_search_glob} for an illustration. Then we exclude those surface Gauss points corresponding to $\boldsymbol{\varphi}^{(1)}$, if ${\widetilde{\boldsymbol{\varphi}}}$ lies outside of the range
\begin{equation}
	\label{crad_1}
	\left\| \boldsymbol{\varphi}_\mathrm{d} \right\|\le{r_\mathrm{c}},
\end{equation}
where ${{\boldsymbol{\varphi }}_{\rm{d}}} \coloneqq {{\boldsymbol{\varphi }}^{(1)}} - {\widetilde{\boldsymbol{\varphi }}}$ defines a \textit{relative position vector} between the axis in slave body and the axis in master body, and $r_\mathrm{c}>0$ is a chosen cutoff radius. Selecting $r_\mathrm{c}$ too small may result in some contact areas being undetected, leading to unphysical penetration. As $r_\mathrm{c}$ increases, the contact search procedure becomes less efficient. It may also depend on cross-sectional dimensions in the current configurations as well, so that it might need to be adjusted, if cross-sectional deformations are very large. In this paper, we use a fixed parameter $r_\mathrm{c} = 3R$ for an initially circular cross-section of radius $R$, and the development of an efficient algorithm to adjust the parameter $r_\mathrm{c}$ remains future work. We further exclude those surface Gauss points that lie outside of the angular range
\begin{subequations}
	\begin{equation}
		\label{gsearch_ang_tht}
		0 \le {\theta _{{\rm{G}}}} \le {\varepsilon _\theta }
	\end{equation}
	with
	\begin{equation}
		\label{ang_tht_g_cr}
		{\theta _{{\rm{G}}}} \coloneqq {\cos ^{ - 1}}\left( {{{\boldsymbol{\varphi }}^{(1)}}-{{\boldsymbol{x}}^{(1)}}} \right) \cdot {{\boldsymbol{\varphi }}_{\rm{d}}},
	\end{equation}
\end{subequations}
and $0 < {\varepsilon _\theta } \le \pi\,\,[\mathrm{rad}]$. As $\varepsilon_\theta$ increases, the more Gauss points in the cross-sections's boundary of the slave body are considered as contact candidates, and the contact search becomes less efficient. On the other hand, too small $\varepsilon_\theta$ may lead to undetected contact areas or oscillation in the active set iteration, see the relevant discussion in Section \ref{ex_btb_slide_cs1}. The subsequent local contact search considers only those surface Gauss points of the slave body which satisfy the criteria of Eqs.\,(\ref{crad_1}) and (\ref{gsearch_ang_tht}).
%
\subsection{Local contact search}
\label{loc_csearch}
The local contact search is an iterative procedure to solve a Gauss point-to-surface minimal distance problem using a local Newton-Raphson iteration.
\subsubsection{Unilateral minimal distance problem}
\label{sec_ulat_mind_p}
We determine the distance between two surface points from
\begin{align}
	d\big({{\boldsymbol{x}}^{(1)}},{{\boldsymbol{x}}^{(2)}}({\boldsymbol{\xi}}^{(2)})\big) \coloneqq \left\| {{{\boldsymbol{x}}^{(1)}} - {{\boldsymbol{x}}^{(2)}}({\boldsymbol{\xi}}^{(2)})} \right\|.
\end{align}
Here and hereafter, we often use ${\boldsymbol{x}}^{(2)}\equiv{\boldsymbol{x}}^{(2)}({\boldsymbol{\xi}})\equiv{{\boldsymbol{x}}^{(2)}({\boldsymbol{\xi}}^{(2)})}$ for brevity. The parametric coordinates of the closest master point ${\bar {\boldsymbol{\xi}}}\equiv{\bar {\boldsymbol{\xi}}}(\boldsymbol{\xi}^{(1)})$ to a given surface point ${{\boldsymbol{x}}^{(1)}} \equiv {{\boldsymbol{x}}^{(1)}}({\boldsymbol{\xi}}^{(1)})$ in the slave body is determined as the solution of the following \textit{unilateral} minimal distance problem:
\begin{subequations}
	\label{min_d_xi}
	\begin{equation}
		\label{contact_min_disp_p_ul}
		{\boldsymbol{\bar \xi }} = \arg \mathop {\min }\limits_{{{\boldsymbol{\xi }}^{(2)}}} d\big({{\boldsymbol{x}}^{(1)}},{{\boldsymbol{x}}^{(2)}}{({\boldsymbol{\xi}}^{(2)})}\big),
	\end{equation}
	where
	\begin{equation}
		\label{side_cnst_xi12}
		\xi _{\min }^\alpha  \le \xi^{\alpha}_{(2)} \le \xi _{\max }^\alpha,\,\,\alpha\in\{1,2\},
	\end{equation}
\end{subequations}
and, in this paper, we normalize both convective coordinates such that $\xi^1_\mathrm{min}=\xi^2_\mathrm{min}=0$ and $\xi^1_\mathrm{max}=\xi^2_\mathrm{max}=1$. As we consider a closed curve for the boundary of the cross-section, the coordinate $\xi^2$ should have a \textit{periodic} property. Therefore, for $\alpha=2$, we employ
\begin{subequations}
	\label{cyclic_xi2}
	\begin{equation}
		\xi^\alpha_{(2)} = p\big(\xi^\alpha_{(2)},{\xi^\alpha_\mathrm{min}},{\xi^\alpha_\mathrm{max}}\big)
	\end{equation}
	with
	\begin{equation}
		{p}(\xi^\alpha,\xi^\alpha_\mathrm{min},\xi^\alpha_\mathrm{max}) = {\xi ^\alpha } - \left\lfloor {\frac{{{\xi ^\alpha } - \xi _{\min }^\alpha }}{{{\xi^\alpha_\mathrm{r}}}}} \right\rfloor {\xi _{\mathrm{r}}^\alpha},
	\end{equation}
	if ${\xi ^\alpha } > \xi _{\max }^\alpha$, and 
	\begin{equation}
		{p}(\xi^\alpha,\xi^\alpha_\mathrm{min},\xi^\alpha_\mathrm{max}) = {\xi ^\alpha } + \left\lfloor {\frac{{\xi _{\max }^\alpha  - {\xi ^\alpha }}}{\xi^\alpha_\mathrm{r}}} \right\rfloor {\xi^\alpha_\mathrm{r}},
	\end{equation}
\end{subequations}
if ${\xi ^\alpha } < \xi _{\min }^\alpha$, where $\xi^\alpha_\mathrm{r}\coloneqq{\xi^\alpha_\mathrm{max}}-{\xi^\alpha_\mathrm{min}}$, and $\left\lfloor {{\xi ^\alpha }} \right\rfloor$ denotes the largest integer smaller than $\xi ^\alpha$. If the axis curve is \textit{closed}, Eq.\,(\ref{cyclic_xi2}) can also be applied for $\alpha=1$. Otherwise, we divide the constraint of Eq.\,(\ref{side_cnst_xi12}) for $\alpha=1$ into two parts: \textit{domain} (i.e., ${\xi^1}\in(\xi^1_\mathrm{min},\xi^1_\mathrm{max})$) and \textit{end edges} (i.e., $\xi^1\in\left\{\xi^1_\mathrm{min},\xi^1_\mathrm{max}\right\}$) and treat those two cases separately. In this paper, we focus on the former as the latter is not relevant to the envisioned examples. But in general both \textit{domain-to-edge} or the \textit{edge-to-edge} contacts might contribute significantly in some applications including an arbitrary arrangement or configuration of beams. For example, in the simulation of biopolymer networks of \citet{meier2016finite}, it turns out that the \textit{end point-to-curve} and \textit{end point-to-end point} contacts play a significant role.
%
The solution of Eq.\,(\ref{contact_min_disp_p_ul}) satisfies the first order necessary condition
\begin{subequations}
	\label{cpp_uni_lat}
	\begin{equation}
		\label{cpp_unit_lat_a_f}
		{\boldsymbol{f}}({\boldsymbol{\xi}^{(2)}})\equiv\left\{ {\begin{array}{*{20}{c}}
				{{f_1}}\\
				{{f_2}}
		\end{array}} \right\}=\boldsymbol{0},
	\end{equation}
	with
	\begin{align}
		\label{cpp_uni_lat_f_a}		
		f_\alpha &\coloneqq \big({{{\boldsymbol{x}}^{(1)}} - {{\boldsymbol{x}}^{(2)}}(\boldsymbol{\xi}^{(2)})}\big) \cdot {\boldsymbol{a}}_\alpha^{(2)}(\boldsymbol{\xi}^{(2)}),
	\end{align}
	$\alpha\in\left\{1,2\right\}$, where
	\begin{equation}
		\label{side_cnst_xi1_fnc}
		\xi _{\min }^1  \lt {\xi _{(2)}^{1}} \lt \xi _{\max }^1.
	\end{equation}
\end{subequations}
Hereafter, for brevity, we often use $\boldsymbol{a}_\alpha^{(2)}\equiv{\boldsymbol{a}_\alpha^{(2)}}(\boldsymbol{\xi}^{(2)})$, $\alpha\in\left\{1,2\right\}$. The first order necessary condition of Eq.\,(\ref{cpp_uni_lat}) in the minimal distance problem finds a local extremum solution. In order to find the closest point, it is required to select an initial guess \textit{sufficiently close} to the closest point. We present a geometrical approach to efficiently determine an initial guess in Section \ref{ig_nr_cpp_range_cs_bd}. Within the Newton-Raphson iteration, for given coordinates $\boldsymbol{\xi}^{(2)}_{(i-1)}$ at the $(i-1)$th iteration, we find the increment $\Delta \boldsymbol{\xi}^{(2)}$ such that
\begin{equation}
	\label{cpp_uni_lat_linearized}
	{\boldsymbol{f}}^*({{\boldsymbol{\xi }}^{(2)}_{{(i - 1)}}})\,\Delta{{\boldsymbol{\xi }}^{(2)}} = -{{\boldsymbol{f}}}({{\boldsymbol{\xi }}^{(2)}_{(i - 1)}}),\,i=1,2,...,
\end{equation}
and the convective coordinates are updated by
\begin{equation}
	\label{cpp_update_add}
	{\boldsymbol{\xi }}^{(2)}_{{(i)}} = {\boldsymbol{\xi }}^{(2)}_{{(i - 1)}} + \Delta {\boldsymbol{\xi }}^{(2)},
\end{equation}
until the condition $\left\| {{\boldsymbol{f}}\big( {{{\boldsymbol{\xi }}^{(2)}_{(i)}}} \big)} \right\| < {\varepsilon _{{\mathrm{cpp}}}}$ is satisfied, where $\varepsilon_\mathrm{cpp}>0$ denotes a chosen tolerance, and $\boldsymbol{\xi}^{(2)}_{(0)}$ is an initial guess. Assuming ${\boldsymbol{f}}^*({{\boldsymbol{\xi }}^{(2)}_{{(i - 1)}}})$ is invertible, we obtain
\begin{equation}
	\label{cpp_uni_lat_linearized_sol_inv}
	\Delta{{\boldsymbol{\xi }}^{(2)}} = -{\boldsymbol{f}}^*({{\boldsymbol{\xi }}^{(2)}_{{(i - 1)}}})^{-1}{{\boldsymbol{f}}}({{\boldsymbol{\xi }}^{(2)}_{(i - 1)}}),	
\end{equation}
where
\begin{subequations}
	\label{f_u_st_1}
	\begin{equation}
		%
		{\boldsymbol{f}}^*({{\boldsymbol{\xi }}^{(2)}})= \left[ {\begin{array}{*{20}{c}}
				{f_{11}^ *}&{f_{12}^ * }\\
				\mathrm{sym.}&{f_{22}^ *}
		\end{array}} \right],
	\end{equation}
	with
	\begin{align}
		{f_{\alpha\beta}^ *}&\coloneqq{\frac{\partial f_\alpha}{\partial{\xi^\beta_{(2)}}}}\nonumber\\
		&={{{({{\boldsymbol{x}}^{(1)}} - {{\boldsymbol{x}}^{(2)}})}} \cdot {\boldsymbol{a}}_{\alpha,\beta}^{(2)} - {a}_{\alpha\beta}^{(2)}},
	\end{align}
\end{subequations}
$\alpha,\beta\in\left\{1,2\right\}$. Hereafter, the notation $\overline {(\bullet)}$ denotes the variable at $\boldsymbol{\xi}^{(2)}=\boldsymbol{\bar \xi}$ in the master body, unless otherwise stated. For example, ${\boldsymbol{\bar x}} \coloneqq {{\boldsymbol{x}}^{(2)}}({\boldsymbol{\bar \xi }})$. The overall procedure of the closest point projection is given in Algorithm \ref{proc_algo_lcsearch}. In the following we explain the determination of an initial guess $\boldsymbol{\xi}^{(2)}_{(0)}\equiv{{\bar {\boldsymbol{\xi}}}}_{{(0)}}=\left[{\bar {\xi}}^{1}_{(0)},{\bar {\xi}}^{2}_{(0)}\right]^\mathrm{T}$.
\subsubsection{Determination of an initial guess}
\label{ig_nr_cpp_range_cs_bd}
We select the intersection point between the relative position vector ${{\boldsymbol{\varphi }}_{\rm{d}}}$ and the cross-section at $\widetilde{\boldsymbol{\varphi}}$ as an initial guess for the local contact search. However, the intersection point, may not always exist on the boundary of the current cross-section at $\widetilde{\boldsymbol{\varphi}}$ due to the following two aspects:
\begin{itemize}
	\item The relative position vector is always normal to the axis of the master body; however, it may not be on the cross-section at $\widetilde{\boldsymbol{\varphi}}$ due to the transverse shear deformation.
	\item For $N>1$, the cross-section in the current configuration may not be planar.
\end{itemize}
Thus, we first consider a projection of the domain of the current cross-section ${\mathcal{A}^{(2)}_t}$ onto its tangent plane. This projected domain $\widetilde{\mathcal{A}_t}$ can be expressed by 
\begin{equation}
	{\widetilde{{\mathcal{A}_t}}} \coloneqq \left\{ {\left. {{\boldsymbol{x}} \in {{\Bbb{R}}^3}} \right\rvert{\boldsymbol{x}} = {\zeta ^\gamma }{{\boldsymbol{d}}^{(2)}_\gamma },\,({\zeta ^1},{\zeta ^2}) \in \mathcal{A}^{(2)}} \right\}.
\end{equation}
Note that for $N=1$, $\mathcal{A}^{(2)}_t\equiv{\widetilde{\mathcal{A}_t}}$. For convenience, we define $\boldsymbol{d}_3\coloneqq{{\boldsymbol{d}}_1} \times {{\boldsymbol{d}}_2}/\left\| {{{\boldsymbol{d}}_1} \times {{\boldsymbol{d}}_2}} \right\|$, and reciprocal base vectors $\boldsymbol{d}^i$ on the tangent plane such that ${{\boldsymbol{d}}_i} \cdot {{\boldsymbol{d}}^j} = \delta _i^j$ $(i,j\in\left\{1,2,3\right\})$, as
\begin{equation}
	\left\{ {\begin{array}{*{20}{c}}
			\begin{aligned}
				{{\boldsymbol{d}}^1} &\coloneqq {{\boldsymbol{d}}_2} \times {{\boldsymbol{d}}_3}/{{\boldsymbol{d}}_3} \cdot \left( {{\boldsymbol{d}}_1} \times {{\boldsymbol{d}}_2} \right),\\
				{{\boldsymbol{d}}^2} &\coloneqq {{\boldsymbol{d}}_3} \times {{\boldsymbol{d}}_1}/{{\boldsymbol{d}}_3} \cdot \left( {{\boldsymbol{d}}_1} \times {{\boldsymbol{d}}_2} \right),\\
				{{\boldsymbol{d}}^3} &\equiv {{\boldsymbol{d}}_3}.
			\end{aligned}
	\end{array}} \right.
\end{equation}
We also project the relative position vector $\boldsymbol{\varphi}_\mathrm{d}$ onto the tangent plane, as\footnote{The projection of the relative position vector was also used for the contact search in \citet{durville2012contact}.}
\begin{equation}
	\widetilde {{{\boldsymbol{\varphi }}_{\mathrm{d}}}} \coloneqq \left( {{\boldsymbol{1}} - {\boldsymbol{d}}_3^{(2)} \otimes {\boldsymbol{d}}_3^{(2)}} \right){{\boldsymbol{\varphi }}_{\mathrm{d}}}.
\end{equation}
Then, $\widetilde{\boldsymbol{\varphi}_\mathrm{d}}$ always intersects the boundary of $\widetilde{\mathcal{A}_t}$, denoted as $\partial{\widetilde{\mathcal{A}_t}}$, and we find the intersection between $\partial\widetilde{{\mathcal{A}}_t}$ and the projected vector $\alpha_\mathrm{ig}{\widetilde{\boldsymbol{\varphi}_\mathrm{d}}}$ with the length adjusted by a parameter $\alpha_\mathrm{ig}>0$, and the intersection point is selected as an initial guess for the closest point in the local contact search. That is, for a given $\boldsymbol{d}^\beta_{(2)}$ and $\widetilde{\boldsymbol{\varphi}_\mathrm{d}}$, we find the convective circumferential coordinate ${\xi^{2}_{(2)}}={\bar \xi^{2}_{(0)}}\in\left[\xi^2_\mathrm{min},\xi^2_\mathrm{max}\right]$, and $\alpha_\mathrm{ig}>0$ such that
\begin{equation}
	\label{intersect_cond_ipd}	
	e^\beta\coloneqq{\zeta ^\beta}\big({\bar \xi^{2}_{(0)}}\big) - {\boldsymbol{d}}^\beta_{(2)}\cdot\alpha_\mathrm{ig}{\widetilde{\boldsymbol{\varphi}_\mathrm{d}}}=0,\,\,\beta\in\left\{1,2\right\}.
\end{equation}
Since Eq.\,(\ref{intersect_cond_ipd}) is nonlinear with respect to ${\bar \xi^{2}_{(0)}}$, we need an iterative solution process. It should be noted that we use ${\bar \xi^{1}_{(0)}}={\widetilde \xi^1}$. Further details can be found in Appendix \ref{lsearch_ig_det}.
\begin{remark}
For typical shapes of the initial cross-section, $\zeta^\gamma$ satisfying Eq.\,(\ref{intersect_cond_ipd}) can be analytically found. For example, for an initially circular cross-section of radius $R$ such that
\begin{equation}
	\label{cir_cs_zta12_rel}
\big(\zeta^1\big)^2+\big(\zeta^2\big)^2=R^2,
\end{equation}
we have
\begin{equation}
	\label{sol_circ_zta12}
	{\zeta ^\gamma }\big({\bar \xi^{2}_{(0)}}\big) = \frac{{R\,{\boldsymbol{d}}_{(2)}^\gamma  \cdot {\widetilde{{{{\boldsymbol{\varphi }}}_{\rm{d}}}}}}}{{\sqrt {{{\big( {{\boldsymbol{d}}_{(2)}^1 \cdot {\widetilde{{{\boldsymbol{\varphi }}}_{\rm{d}}}}} \big)}^2} + {{\big( {{\boldsymbol{d}}_{(2)}^2 \cdot {\widetilde{{{\boldsymbol{\varphi }}}_{\rm{d}}}}} \big)}^2}} }},
\end{equation}
$\gamma\in\left\{1,2\right\}$. The process to find ${\bar \xi_{(0)}^{2}}$ satisfying Eq.\,(\ref{sol_circ_zta12}) requires an iterative process due to the parameterization of the coordinates $\big(\zeta^1,\zeta^2\big)$ in Eq.\,(\ref{plane_zeta12_curve_nurbs}). Instead, in this paper, we directly solve Eq.\,(\ref{intersect_cond_ipd}), where both ${\bar \xi_{(0)}^{2}}$ and $\alpha_\mathrm{ig}$ are determined simultaneously.
\end{remark}
\begin{figure}	
	\centering
	\begin{subfigure}[b] {0.3125\textwidth} \centering
		\includegraphics[width=\linewidth]{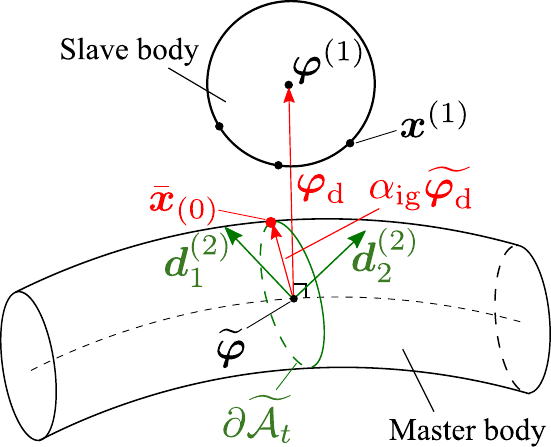}
	\end{subfigure}
	\caption{A geometrical approach of selecting an initial guess ${\bar {\boldsymbol{x}}}_{(0)}\equiv{\boldsymbol{x}}^{(2)}({\bar {\boldsymbol{\xi}}_{(0)}})$ in the local contact search as the intersection of the vector $\alpha_\mathrm{ig}\widetilde{\boldsymbol{\varphi}_\mathrm{d}}$ with the boundary of the projected cross-section $\partial{\widetilde{\mathcal{A}_t}}$.}
	\label{contact_lsearch_ig}
\end{figure}
\subsection{Variational formulation}
The first variation of the position vector of the closest point on the master surface, i.e., ${\bar{\boldsymbol{x}}} \equiv{\bar{\boldsymbol{x}}}(\bar {\boldsymbol{\xi}})$ is obtained by \citep{wriggers2006}
\begin{equation}
	\label{fvar_pos_cmaster_p}
	\delta {\bar{\boldsymbol{x}}} \coloneqq\frac{d}{{d\varepsilon }}{\left. {{{{\boldsymbol{\bar x}}}_\varepsilon }({\bar {\boldsymbol{\xi}}}_\varepsilon)} \right\vert_{\varepsilon  = 0}}= \delta{\bar{\boldsymbol{u}}} + {{\bar{\boldsymbol{a}}_{\alpha }}}\,\delta {\bar \xi ^\alpha },
\end{equation}
where $\delta{\bar{\boldsymbol{u}}}$ defines the directional derivative of $\bar{\boldsymbol{x}}$ with the dependence of the parametric coordinate $\bar {\boldsymbol{\xi}}$ on $\varepsilon$ suppressed, i.e., 
\begin{equation}
	\delta{\bar{\boldsymbol{u}}}\coloneqq\frac{d}{{d\varepsilon }}{\left. {{{{\boldsymbol{\bar x}}}_\varepsilon }({\bar {\boldsymbol{\xi}}})} \right\rvert_{\varepsilon  = 0}}.
\end{equation}
For beams, evaluating Eq.\,(\ref{compact_xt_pi}) at ${\boldsymbol{\xi}}={\bar{\boldsymbol{\xi}}}$ gives the position vector of the closest point on the surface $\mathcal{R}^{(2)}_t$
\begin{equation}
	\label{cpp_x_pi_q}
	{\bar{\boldsymbol{x}}} = {\bar{\boldsymbol{\Pi}}}^\mathrm{T}{\bar{\boldsymbol{q}}},
\end{equation}
where $\bar{\boldsymbol{x}}\coloneqq{\boldsymbol{x}}^{(2)}(\bar{\boldsymbol{\xi}})$, ${\bar {\boldsymbol{\Pi}}}\coloneqq{{\boldsymbol{\Pi}}}^{(2)}(\zeta^1(\bar{\boldsymbol{\xi}}),\zeta^2(\bar{\boldsymbol{\xi}}))$, and ${\bar{\boldsymbol{q}}}\coloneqq{\boldsymbol{q}^{(2)}(s({\bar{\xi}^1}))}$. Eq.\,(\ref{fvar_pos_cmaster_p}) can thus be rewritten as
\begin{equation}
	\label{del_x_bar_fvar_xbar}
	\delta {\boldsymbol{\bar x}} = {\bar{\boldsymbol{\Pi}}}^\mathrm{T}\delta{\widetilde{\boldsymbol{q}}} + {{\bar{\boldsymbol{a}}_{\alpha }}}\,\delta {\bar \xi ^\alpha },
\end{equation}
where we define
\begin{equation}
	\label{del_delta_w_bar_sup}
	\delta {\boldsymbol{\widetilde q}} \coloneqq \frac{d}{{d\varepsilon }}{\left. {{{{\boldsymbol{\bar q}}}_\varepsilon }({\bar \xi^1})} \right\rvert_{\varepsilon  = 0}}.
\end{equation}
For a slave point, we simply obtain
\begin{equation}
	\label{del_x1_fvar_x_slave}
	\delta {\boldsymbol{x}}^{(1)} = {{\boldsymbol{\Pi}}_{(1)}^{\mathrm{T}}}\delta{{\boldsymbol{q}}}^{(1)}.
\end{equation}
At ${\boldsymbol{\xi}}^{(2)}={\bar{\boldsymbol{\xi}}}$ for a given $\boldsymbol{x}^{(1)}$, by the definition of $g_\mathrm{N}$ in Eq.\,(\ref{cond_impenetrate_normal}), we have 
\begin{equation}
	\label{diff_x1_x2_at_x_bar}
	{{\boldsymbol{x}}^{(1)}} - {{\boldsymbol{\bar x}}} = {g_{\mathrm{N}}}\,{\boldsymbol{\bar \nu }}_t.
\end{equation}
Taking the first variation of Eq.\,(\ref{diff_x1_x2_at_x_bar}), and substituting Eqs.\,(\ref{del_x_bar_fvar_xbar}) and (\ref{del_x1_fvar_x_slave}), we obtain
\begin{align}
	\label{var_normal_gap_prv}
	\delta {g_{\mathrm{N}}}{\boldsymbol{\bar \nu }}_t &= {{\boldsymbol{\Pi}}_{(1)}^{\mathrm{T}}}\delta {{\boldsymbol{q}}^{(1)}} - {\bar{\boldsymbol{\Pi}}}^\mathrm{T}\delta {{\boldsymbol{\widetilde q}}} - {\boldsymbol{\bar a}}_{\alpha }\,\delta {\bar \xi}^{\alpha} \nonumber\\
	&- {g_{\rm{N}}}\,\delta {\boldsymbol{\bar \nu}}_t.
\end{align}
Then, by taking the inner product of Eq.\,(\ref{var_normal_gap_prv}) and ${\boldsymbol{\bar \nu}}_t$, and using ${{\boldsymbol{\bar \nu}}}_t \cdot {{\boldsymbol{\bar \nu }}}_t = 1$ and $\delta {{\boldsymbol{\bar \nu}}}_t \cdot {{\boldsymbol{\bar \nu }}}_t = {{\boldsymbol{\bar \nu}}}_t \cdot {{\boldsymbol{\bar a}}_\alpha} = 0$, we have
\begin{equation}
	\label{var_normal_gap}
	\delta {g_{\rm{N}}} = {\delta {{\boldsymbol{q}}_{(1)}^{\mathrm{T}}}}{{{\boldsymbol{\Pi }}^{(1)}}}{{\boldsymbol{\bar \nu }}}_t
	- {\delta {\boldsymbol{\widetilde q}}}^\mathrm{T}{{\bar{\bf \Pi }}}\,{{\boldsymbol{\bar \nu }}}_t.
\end{equation}
The internal virtual work due to the contact pressure can be written as
\begin{align}
	{G_{\rm{N}}}({\boldsymbol{q}},\delta {\boldsymbol{q}}) &\coloneqq \int_{\mathcal{R}_0} { - {p_{\rm{N}}}\,\delta {g_{\rm{N}}}\,{\rm{d}}\mathcal{R}_0} \nonumber\\
	&= \int_{\mathcal{S}_0^\mathrm{L}} { - {p_{\rm{N}}}\,\delta {g_{\rm{N}}}\,{\rm{d}}\mathcal{S}_0^\mathrm{L}} 	\label{vwork_ncontact_gn}.		
\end{align}
Using Eqs.\,(\ref{inf_area_init_config}) and (\ref{var_normal_gap}), Eq.\,(\ref{vwork_ncontact_gn}) can be rewritten into
\begin{subequations}
	\label{contvar_form_gn}
	\begin{align}
		{G_{\rm{N}}}({\boldsymbol{q}},\delta {\boldsymbol{q}}) =  \int_0^{{L^{(1)}}} {\left\{ {\setlength{\arraycolsep}{0.1pt}
				\renewcommand{\arraystretch}{1}\begin{array}{*{20}{c}}
					{\delta {{\boldsymbol{q}}^{(1)}}}\\
					{\delta {\boldsymbol{\widetilde q}}}
			\end{array}} \right\} \cdot {{{\boldsymbol{r}}_{\rm{N}}}}\,{\rm{d}}s},\label{reg_first_var_tot_pot_2}
	\end{align}
with the \textit{normal contact stress resultant} vector
	\begin{align}
		\label{cforce_dist_rn_slave}
		{{\boldsymbol{r}}_{\rm{N}}} &\coloneqq -\frac{1}{{{{\tilde j}^{(1)}}}}\int_{{\varXi_{(1)}^2}}\!{{{{p}}_{\rm{N}}}\,{\tilde J}^{(1)}{\boldsymbol{\Pi }}^{(1)}{{\boldsymbol{\bar \nu }}_t}\,{\rm{d}}{\xi ^2}}\nonumber\\
		&+ \frac{1}{{{{\tilde j}^{(1)}}}}\int_{{\varXi_{(1)}^2}}\!{{{{p}}_{\rm{N}}}\,{\tilde J}^{(1)}{\bar {\boldsymbol{\Pi}}}\,{{\boldsymbol{\bar \nu }}_t}\,{\rm{d}}{\xi ^2}},
	\end{align}
\end{subequations}
where $\varXi^2\coloneqq[\xi^2_\mathrm{min},\xi^2_\mathrm{max})$ denotes the domain of the parametric coordinate $\xi^2$ in the cross-section's boundary curve.
\begin{remark} \textit{Units of the contact stress resultant}. Each component of the stress resultant $\boldsymbol{r}_\mathrm{N}$, which is energy conjugate to the $n$th order director, has units of ${F_0}{L_0}^{n-1}$, where $F_0$ and $L_0$ denote a unit force and unit length, respectively. For example, those components of $\boldsymbol{r}_\mathrm{N}$ corresponding to $n=0$ represent a resultant force per unit undeformed arc-length, and those corresponding to $n=1$ represent a resultant director moment per unit undeformed arc-length. 
\end{remark}
\subsection{Linearization of the contact variational form}
In order to solve the nonlinear variational equation, we linearize the contact variational form of Eq.\,(\ref{vwork_ncontact_gn}). Taking the directional derivative of Eq.\,(\ref{vwork_ncontact_gn}) gives
\begin{align}
	\label{del_gn_inc_del}
	&\Delta{G_{\rm{N}}}({\boldsymbol{q}};\!\delta {\boldsymbol{q}},\Delta {\boldsymbol{q}}) \coloneqq \nonumber\\
	&\int_{\mathcal{S}_0^\mathrm{L}} {{\epsilon}_{\mathrm{N}}\left(\delta {g_{\mathrm{N}}}\,\Delta {g_{\mathrm{N}}}+{g}_{\mathrm{N}}\,\Delta \delta {g_{\mathrm{N}}}\,\right)\omega\,\mathrm{d}{\mathcal{S}_0^\mathrm{L}}},
\end{align}
where we employ the Heaviside function $\omega\equiv\omega(\boldsymbol{X})$, defined as
\begin{equation}\label{gn_bracket}
	{{{\omega}}}\coloneqq \left\{ {\begin{array}{*{20}{c}}
			\begin{array}{l}
				1\\
				0
			\end{array}&\begin{array}{l}
				{\rm{if}}\,\,\,{\boldsymbol{X}} \in \mathcal{R}_0,\\
				{\rm{if}}\,\,\,{\boldsymbol{X}} \in {\mathcal{S}_0^\mathrm{L}}\setminus\mathcal{R}_0.
			\end{array}
	\end{array}} \right.
\end{equation}
Further, using Eq.\,(\ref{var_normal_gap}), we have
\begin{align}
	\label{del_gn_del_gn}
	\delta {g_{\rm{N}}}\Delta {g_{\rm{N}}} = {\left\{ {\begin{array}{*{20}{c}}
				{\delta {{\boldsymbol{q}}^{(1)}}}\\
				{\delta {\boldsymbol{\widetilde q}}}
		\end{array}} \right\}^{\rm{T}}}\boldsymbol{k}^\mathrm{M}_\mathrm{N}\left\{ {\begin{array}{*{20}{c}}
			{\Delta {{\boldsymbol{q}}^{(1)}}}\\
			{\Delta {\boldsymbol{\widetilde q}}}
	\end{array}} \right\},
\end{align}
where
\begin{align}
	\boldsymbol{k}^\mathrm{M}_\mathrm{N}&\coloneqq\left[ {\begin{array}{*{20}{c}}
			{{{\boldsymbol{\Pi }}^{(1)}}{{{\boldsymbol{\bar \nu }}}_t} \otimes {{\boldsymbol{\Pi }}^{(1)}}{{{\boldsymbol{\bar \nu }}}_t}}&{ - {{\boldsymbol{\Pi }}^{(1)}}{{{\boldsymbol{\bar \nu }}}_t} \otimes {\boldsymbol{\bar \Pi }}{{{\boldsymbol{\bar \nu }}}_t}}\\
			{\mathrm{sym.}}&{{\boldsymbol{\bar \Pi }}{{{\boldsymbol{\bar \nu }}}_t} \otimes {\boldsymbol{\bar \Pi }}{{{\boldsymbol{\bar \nu }}}_t}}
	\end{array}} \right].
\end{align}
We obtain the increment of $\delta {g_\mathrm{N}}$ by taking the directional derivative of Eq.\,(\ref{var_normal_gap_prv})\footnote{The directional derivative of the vanishing terms in Eq.\,(\ref{var_normal_gap_prv}) can also contribute to the tangent stiffness \citep{wriggers2006}.} and applying Eq.\,(\ref{del_x_bar_fvar_xbar}), as (see Appendix \ref{app_lin_contact_var_form} for details)
\begin{align}
	\label{del_del_g_N_beam}
	\Delta \delta {g_{\mathrm{N}}} = {\left\{ {\begin{array}{*{20}{c}}
				{\delta {{\boldsymbol{q}}^{(1)}}}\\
				{\delta \widetilde{\boldsymbol{q}}}
		\end{array}} \right\}^{\mathrm{T}}}{{\boldsymbol{k}}^\mathrm{G}_{\mathrm{N}}}\left\{ {\begin{array}{*{20}{c}}
			{\Delta {{\boldsymbol{q}}^{(1)}}}\\
			{\Delta \widetilde{\boldsymbol{q}}}
	\end{array}} \right\},
\end{align}
where the symmetric matrix $\boldsymbol{k}^\mathrm{G}_\mathrm{N}$ is given by Eq.\,(\ref{deriv_kg_n}). Finally, substituting Eqs.\,(\ref{del_gn_del_gn}) and (\ref{del_del_g_N_beam}) into Eq.\,(\ref{del_gn_inc_del}), we have
\begin{align}
	\label{cont_form_inc_cont_var}
	&\Delta{G_{\rm{N}}}({\boldsymbol{q}};\!\delta {\boldsymbol{q}},\Delta {\boldsymbol{q}})= \nonumber\\
	&\int_0^{{L^{(1)}}} {{{\left\{\setlength{\arraycolsep}{0.125pt}{\begin{array}{*{20}{c}}
						{\delta {{\boldsymbol{q}}^{(1)}}}\\
						{\delta \widetilde{\boldsymbol{q}}}
				\end{array}} \right\}}^{\mathrm{T}}}{{\boldsymbol{k}}}_\mathrm{N}\left\{ {\setlength{\arraycolsep}{0.125pt}\begin{array}{*{20}{c}}
				{\Delta {{\boldsymbol{q}}^{(1)}}}\\
				{\Delta {\widetilde{\boldsymbol{q}}}}
		\end{array}} \right\}\,\mathrm{d}s},
\end{align}
with
\begin{align}
	\label{def_kn_stiff_ncont}
	{{\boldsymbol{k}}}_\mathrm{N}\coloneqq\frac{1}{{{{\tilde j}^{(1)}}}}\int_{{{\varXi^{2}_{(1)}}}}\!{\epsilon_\mathrm{N}}{({{\boldsymbol{k}}^\mathrm{M}_{\mathrm{N}}}+{{g_{\mathrm{N}}}\,{{\boldsymbol{k}}^\mathrm{G}_{\mathrm{N}}}})\,{\tilde J}^{(1)}\!\omega\,\mathrm{d}{\xi^2}},
\end{align}
which is symmetric.
\section{Isogeometric finite element discretization}
\label{sec_spatial_disc_iga}
\subsection{Beam formulation}
We discretize the contact variational form and its increment using NURBS basis functions. A discussion on the crucial properties of NURBS in isogeometric analysis can be found in \cite{hughes2005isogeometric}. The geometry of the beam's initial axis is described by a NURBS curve, as
\begin{equation}
	{\boldsymbol{\varphi}}_0({\xi ^1}) = \sum\limits_{I = 1}^{{n_{{\mathrm{cp}}}}} {{N_I}({\xi ^1})\,{{\boldsymbol{P}}_I}},
\end{equation}
where ${n_{{\mathrm{cp}}}}$ denotes the total number of basis functions (or control points) of the axis, and ${{\boldsymbol{P}}_I}\in\Bbb{R}^3$ denotes the position vector of the $I$th control point. Using the NURBS basis functions $N_I=N_I(\xi^1)$, the
finite element approximation $\delta {{\boldsymbol{q}}^h}=\delta{{\boldsymbol{q}}^h}(s(\xi^1))$ is expressed by
\begin{align}
	\label{nurbs_disc_general_dir_var}
	\delta {{\boldsymbol{q}}^h} &= \left[ {\begin{array}{*{20}{c}}
			{N_1{{\boldsymbol{1}}_{{n_\mathrm{cs}}}}},\cdots,{N_{{n_{e}}}{{\boldsymbol{1}}_{{n_\mathrm{cs}}}}}
	\end{array}} \right]\left\{ {\setlength{\arraycolsep}{0.1pt}
		\renewcommand{\arraystretch}{1}\begin{array}{*{20}{c}}
			{\delta {\bf{q}}_1}\\
			\vdots \\
			{\delta {\bf{q}}_{{n_{\rm{e}}}}}
	\end{array}} \right\}\nonumber\\
	&\eqqcolon{\Bbb{N}_e}(\xi^1)\,\delta{{\bf{q}}^{e}},
\end{align}
where ${\delta{\bf{q}}_I} \in {{\Bbb{R}}^{n_{{\rm{cs}}}}}$ denotes the coefficient vector of the $I$th control point ($I\in\left\{1,...,n_e\right\}$), and ${n_e}$ denotes the number of basis functions having local support in the knot span ${\varXi _e}$ with $e \in \left\{ {1,...,{n_{{\mathrm{el}}}}} \right\}$, and $n_\mathrm{el}$ denotes the total number of nonzero knot spans. Substituting Eq.\,(\ref{nurbs_disc_general_dir_var}) into Eq.\,(\ref{intv_work}), and using the standard finite element assembly operator $\bf{A}$, we have
\begin{equation}
	\label{disc_int_vwork_q}
	{G_{{\mathop{\rm int}} }}({{\boldsymbol{q}}^h},\delta {{\boldsymbol{q}}^h}) = \delta {{\bf{q}}^{\rm{T}}}{{\bf{F}}_{{\mathop{\rm int}} }}
\end{equation}
with ${{\bf{F}}_{{\mathop{\rm int}} }} \coloneqq \mathop {\bf{A}}_{e = 1}^{{n_{{\rm{el}}}}} {\bf{F}}_{{\mathop{\rm int}} }^e$, and $\delta\bf{q}$ denotes the global coefficient vector of the generalized directors. The element internal load vector is obtained by
\begin{equation}
	{\bf{F}}_{{\mathop{\rm int}} }^e \coloneqq \int_{{\varXi _e}} {{{\bf{R}}^e}\,{\tilde j}\,\mathrm{d}\xi^1},
\end{equation}
where we define
\begin{equation}
	{{\bf{R}}^e} \coloneqq \int_\mathcal{A} {{\bf{\bar \Xi }}{{_e^h}^{\rm{T}}}{\boldsymbol{S}}\,{j_0}\,\mathrm{d}\mathcal{A}},
\end{equation}
with ${\bf{\bar \Xi }}_e^h \coloneqq {\left[ {{\bf{\Xi }}_1^h, \cdots ,{\bf{\Xi }}_{{n_e}}^h} \right]_{6 \times {n_{{\rm{cs}}}{n_e}}}}$, and
\begin{equation}
	{\bf{\Xi }}_I^h \coloneqq {\left[ \setlength{\arraycolsep}{1pt}
		\renewcommand{\arraystretch}{1.25}{\begin{array}{*{20}{c}}
				{{{\boldsymbol{q}}^{\rm{T}}}{\boldsymbol{\Pi }}_{,{\zeta ^1}}\,{{\boldsymbol{\Pi }}^{\rm{T}}_{,{\zeta ^1}}}{N_I}}\\
				{{{\boldsymbol{q}}^{\rm{T}}}{\boldsymbol{\Pi }}_{,{\zeta ^2}}\,{{\boldsymbol{\Pi }}^{\rm{T}}_{,{\zeta ^2}}}{N_I}}\\
				{{\boldsymbol{q}}_{,s}^{\rm{T}}\,{{\boldsymbol{\Pi }}}{\boldsymbol{\Pi }}^{\rm{T}}\,{N_{I,s}}}\\
				{{{\boldsymbol{q}}^{\rm{T}}}\left({\boldsymbol{\Pi }}_{,{\zeta ^2}}\,{{\bf{\Pi }}^{\rm{T}}_{,{\zeta ^1}}}+{\boldsymbol{\Pi }}_{,{\zeta ^1}}\,{{\bf{\Pi }}^{\rm{T}}_{,{\zeta ^2}}}\right){N_I}}\\
				{{\boldsymbol{q}}^{\rm{T}}_{,s}\,{{\bf{\Pi }}}{{\bf{\Pi }}_{,{\zeta ^1}}^{\rm{T}}}{N_I} + {{\boldsymbol{q}}^{\rm{T}}}{\bf{\Pi }}_{,{\zeta ^1}}\,{\boldsymbol{\Pi }}^{\rm{T}}\,{N_{I,s}}}\\
				{{\boldsymbol{q}}_{,s}^{\rm{T}}\,{{\boldsymbol{\Pi }}}\,{{\boldsymbol{\Pi }}^{\rm{T}}_{,{\zeta ^2}}}{N_I} + {{\boldsymbol{q}}^{\rm{T}}}{\boldsymbol{\Pi }}_{,{\zeta ^2}}\,{\bf{\Pi }}^{\rm{T}}\,{N_{I,s}}}
		\end{array}} \right]}.
\end{equation}
It is noted that, for brevity, we often use \citep{choi2021isogeometric}
\begin{equation}
	N_{I,s}\coloneqq{N_{I,1}}\frac{\mathrm{d}\xi^1}{\mathrm{d}s}=\frac{1}{\tilde j}{N_{I,1}}.
\end{equation}
The external virtual work of Eq.\,(\ref{ext_vwork_q}) is discretized as
\begin{equation}
	\label{disc_ext_vwork_q}
	{G_{{\rm{ext}}}}(\delta {{\boldsymbol{q}}^h}) = \delta {{\bf{q}}^{\rm{T}}}{{\bf{F}}_{{\rm{ext}}}},
\end{equation}
with ${{\bf{F}}_{{\rm{ext}}}} \coloneqq \mathop {\bf{A}}_{e = 1}^{{n_{{\rm{el}}}}} {\bf{F}}_{{\rm{ext}}}^e$, where
\begin{equation}
	{\bf{F}}_{{\rm{ext}}}^e \coloneqq \int_{{\varXi _e}} {{\Bbb{N}}_e^{\rm{T}}{\boldsymbol{\bar R}}\,\,{\tilde j}\,\mathrm{d}\xi^1}.
\end{equation}
Similarly, the increment of the internal virtual work is discretized as
\begin{equation}
	\Delta {G_{{\mathop{\rm int}} }}({{\boldsymbol{q}}^h};\delta {{\boldsymbol{q}}^h},\Delta {{\boldsymbol{q}}^h}) = \delta {{\bf{q}}^{\rm{T}}}{{\bf{K}}_{{\mathop{\rm int}} }}\Delta {\bf{q}},
\end{equation}
with ${{\bf{K}}_{{\mathop{\rm int}} }}\coloneqq\mathop {\bf{A}}_{e = 1}^{{n_{{\rm{el}}}}} {\bf{K}}_{{\mathop{\rm int}} }^e$.
The element tangent stiffness matrix is
\begin{equation}
	{\bf{K}}_{{\mathop{\rm int}} }^e \coloneqq \int_{{\varXi _e}} {\left( {{{\Bbb{C}}_e} + {\Bbb{Y}}_e^{\rm{T}}{{\boldsymbol{k}}_{\rm{G}}}{{\Bbb{Y}}_e}} \right)\tilde j{\kern 1pt} {\rm{d}}\xi^1},
\end{equation}
where
\begin{equation}
	{{\Bbb{C}}_e} \coloneqq \int_\mathcal{A} {{\bf{\bar \Xi }}{{_e^h}^{\rm{T}}}{\boldsymbol{\munderbar{\munderbar{\mathcal{C}}}}}\,\,{\bf{\bar \Xi}}_e^h\,{j_0}\,\mathrm{d}\mathcal{A}},
\end{equation}
and ${{\Bbb{Y}}_e} \coloneqq {\left[ {\begin{array}{*{20}{c}}
			{{{\bf{Y}}_1}}, \cdots, {{{\bf{Y}}_{{n_e}}}}
	\end{array}} \right]_{2n_{{\rm{cs}}} \times n_{{\rm{cs}}}{n_e}}}$ with ${{\bf{Y}}_I} \coloneqq {\left[ {\begin{array}{*{20}{c}}
			{N_I}{{{\bf{1}}_{{n_{{\rm{cs}}}}}}}, {N_{I,s}}{{{\bf{1}}_{{n_{{\rm{cs}}}}}}}
	\end{array}} \right]^{\rm{T}}}$. It is noted that the global tangent stiffness matrix ${{\bf{K}}_{{\mathop{\rm int}} }}$ is symmetric, since the matrix ${{\Bbb{C}}_e}$ and ${\boldsymbol{k}}_\mathrm{G}$ are symmetric.
\subsection{Beam contact formulation}
Let $\mathcal{W}_k$ denote the active (working) set of indices of the surface Gauss point at the $k$th iteration of the outer loop, defined by 
\begin{align}
	{\mathcal{W}_k}\!\coloneqq\!\left\{ {\left. {i}\!\in\!\{ {1,...,n^\mathrm{L}_\mathrm{G}}\} \right\rvert{p^i_{\rm{N}}} \coloneqq-\epsilon_\mathrm{N}{g^i_{\rm{N}}} \gt 0 }\right\}.
\end{align}
Hereafter, we often use ${{\boldsymbol{\xi}}_i}\coloneqq[\xi_i^{1},\xi_i^{2}]^\mathrm{T}$, which denotes the parametric coordinates of the $i$th surface Gauss integration point of the slave body, and ${\bar{\boldsymbol{\xi}}}_i\coloneqq{\bar{\boldsymbol{\xi}}}({{\boldsymbol{\xi}}_i})$ with ${\bar{\boldsymbol{\xi}}}_i\coloneqq[{\bar{\xi}^1_i},{\bar{\xi}^2_i}]^\mathrm{T}$. We also use $g^i_{\rm{N}}\!\coloneqq\!{g_{\rm{N}}}({\boldsymbol{x}}^{(1)}\!({{\boldsymbol{\xi }}_i}),{{\bar{\boldsymbol{x}}}}({\bar{\boldsymbol{\xi}}}_i))$, and let $n^\mathrm{L}_\mathrm{G}$ denote the total number of Gauss integration points on the lateral surface of the slave body. $w^\alpha_i$ denotes the weight of the Gauss integration point in the domain ${\varXi_{(1)}^\alpha}\ni{\xi^\alpha}$, $\alpha\in\left\{1,2\right\}$. Then, from Eq.\,(\ref{cforce_dist_rn_slave}), we define 
\begin{equation}
	{\boldsymbol{r}^i_\mathrm{N}} \coloneqq {\frac{{p_\mathrm{N}^i}{w^2_i}\,{\tilde J_i^{(1)}}}{{{{\tilde j}_i^{(1)}}}}}\left\{ {\setlength{\arraycolsep}{0.01pt}
		\renewcommand{\arraystretch}{1.3}\begin{array}{*{20}{c}}
			{ - {{{\boldsymbol{\Pi }}^{(1)}}}\!(\zeta_i^1,\zeta_i^2)\,{\bar{\boldsymbol{\nu}}}^i_t}\\
			{{\bf{\Pi }}}^{(2)}\!({\bar \zeta}_i^1,{\bar \zeta}_i^2)\,{\bar{\boldsymbol{\nu}}}^i_t
	\end{array}} \right\},
\end{equation}
$i\!\in\!{\mathcal{W}_k}$, with ${\tilde J^{(1)}_i}\coloneqq{\tilde J^{(1)}}{({{\boldsymbol{\xi}}_i})}$, ${{{\tilde j_i^{(1)}}}}\coloneqq{{{\tilde j^{(1)}}\!(\xi_i^1)}}$, ${{\zeta}}^\alpha_i\coloneqq\zeta^\alpha({\boldsymbol{\xi}}_i)$, ${\bar\zeta^\alpha_i}\coloneqq\zeta^\alpha({\bar{\boldsymbol{\xi}}}_i)$, $\alpha\in\left\{1,2\right\}$, and ${\bar{\boldsymbol{\nu}}_t^{i}}\coloneqq{\boldsymbol{\nu}_t^{(2)}}\big({\bar{\boldsymbol{\xi}}_i}\big)$. Substituting Eq.\,(\ref{nurbs_disc_general_dir_var}) into the contact variational form of Eq.\,(\ref{contvar_form_gn}) gives
\begin{equation}
	\label{disc_var_pi_cont_var_form}
	{G_{\rm{N}}}({\boldsymbol{q}}^h,\delta {\boldsymbol{q}}^h) = \delta {{\bf{q}}^{\mathrm{T}}}{{\bf{F}}_{\mathrm{N}}},
\end{equation}
with
\begin{equation}
	{{\bf{F}}_{\mathrm{N}}} \coloneqq \mathop {\blds{\bf{A}}}_{i \in {\mathcal{W}^k}} {{{\mathbbm{r}}^i_{\rm{N}}}\,{\tilde j^{(1)}_i}\,{w^1_i}},
\end{equation}
where
\begin{equation}
	{\mathbbm{r}^i_\mathrm{N}} \coloneqq {\frac{{p_\mathrm{N}^i}{w^2_i}\,{\tilde J_i^{(1)}}}{{{{\tilde j}_i^{(1)}}}}}\left\{ {\setlength{\arraycolsep}{0.01pt}
		\renewcommand{\arraystretch}{1.3}\begin{array}{*{20}{c}}
			{ - {\Bbb{N}_{e}^{(1)\,\mathrm{T}}}\,{{{\boldsymbol{\Pi }}^{(1)}}}\!(\zeta_i^1,\zeta_i^2)\,{\bar{\boldsymbol{\nu}}}^i_t}\\
			{{\bar {\Bbb{N}}}^\mathrm{T}_{\bar e}}\,{{\bf{\Pi }}}^{(2)}\!({\bar \zeta}_i^1,{\bar \zeta}_i^2)\,{\bar{\boldsymbol{\nu}}}^i_t
	\end{array}} \right\},
\end{equation}
with ${\Bbb{N}_{e}^{(1)}}\equiv{{\Bbb{N}_{e}^{(1)}}}(\xi^1_i)$, and ${{\bar {\Bbb{N}}}_{\bar e}}\coloneqq{{{\Bbb{N}}}^{(2)}_{\bar e}}({\bar \xi}^1_i)$. $e$ and $\bar e$ represent the indices of the nonzero knot spans including the coordinates ${\xi}_i^{1}$ and ${\bar \xi}_i^1$, respectively, i.e., $\varXi^{1(1)}_{e}\ni{{\xi}_i^1}$ and $\varXi^{1(2)}_{\bar e}\ni{{\bar \xi}_i^1}$. Similarly, the increment of the contact variational form of Eq.\,(\ref{cont_form_inc_cont_var}) can be discretized as
\begin{equation}
	\Delta{G_{\rm{N}}}({\boldsymbol{q}}^h;\!\delta {\boldsymbol{q}}^h,\Delta {\boldsymbol{q}}^h) = \delta {{\bf{q}}^{\mathrm{T}}}{{\bf{K}}_{\mathrm{N}}}\Delta {\bf{q}},
\end{equation}
where
\begin{equation}
	{{\bf{K}}_{\mathrm{N}}} \coloneqq \mathop {\blds{\bf{A}}}_{i \in {\mathcal{W}^k}}{{\Bbbk^i_\mathrm{N}}{({\boldsymbol{\xi}}_i})\,{{{\tilde j_i^{(1)}}}}}\,{w^1_i}.
\end{equation}
From Eq.\,(\ref{def_kn_stiff_ncont}), we have
\begin{equation}
	{\Bbbk^i_\mathrm{N}}({\boldsymbol{\xi}}_i) \coloneqq {\frac{{\epsilon_\mathrm{N}}{w^2_i}{\tilde J_i^{(1)}}}{{{{\tilde j}_i^{(1)}}}}}{\left\{{{{\Bbbk}}^\mathrm{M}_{\mathrm{N}}({\boldsymbol{\xi}}_i)}+{g^i_{\mathrm{N}}}\,{{{\Bbbk}}^\mathrm{G}_{\mathrm{N}}}({\boldsymbol{\xi}}_i)\right\}},
\end{equation}
where $i\!\in\!{\mathcal{W}_k}$, and ${\Bbbk}^{\mathrm{M}}_\mathrm{N}$ and $\Bbbk^{\mathrm{G}}_\mathrm{N}$ are defined by Eqs.\,(\ref{disc_km_n}) and (\ref{disc_kg_n}), respectively.\\
\section{Numerical examples}
\label{num_ex_contact}
We verify the presented isogeometric finite element formulation for beams and frictionless contact by comparison with reference solutions based on isogeometric brick elements. In the latter, we employ an active set strategy for Gauss point-to-surface contact combined with a penalty regularization. For the brick element formulation, we denote the degrees of basis functions in each parametric coordinate direction by $\mathrm{deg.}=(p_\mathrm{L},p_\mathrm{W},p_\mathrm{H})$, where $p_\mathrm{L}$, $p_\mathrm{W}$, and $p_\mathrm{H}$ denote the degrees of basis functions along the length (L), width (W), and height (H), respectively. Further, the number of elements in each of those directions are indicated by ${n_\mathrm{el}}={n_\mathrm{el}^\mathrm{L}}\times{{n_\mathrm{el}^\mathrm{W}}}\times{{n_\mathrm{el}^\mathrm{H}}}$. In both beam and brick element formulations, for the contact integral, we employ a subdivision of elements in order to increase the accuracy of numerical integration. The number of sub-elements are given in each numerical example, and we always use a single Gauss integration point per sub-element, so that the number of sub-elements is the same as the number of Gauss integration points for the contact integral. In the beam formulation, for the evaluation of the contact integral along the circumferential direction, we subdivide the whole domain $\left[\xi^2_\mathrm{min},\xi^2_\mathrm{max}\right]\ni{\xi^2_{(1)}}$ of each curve patch into $m_\mathrm{el}^\mathrm{sub}$ sub-elements with uniform intervals. In the brick formulation, $m^\mathrm{sub}_\mathrm{el}$ denotes the number of sub-elements within each nonzero knot span (element) along the circumferential direction. In both beam and brick element formulations, the number of sub-elements in each element (nonzero knot span) along the axial direction is denoted by $n^\mathrm{sub}_\mathrm{el}$. Thus, in the beam formulation, the total numbers of sub-elements for the evaluation of contact integral is ${n_\mathrm{el}}\times{n_\mathrm{el}^\mathrm{sub}}$ and $m_\mathrm{el}^\mathrm{sub}$ in the axial and circumferential directions, respectively. In the brick element formulation, the total numbers of sub-elements in two transverse directions are ${n^\mathrm{W}_\mathrm{el}}\times{m^\mathrm{sub}_\mathrm{el}}$, and ${n^\mathrm{H}_\mathrm{el}}\times{m^\mathrm{sub}_\mathrm{el}}$, respectively. In the numerical examples, we consider a compressible Neo-Hookean material, which is the simplest type of hyperelastic material for arbitrarily large deformations. The St.\,Venant-Kirchhoff material, which provides a simpler formulation, is not valid for arbitrarily large deformations. Especially in the example of Section \ref{num_ex_twist}, we apply a large pre-strain in the axial direction, in which the St.\,Venant-Kirchhoff material might show unphysical decrease of volume for positive Poisson's ratio, as discussed in \cite{choi2021isogeometric}. In Sections \ref{num_ex_slide_2beam} and \ref{num_ex_twist}, the boundary curves of initial (circular) cross-sections are represented by a single patch of periodic NURBS with $p=3$, and $m_\mathrm{cp}=67$.
\subsection{Lateral contact of a straight beam under end moments}
\begin{figure}
	\centering
	\begin{subfigure}[b] {0.475\textwidth} \centering
		\includegraphics[width=\linewidth]{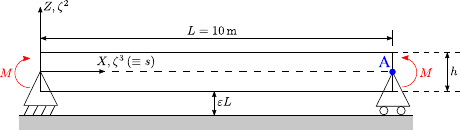}
		\caption{Undeformed configuration}
		\label{lat_contact_str_beam_undeform}
	\end{subfigure}		
	\begin{subfigure}[b] {0.45\textwidth} \centering			\includegraphics[width=\linewidth]{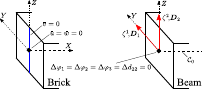}
		\caption{Hinge condition at the left end}
		\label{lat_contact_str_beam_bdc_brick_hinge}
	\end{subfigure}
	\begin{subfigure}[b] {0.375\textwidth} \centering				\includegraphics[width=\linewidth]{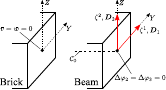}
		\caption{Roller condition at the right end}
		\label{lat_contact_str_beam_bdc_brick_roller}
	\end{subfigure}
	\caption{Lateral contact of a straight beam: (a) Undeformed configuration and boundary conditions, (b) hinge condition at the left end ($s=0$) for the brick and beam element models, (c) roller support condition at the right end ($s=L$) for the brick and beam element models. Here, $\bar u$, $\bar v$, and $\bar w$ denote the prescribed $X$-, $Y$-, and $Z$-displacements in the brick element model, respectively. The dashed line $\mathcal{C}_0$ shows the initial center axis of beam.}
	\label{lat_contact_str_beam_undeform_bdc}
\end{figure}
This example investigates the alleviation of Poisson locking in the beam formulation, and the verification of the contact pressure distribution in lateral contact to a rigid flat surface. We consider a straight beam with length $L=10\,\mathrm{m}$ and a rectangular cross-section of height $h=0.2\,\mathrm{m}$, and width $w=1\,\mathrm{m}$, and choose a Neo-Hookean material with Young's modulus $E=1.2\times{10^7}\,\mathrm{Pa}$, and the Poisson's ratio $\nu\!=\!0.25$. As shown in Fig.\,\ref{lat_contact_str_beam_undeform}, the straight beam is aligned in $X$-direction and it has an initial lateral distance ${\epsilon}L$ from the rigid flat surface, with chosen $\epsilon=2\times10^{-2}$. An end moment $M\!=\!\alpha\lambda_\mathrm{load}EI/L$ is applied, where $\lambda_\mathrm{load}$ denotes the load parameter, satisfying $0\le\lambda_\mathrm{load}\le1$, and $I$ denotes the second moment of inertia of the rectangular cross-section, obtained by $I=w{h^3}/12$. It is noted that for $\alpha=2\pi$ the two end points of the beam's axis meet, leading to an indeterminacy for rotation around the $Y$-axis, which is why we choose $\alpha=1.8\pi$. In this example, we use uniform load increments, and the load parameter is obtained by $\lambda_\mathrm{load}=n/{n_\mathrm{load}}$, where $n_\mathrm{load}$ denotes the total number of load step, and $n=0,1,...,n_\mathrm{load}$ denotes the load step number. We use $n_\mathrm{load}=20$ for both beam and brick elements. In both beam and brick element solutions, we use $n_\mathrm{el}^\mathrm{sub}=10$ and $m_\mathrm{el}^\mathrm{sub}=20$. For the beam formulation, a hinge condition is imposed at the left end ($s=0$), as $\Delta{\varphi}_{1} = \Delta{\varphi}_{2} = \Delta{\varphi}_{3} = 0$, and $\Delta{d}_{22}=0$ is further imposed in order to avoid rigid body rotation around the $X$-axis, where $\Delta\varphi_i\coloneqq{\Delta \boldsymbol{\varphi}\cdot{\boldsymbol{e}}_i}$ and ${\Delta}d_{{\alpha}i}\coloneqq\Delta{{\boldsymbol{d}}_{\alpha}}\cdot\boldsymbol{e}_i$ ($i\in\left\{1,2,3\right\}$, $\alpha\in\left\{1,2\right\}$). In the brick element formulation, at the left end, we also impose the condition ${\bar v}=0$ along the vertical line $Y=0$ in the initial configuration (see the vertical blue line in Fig.\,\ref{lat_contact_str_beam_bdc_brick_hinge}) in order to avoid the rigid body rotation around $X$-axis. The detailed formulation of the end moment condition can be found in Remark \ref{remark_imp_dfol}. In the brick element model, we apply the same traction boundary condition of Remark\,\ref{remark_imp_dfol} (see Eq.\,(\ref{end_moment_first_pk})) with ${\bar p}=-{MZ/I}$ ($- h/2\!\le\!Z\!\le\!h/2$) at both end faces at $X=0$ and $X=L$. At the right end of the axis ($s=L$), roller conditions are applied by ${\bar v}={\bar w}=0$, and $\Delta{\varphi_2}=\Delta{\varphi_3}=0$ in the brick and beam element formulations, respectively, see Fig.\,\ref{lat_contact_str_beam_bdc_brick_roller}. Since the surface basis functions do not have an interpolatory property in the domain of the end faces, these displacement boundary conditions along the selected points or lines in the initial configuration of the brick element model are imposed by using a penalty method with penalty parameter $\epsilon_\mathrm{D}=10^{7}E/{L_0}$, where $L_0$ denotes the unit length. In the beam formulation, however, those displacement boundary conditions can be imposed exactly since the basis functions using the clamped knot vectors in the axis satisfies the Kronecker-delta property at the ends. 
\subsubsection{Deformation without contact}
\begin{figure}[htp]
	\centering
	\begin{subfigure}[b] {0.475\textwidth} \centering
	\includegraphics[width=1\linewidth]{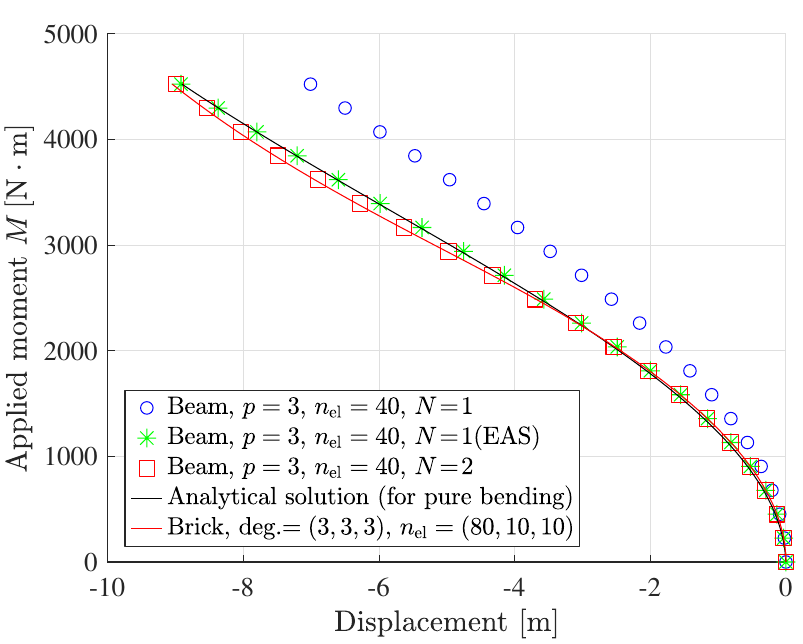}
		\caption{$\nu=0.25$}
		\label{ex1_tipa_disp_analytic_comp_ux_pr025}
	\end{subfigure}		
	\begin{subfigure}[b] {0.475\textwidth} \centering
		\includegraphics[width=1\linewidth]{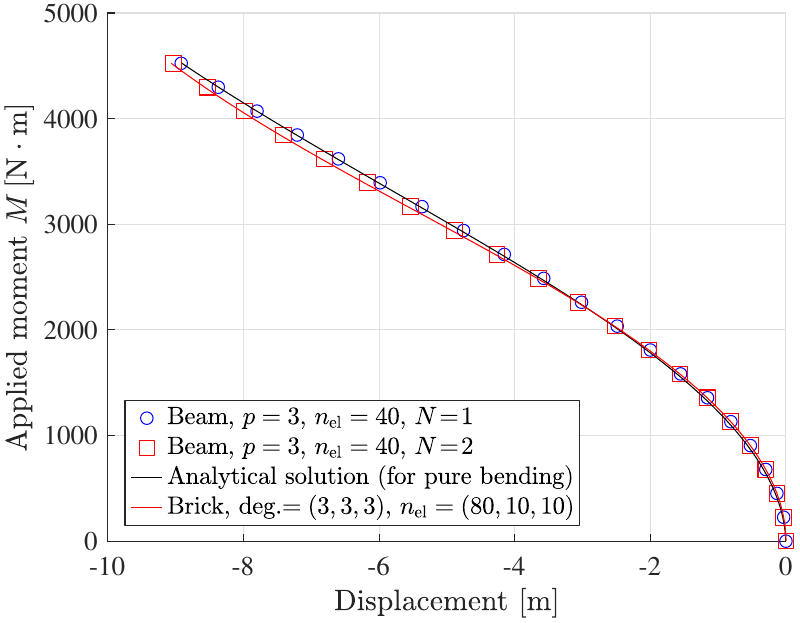}
		\caption{$\nu=0$}
		\label{ex1_tipa_disp_analytic_comp_ux_pr0}
	\end{subfigure}		
	\caption{Lateral contact of a straight beam: Comparison of the $X$-displacement at the point A (marked in Fig.\,\ref{lat_contact_str_beam_undeform}) for two values of Poisson's ratio: $\nu=0.25$ and $\nu=0$. No contact condition is imposed.}
	\label{ex1_tipa_disp_analytic_comp_ux}
\end{figure}

\begin{figure} \centering	
	\centering
	\begin{subfigure}[b] {0.5\textwidth} \centering
		\includegraphics[width=0.9975\linewidth]{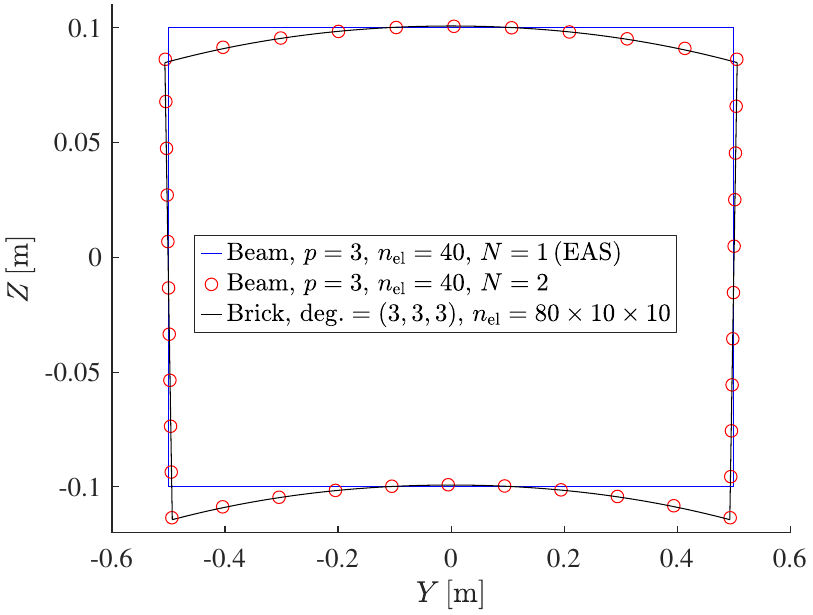}
		\caption{$\nu=0.25$}
	\label{ex_socket_cs_shape_comp_nu025}
	\end{subfigure}
	\begin{subfigure}[b] {0.5\textwidth} \centering
		\includegraphics[width=0.9975\linewidth]{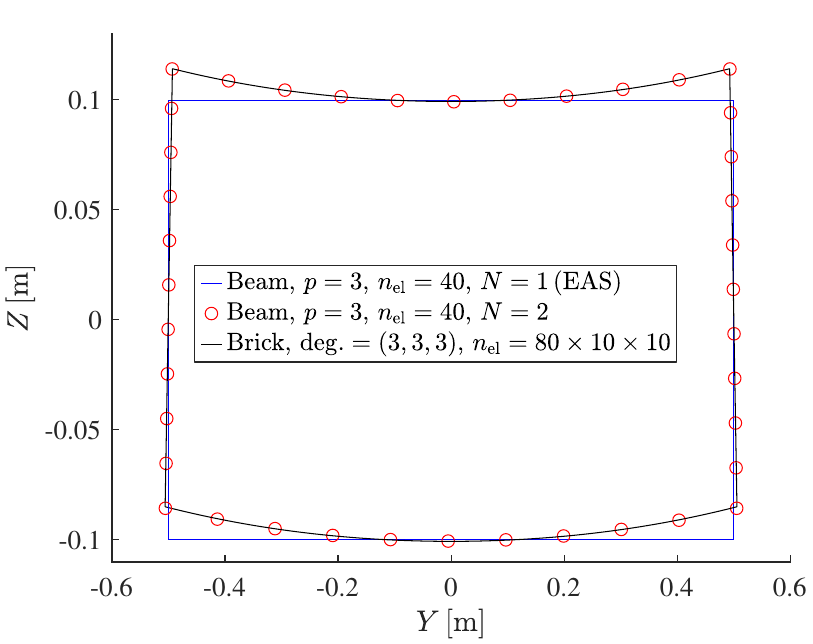}
		\caption{$\nu=-0.25$}
		\label{ex_socket_cs_shape_comp_nu-025}
	\end{subfigure}
	\caption{Lateral contact of a straight beam: Comparison of the deformed cross-section shape at the center of the beam ($s=0.5L$) in the final deformed configuration for the two cases $\nu=0.25$ and $\nu=-0.25$. No contact condition is imposed.}
	\label{ex_socket_cs_shape_comp}
\end{figure}

We first consider a problem without the impenetrability condition. Fig.\,\ref{ex1_tipa_disp_analytic_comp_ux} compares the $X$-displacement at the point A (marked in Fig.\,\ref{lat_contact_str_beam_undeform}) between beam element solutions and reference solutions for two different values of Poisson's ratios: $\nu=0.25$ and $\nu=0$. We have the following two reference solutions: 
\begin{itemize}
	\item Analytical solution: Under the assumption of \textit{pure bending}, the applied moment $M$ at both ends deforms the beam axis into a circle with radius $R=EI/M$, where the $X$-displacement at point A can be obtained as
	\begin{equation}
		u_\mathrm{A} = L\left\{\frac{\sin\left({{\alpha{\lambda}_\mathrm{load}/2}}\right)}{\alpha{\lambda}_\mathrm{load}/2} - 1\right\},
	\end{equation}
	for $0<\lambda_\mathrm{load}\le1$, and $u_\mathrm{A}=0$ if $\lambda_\mathrm{load}=0$,
	\item The brick element solution using B-spline basis functions of $\mathrm{deg.}=(3,3,3)$ and $n_\mathrm{el}=80\times10\times10$.
\end{itemize} 
The results using $N=1$ suffers from Poisson locking such that the bending stiffness is artificially increased due to the inability to represent linear in-plane strains of the cross-section. By employing the EAS method in \citet{choi2021isogeometric}, it is shown that the beam solution agrees very well with the analytical solution. However, the beam solution for $N=1$ (EAS) does not consider cross-sectional warping properly, so that it is shown to slightly deviate from the brick element solution. By increasing the order of approximation in the transverse direction to $N=2$, which enables to represent cross-sectional warping properly, the beam solution comes very close to the brick element solution. In case of $\nu=0$, there is no Poisson effect, so that the beam solution for $N=1$ already agrees very well with the analytical solution, but still slightly deviates from the brick element solution. Similar to the results in the case of $\nu=0.25$, the beam solution of $N=2$ agrees very well with the brick element solution. Fig.\,\ref{ex_socket_cs_shape_comp} compares the deformed cross-section shapes for the two different Poisson's ratios: $\nu=0.25$ and $\nu=-0.25$. It is shown that the beam solution for $N=1$ (EAS) always maintains straight boundaries, while the beam solution for $N=2$ properly represents the curved cross-section boundary due to the Poisson effect in the cases of $\nu=0.25$, and $\nu=-0.25$. This is essential for accurately capturing contact, as is shown next.
%
%
\subsubsection{Contact pressure distribution on the bottom surface}
\begin{figure*}[htp]
	\centering
	\begin{subfigure}[b] {0.31\textwidth} \centering
		\includegraphics[width=1\linewidth]{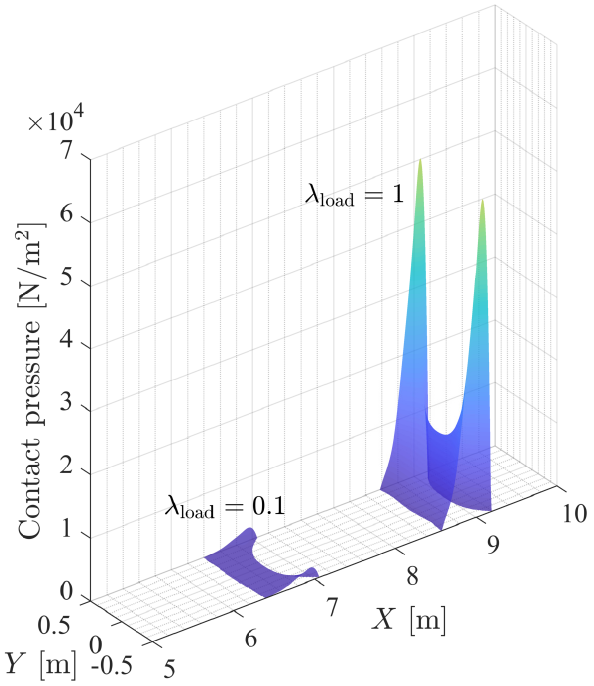}
		\caption{$\nu=0.25$, Brick}
		\label{ex1_socket_prs_surf_pr025_brick}
	\end{subfigure}		
	\begin{subfigure}[b] {0.31\textwidth} \centering
		\includegraphics[width=1\linewidth]{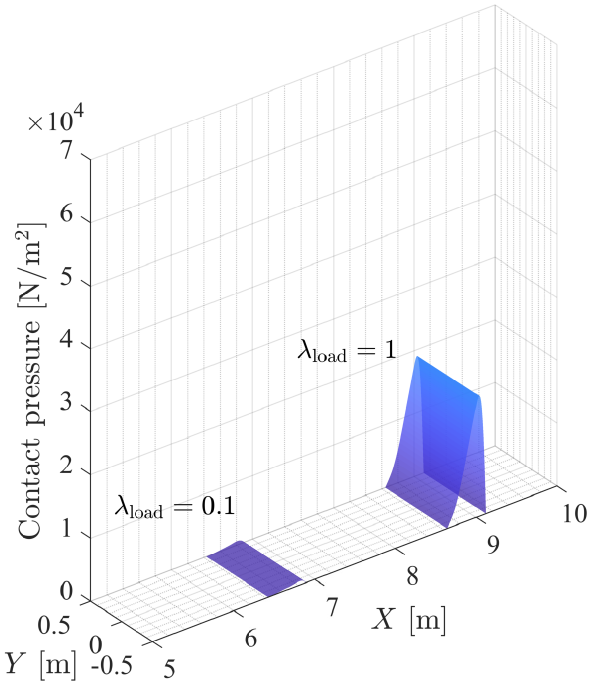}
		\caption{$\nu=0.25$, Beam, $N\!=\!1$ (EAS)}
		\label{ex1_socket_prs_surf_pr025_beam_n1eas}
	\end{subfigure}	
	\begin{subfigure}[b] {0.31\textwidth} \centering
		\includegraphics[width=1\linewidth]{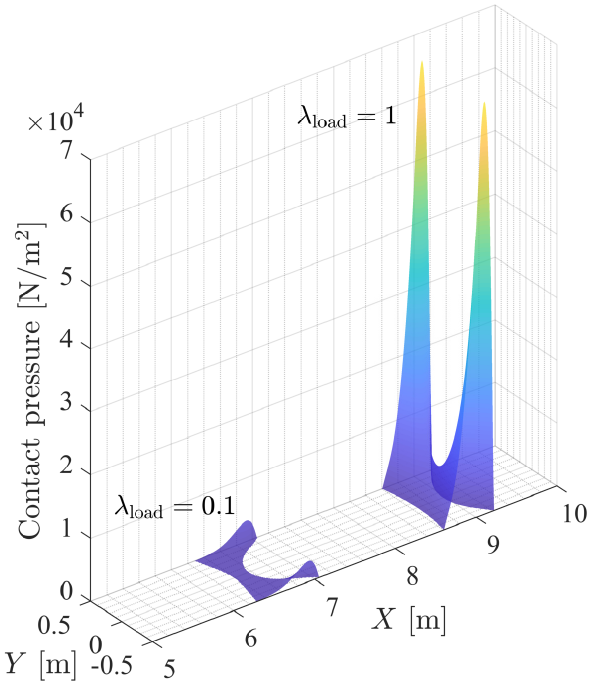}
		\caption{$\nu=0.25$, Beam, $N\!=\!2$}
		\label{ex1_socket_prs_surf_pr025_beam_n2}
	\end{subfigure}			
	\begin{subfigure}[b] {0.31\textwidth} \centering
		\includegraphics[width=1\linewidth]{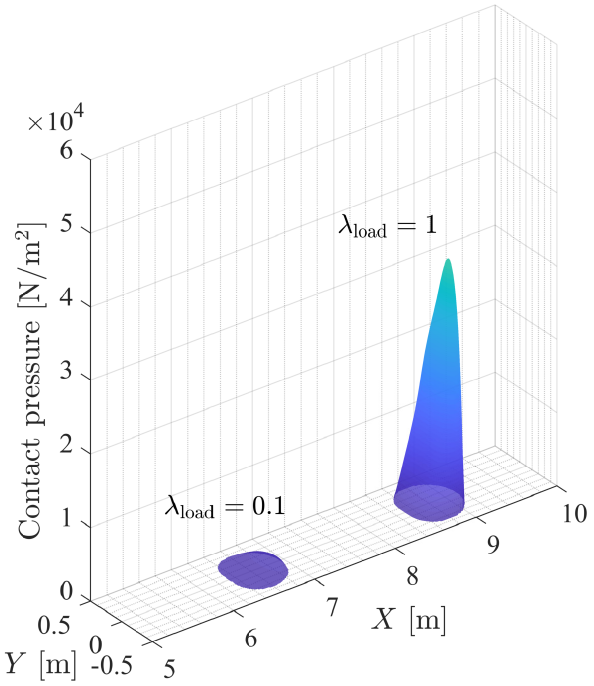}
		\caption{$\nu=-0.25$, Brick}
		\label{ex1_socket_prs_surf_pr-025_brick}
	\end{subfigure}		
	\begin{subfigure}[b] {0.31\textwidth} \centering
		\includegraphics[width=1\linewidth]{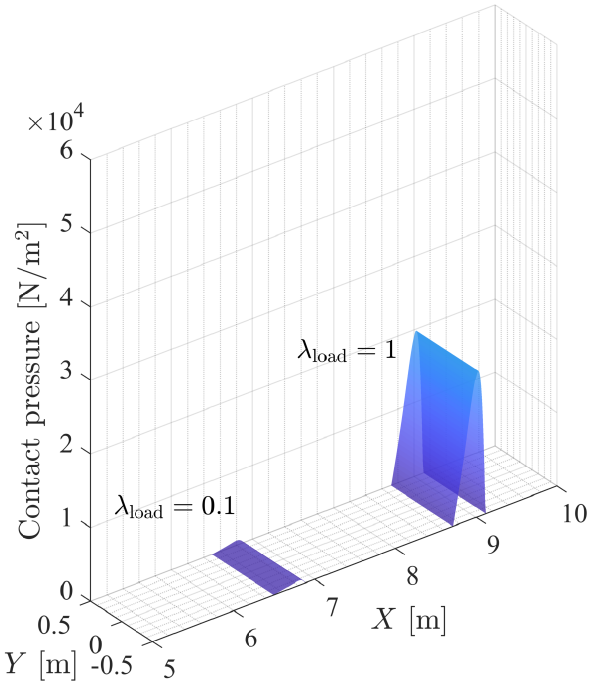}
		\caption{$\nu=-0.25$, Beam, $N\!=\!1$ (EAS)}
		\label{ex1_socket_prs_surf_pr-025_beam_N1EAS}
	\end{subfigure}	
	\begin{subfigure}[b] {0.31\textwidth} \centering
		\includegraphics[width=1\linewidth]{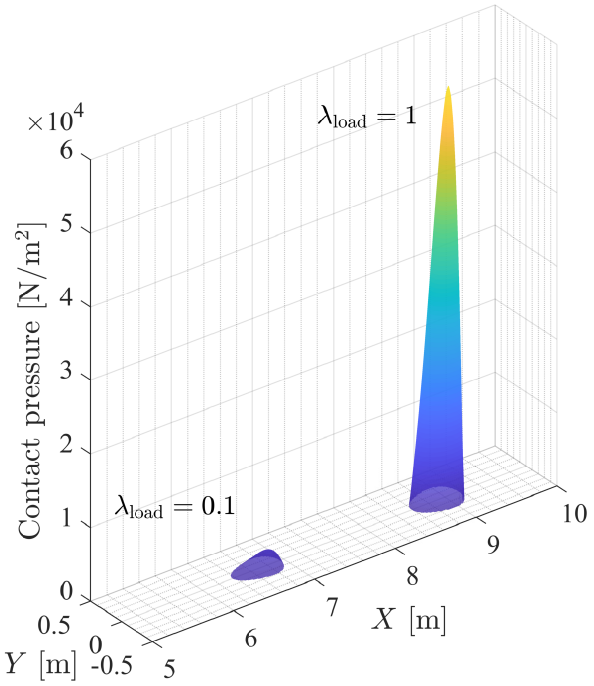}
		\caption{$\nu=-0.25$, Beam, $N\!=\!2$}
		\label{ex1_socket_prs_surf_pr-025_beam_N2}
	\end{subfigure}			
\begin{subfigure}[b] {0.31\textwidth} \centering
	\includegraphics[width=1\linewidth]{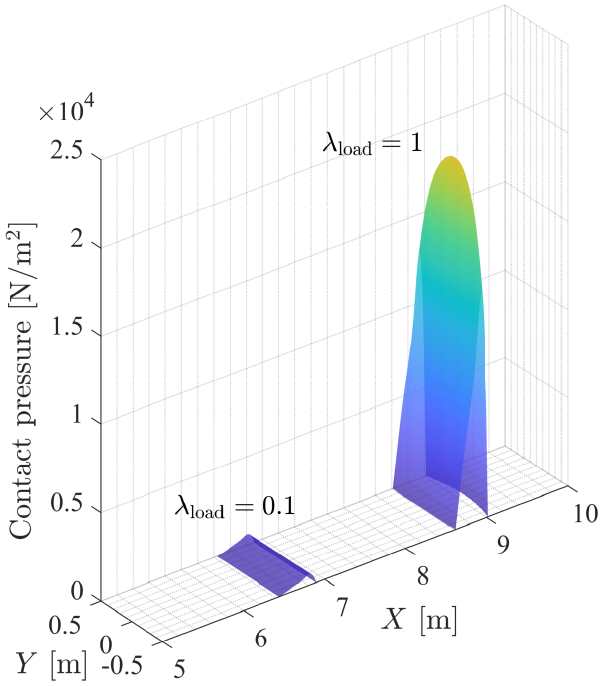}
	\caption{$\nu=0$, Brick}
	\label{ex1_socket_prs_surf_pr0_brick}
\end{subfigure}		
\begin{subfigure}[b] {0.31\textwidth} \centering
	\includegraphics[width=1\linewidth]{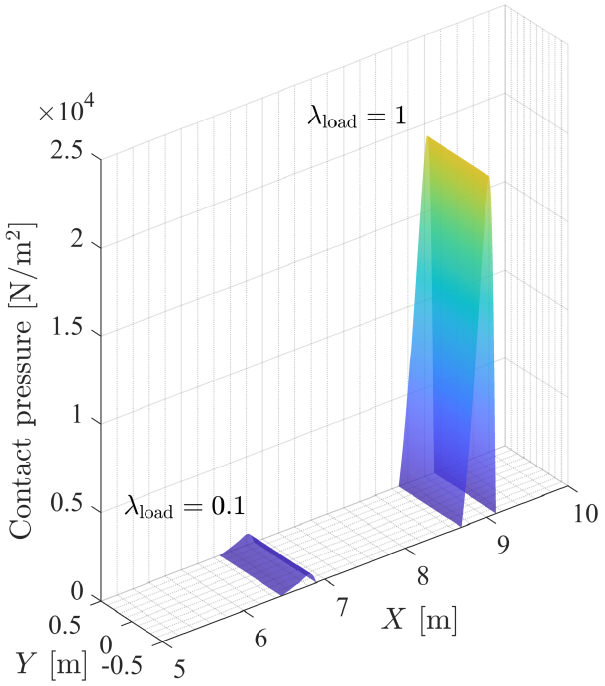}
	\caption{$\nu=0$, Beam, $N\!=\!1$}
	\label{ex1_socket_prs_surf_pr0_n1eas}
\end{subfigure}	
\begin{subfigure}[b] {0.31\textwidth} \centering
	\includegraphics[width=1\linewidth]{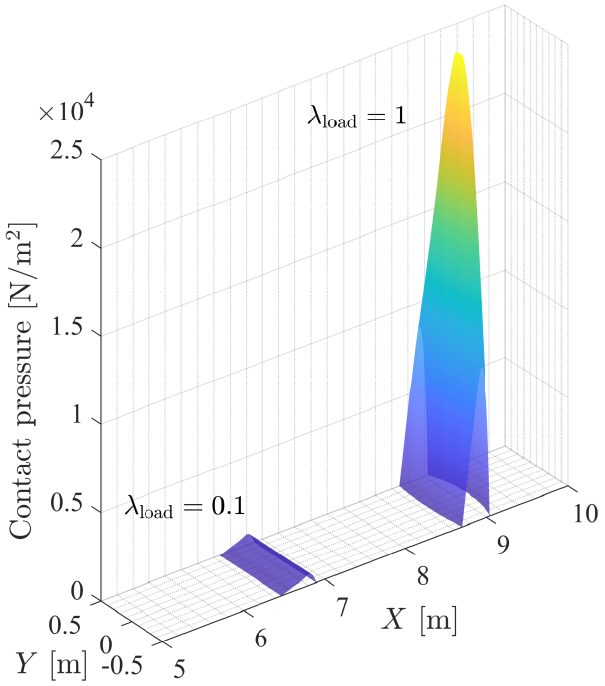}
	\caption{$\nu=0$, Beam, $N\!=\!2$}
	\label{ex1_socket_prs_surf_pr0_n2}
\end{subfigure}			
	\caption{Lateral contact of a straight beam: Distribution of the contact pressure on the bottom surface for three different values of Poisson's ratio. The penalty parameter is $\epsilon_\mathrm{N}=10^1E/{L_0}$ for all cases. The contact pressure is plotted in the initial configuration for two different load steps: $\lambda_\mathrm{load}=0.1$, and $1$ (not showing $p=0$). The colors correspond to the contact pressure values. We use brick elements with deg.=$(3,3,3)$, ${n_\mathrm{el}}=80\times10\times10$, and beam elements with $p=3$, ${n_\mathrm{el}}=80$.}
	\label{ex1_socket_pres_surf_diff_pr_case}
\end{figure*}

\begin{figure}
	\centering
	\includegraphics[width=0.95\linewidth]{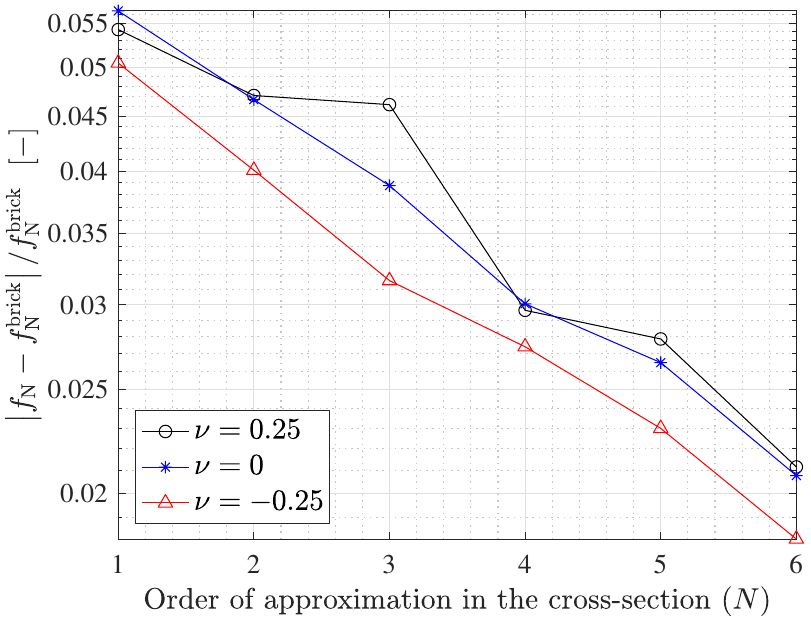}
	\caption{Lateral contact of a straight beam: convergence of the relative error in the total contact force between the beam $(f_\mathrm{N})$ and the brick $(f_\mathrm{N}^\mathrm{brick})$ element solutions for three different cases of Poisson's ratio: $\nu=0.25, 0, -0.25$. We use brick elements of $\mathrm{deg.}=(3,3,3)$, $n_\mathrm{el}=80\times10\times10$, and beam elements of $p=3$, and $n_\mathrm{el}=80$. Note that, in the results of $\nu=0.25$ and $\nu=-0.25$, we use the EAS method for $N=1$.}
	\label{ex_socket_conv_cforce_N}
\end{figure}

Next, we verify the contact formulation by comparing the pressure distribution between the presented beam contact formulation and the brick element solutions. Fig.\,\ref{ex1_socket_pres_surf_diff_pr_case} shows the comparison for the three different values of Poisson's ratios: $\nu=0.25$, $\nu=0$, and $\nu=-0.25$. For a positive Poisson's ratio, the bottom surface deforms into a concave shape, so that the contact pressure is concentrated around the two lateral edges, see Figs.\,\ref{ex1_socket_prs_surf_pr025_brick}-\ref{ex1_socket_prs_surf_pr025_beam_n2}. On the other hand, for a negative Poisson's ratio, the bottom surface deforms into a convex shape, so that the contact pressure is concentrated around the center of the bottom surface, see Figs.\,\ref{ex1_socket_prs_surf_pr-025_brick}-\ref{ex1_socket_prs_surf_pr-025_beam_N2}. For zero Poisson's ratio, the contact pressure is higher in the center region than that around the lateral edge in the results of the brick and beam ($N=2$) element formulations, see Figs.\,\ref{ex1_socket_prs_surf_pr0_brick}-\ref{ex1_socket_prs_surf_pr0_n2}.
\begin{figure*} \centering	
	\centering
	\begin{subfigure}[b] {0.325\textwidth} \centering
		\includegraphics[width=1\linewidth]{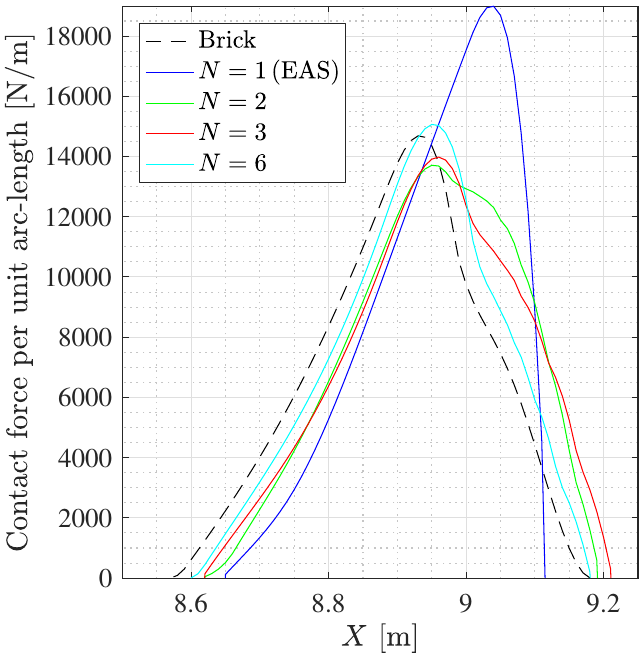}
		\caption{$\nu=0.25$}
		\label{ex_socket_cs_intg_prs_nu025}
	\end{subfigure}
	\begin{subfigure}[b] {0.325\textwidth} \centering
		\includegraphics[width=1\linewidth]{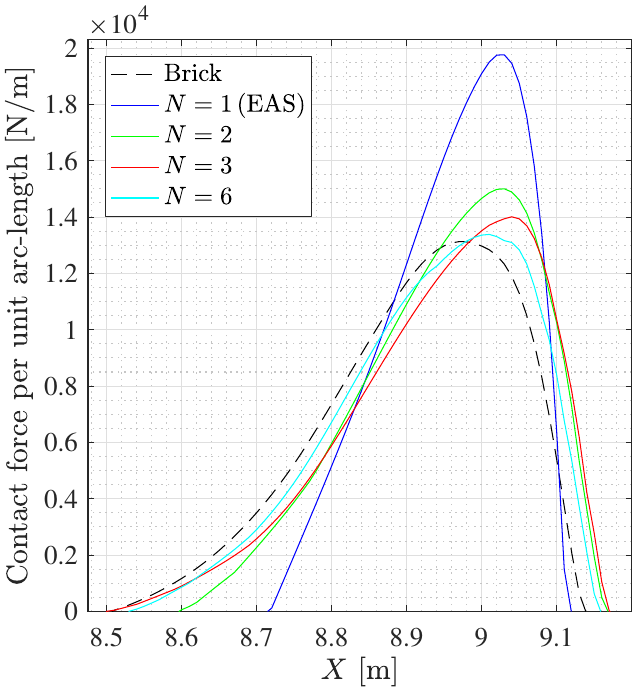}
		\caption{$\nu=-0.25$}
		\label{ex_socket_cs_intg_prs_nu_minus025}
	\end{subfigure}
	\begin{subfigure}[b] {0.325\textwidth} \centering
		\includegraphics[width=1\linewidth]{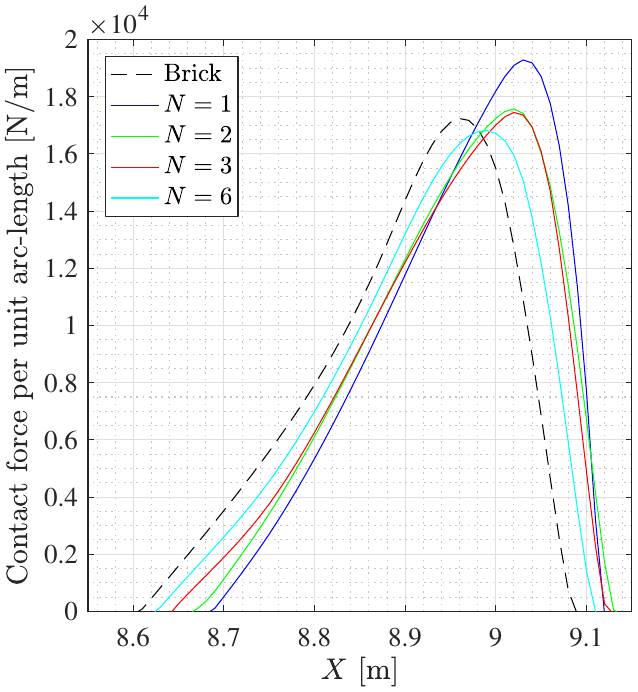}
		\caption{$\nu=0$}
		\label{ex_socket_cs_intg_prs_nu0}
	\end{subfigure}
	\caption{{Lateral contact of a straight beam: Comparison of the distribution of the contact force per unit initial arc-length (i.e., $r_\mathrm{N}$ from Eq.\,(\ref{def_cf_per_uinit_alen})) along the beam axis for the final deformed configuration (i.e., $\lambda_\mathrm{load}=1$) for three values of Poisson's ratio. We use brick elements with deg.=$(3,3,3)$, ${n_\mathrm{el}}=80\times10\times10$, and beam elements with $p=3$ and $n_\mathrm{el}=80$.}}
	\label{ex_socket_cs_cprs_cs_intg}
\end{figure*}
Further, we define the \textit{total contact force} by the integral of the contact pressure in the lateral surface, as
\begin{align}
	\label{tot_cf_intg}
		{f_{\rm{N}}} &\coloneqq  \int_{{\mathcal{R}_0}} {{p_{\rm{N}}}\,{\rm{d}}{\mathcal{R}_0}}=  \int_0^{{L^{(1)}}} {{r_{\rm{N}}}\,{\rm{d}}s}, 
\end{align}
where the \textit{contact force per unit initial arc-length} is obtained by
\begin{equation}
	\label{def_cf_per_uinit_alen}
	{r_{\rm{N}}} \coloneqq  \frac{1}{{{{\tilde j}^{(1)}}}}\int_{{\Xi ^{2(1)}}} {{p_{\rm{N}}}{{\tilde J}^{(1)}}\,{\rm{d}}{\xi ^2}}.
\end{equation}
Fig.\,\ref{ex_socket_cs_cprs_cs_intg} compares the distribution of $r_\mathrm{N}$ along the axis in the final deformed configuration from the beam solution with that of the brick element solution. In the beam solution for $N=1$, the cross-section boundary always remains straight (hence the constant contact pressure in transverse direction in Figs.\,\ref{ex1_socket_prs_surf_pr025_beam_n1eas}, \ref{ex1_socket_prs_surf_pr-025_beam_N1EAS}, and \ref{ex1_socket_prs_surf_pr0_n1eas}), which leads to an overestimation of the contact force. However, it is seen that the agreement improves, as $N$ increases. Fig.\,\ref{ex_socket_conv_cforce_N} also shows that the total contact force of Eq.\,(\ref{tot_cf_intg}) converges to the brick element solution, as $N$ increases.
\subsection{Lateral contact of a circular ring}
\label{contact_num_ex_ring_flat}
\begin{figure}[htb]
	\centering
	\begin{subfigure}[b]{0.295\textwidth} \centering
		\includegraphics[width=\linewidth]{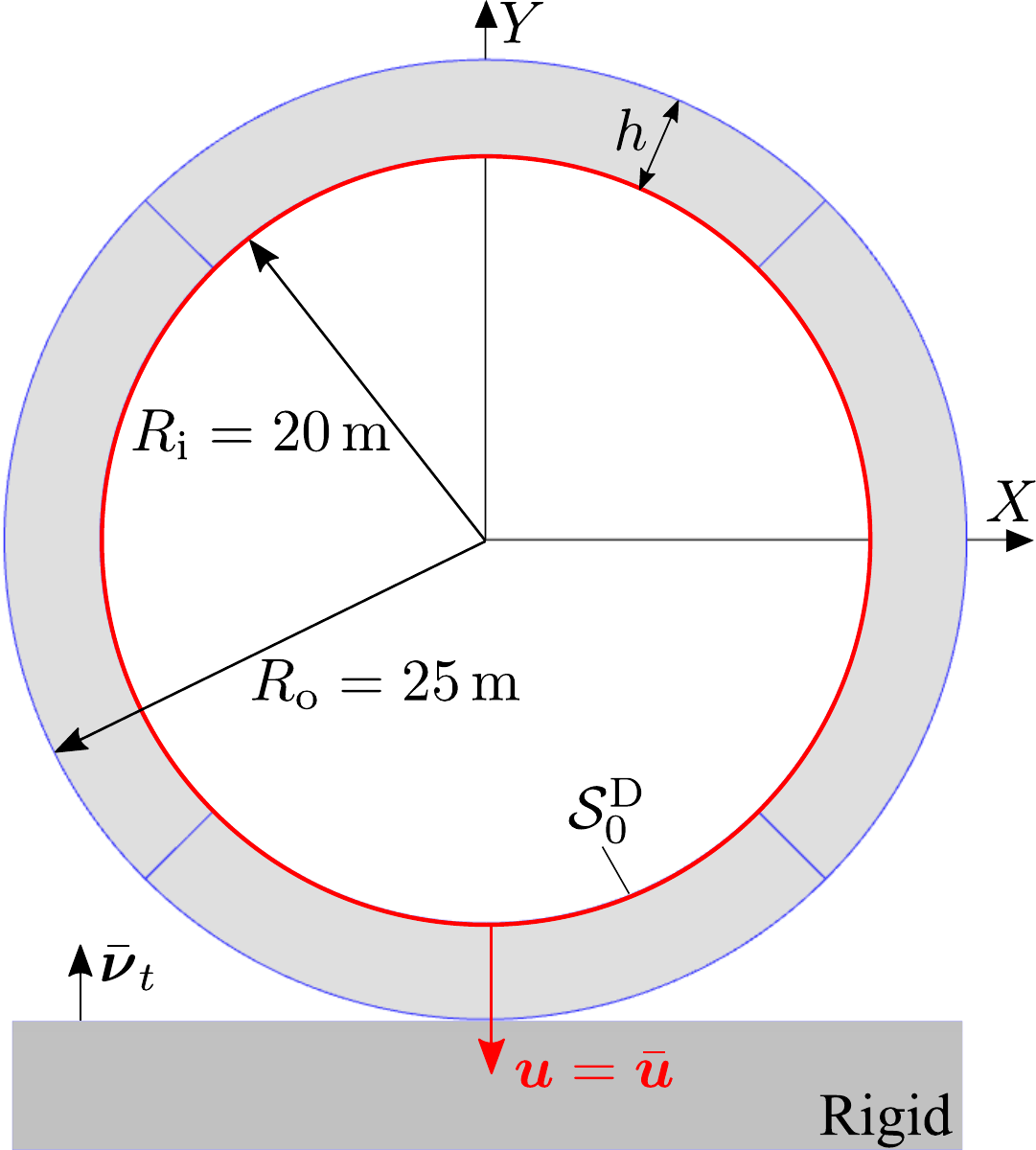}
		\caption{}
		\label{ex_ring_rigid_in_undeform_plane}
	\end{subfigure}
	\begin{subfigure}[b]{0.325\textwidth} \centering
		\includegraphics[width=\linewidth]{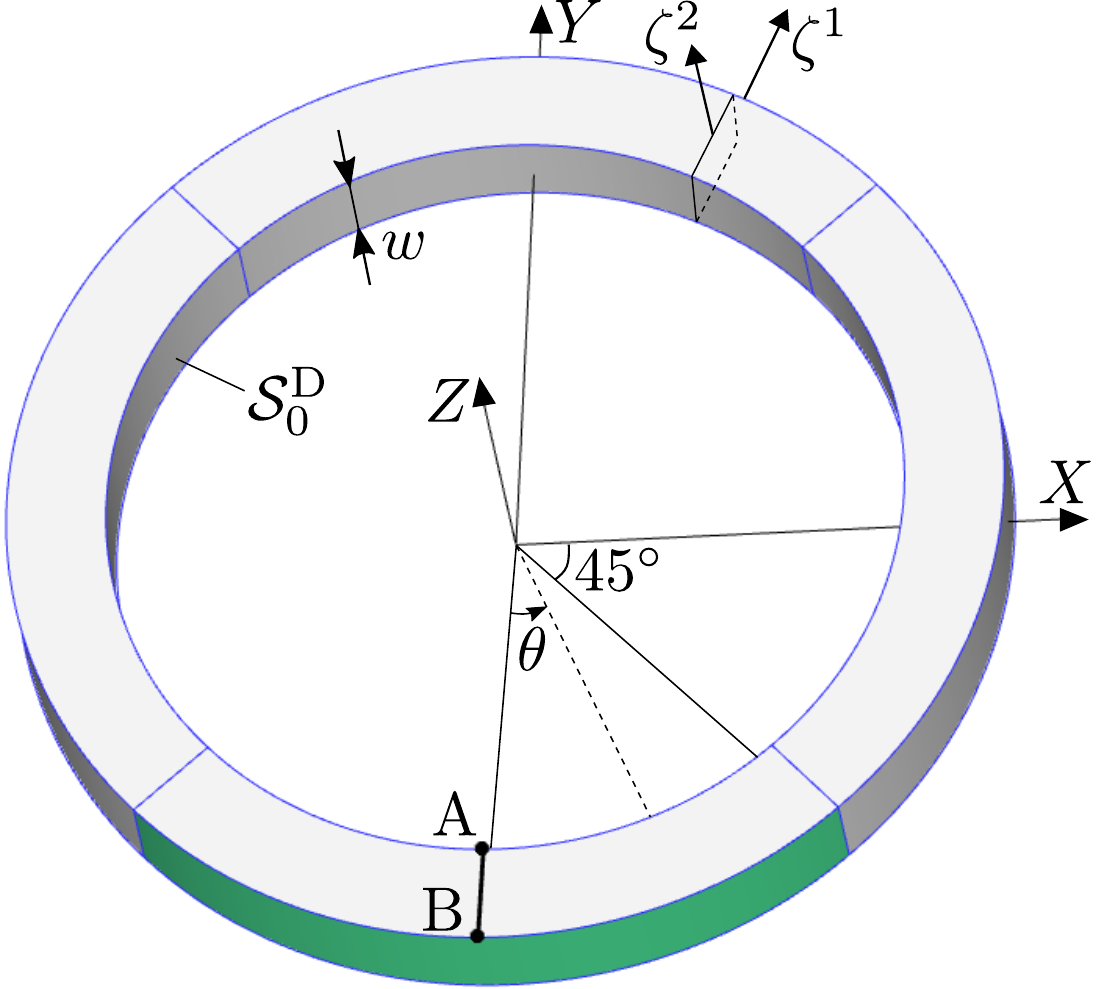}	
		\caption{}
		\label{ex_ring_rigid_in_undeform}
	\end{subfigure}
	\caption{Lateral contact of a circular ring: (a) undeformed configuration and prescribed displacement boundary condition at the inner surface $\mathcal{S}^\mathrm{D}_0$. The blue lines indicate the division of the domain into four parts due to modeling the ring exactly with four NURBS curves. (b) The two cross-sectional coordinates $\zeta^1$ and $\zeta^2$ are aligned with the radial, and $Z$-directions, respectively. In this problem, we consider contact only on the green-colored surface.}
	\label{ex1_cring_rflat_prob}
\end{figure}

\begin{figure} \centering
	\includegraphics[width=0.9999995\linewidth]{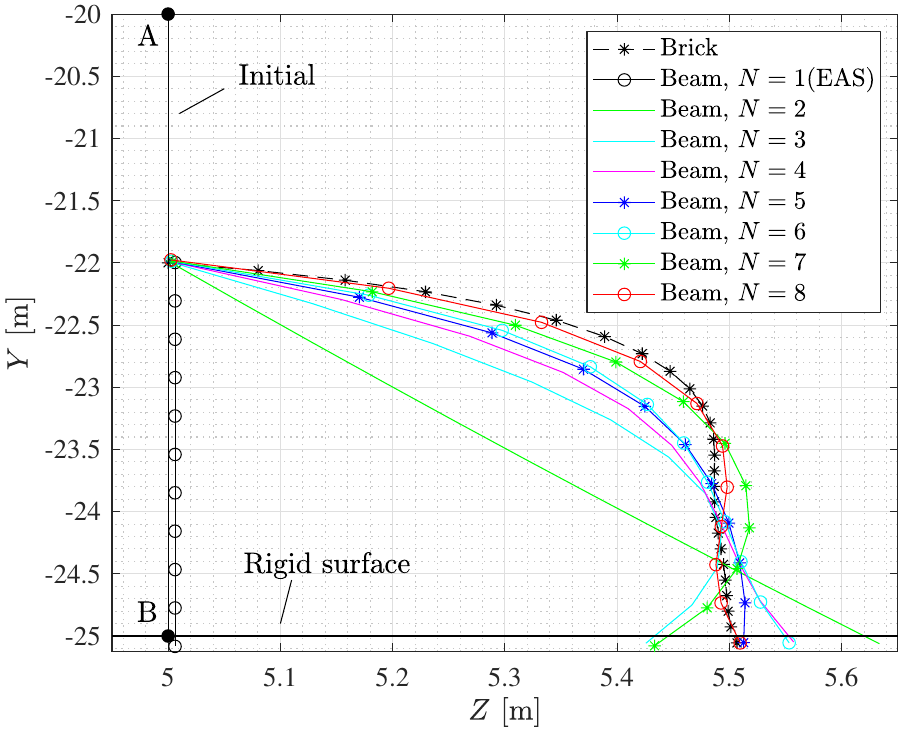}
	\caption{Lateral contact of a circular ring: Comparison of the deformation of the line $\overline{AB}$, defined in Fig.\,\ref{ex_ring_rigid_in_undeform}, for the beam and brick element solutions. Due to the chosen finite penalty parameters, the impenetrability constraints, and the displacement boundary condition on $\mathcal{S}^\mathrm{D}_0$ in the beam formulation are slightly violated. We use brick elements of deg.=$(4,4,4)$, $n_\mathrm{el}=320\times4\times4$, and beam elements of $p=3$, and $n_\mathrm{el}=320$.}
	\label{ex1_cring_rflat_latdisp_compare}
\end{figure}

\begin{figure*}
	\centering
	\begin{subfigure}[b] {0.3125\textwidth} \centering
		\includegraphics[width=\linewidth]{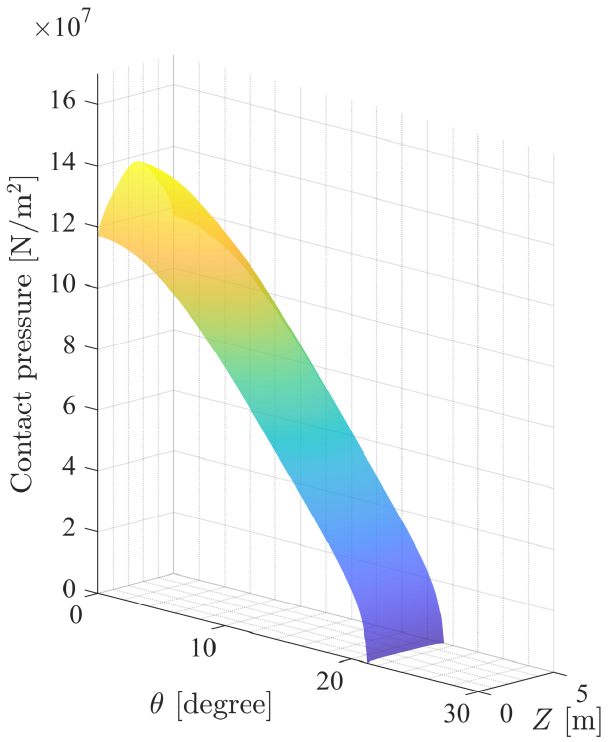}
		\caption{Brick}
	\end{subfigure}		
	\begin{subfigure}[b] {0.3125\textwidth} \centering				
		\includegraphics[width=\linewidth]{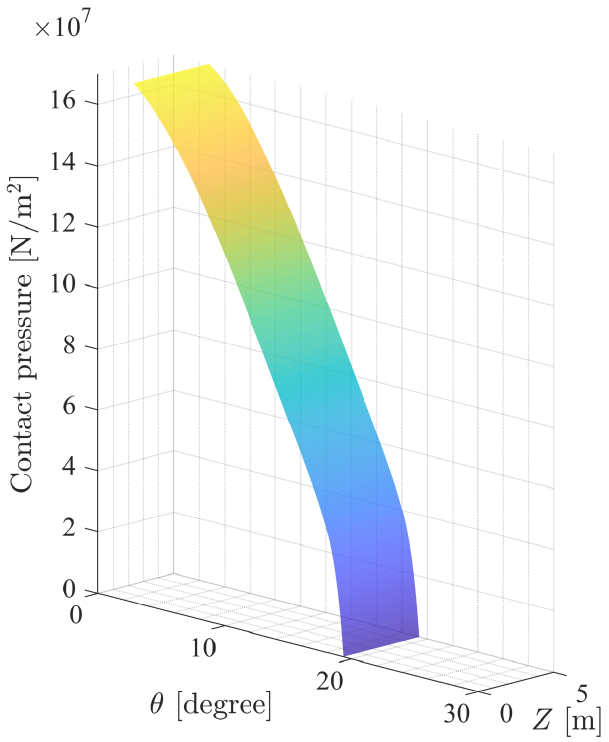}
		\caption{Beam, $N=1$ (EAS)}
	\end{subfigure}	
	\begin{subfigure}[b] {0.3125\textwidth} \centering				
		\includegraphics[width=\linewidth]{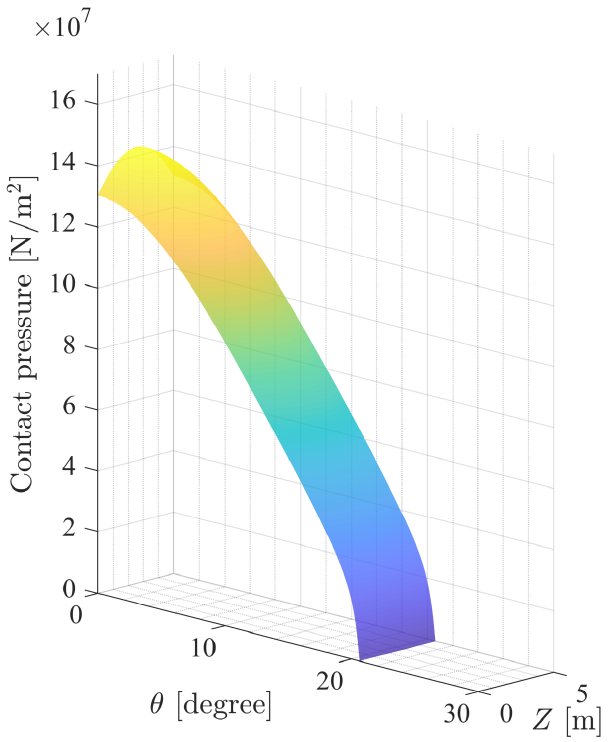}
			\caption{Beam, $N=2$}
	\end{subfigure}	
	\caption{Lateral contact of a circular ring: Comparison of contact pressure distributions on the bottom surface in beam and brick element models. (a) We use brick elements of $\mathrm{deg.}=(4,4,4)$, $n_\mathrm{el}=320\times4\times4$, (b,c) and beam elements of $p=3$, and $n_\mathrm{el}=320$.}
	\label{ex1_cring_rflat_prs_dist_surf}
\end{figure*}

\begin{figure}
	\centering
	\begin{subfigure}[b] {0.475\textwidth} \centering
		\includegraphics[width=\linewidth]{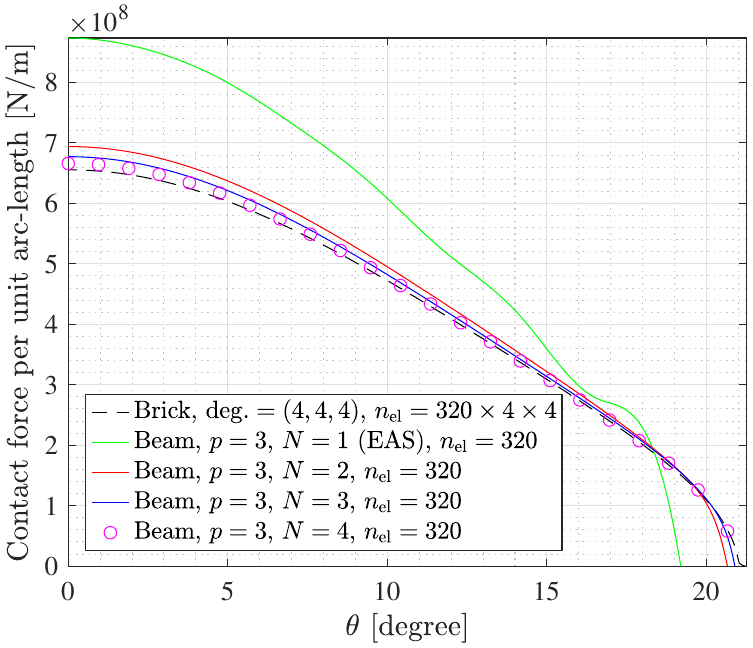}
	\end{subfigure}		
	\caption{Lateral contact of a circular ring: Comparison of the distribution of the contact force per unit initial arc-length, i.e., $r_\mathrm{N}$ from Eq.\,(\ref{def_cf_per_uinit_alen}), at the bottom surface in beam and brick element solutions. We use brick elements of $\mathrm{deg.}=(4,4,4)$, $n_\mathrm{el}=320\times4\times4$.}
	\label{ex1_cring_rflat_prs_d}
\end{figure}

\begin{figure}
	\centering
	\begin{subfigure}[b] {0.475\textwidth} \centering
		\includegraphics[width=\linewidth]{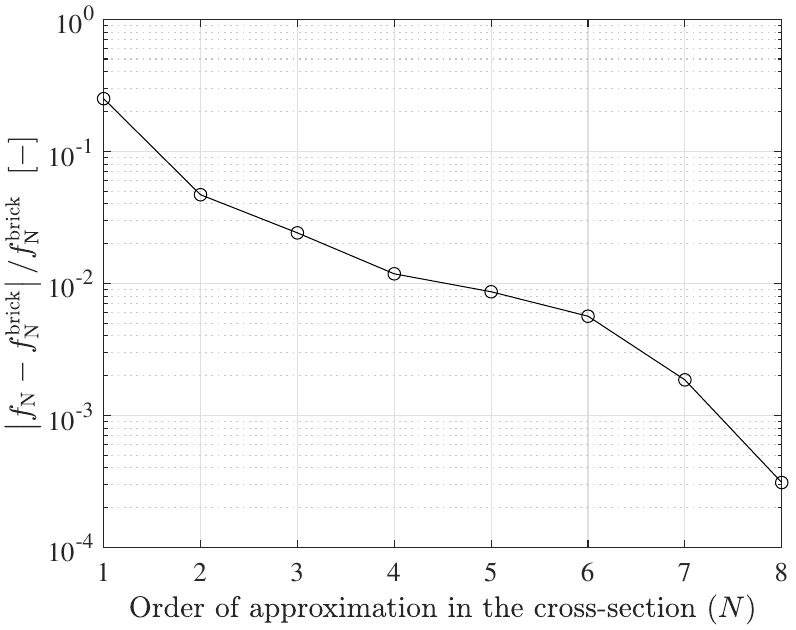}
	\end{subfigure}		
	\caption{Lateral contact of a circular ring: Convergence of the relative error in total contact force between the beam ($f_\mathrm{N}$) and brick ($f_\mathrm{N}^\mathrm{brick}$) element solutions, with increasing the order of approximation in the cross-section ($N$). We use brick elements of $\mathrm{deg.}=(4,4,4)$, $n_\mathrm{el}=320\times4\times4$, and beam elements of $p=3$, and $n_\mathrm{el}=320$.}
	\label{ex1_cring_rflat_prs_conv_rel_tot_f}
\end{figure}

In this example, we consider contact between an elastic circular ring and a rigid flat surface. The circular ring has a square cross-section with dimension $h\!=w\!=\!5\,\mathrm{m}$, and the inner and outer radii of the ring are $R_\mathrm{i}\!=\!20\,\mathrm{m}$ and $R_\mathrm{o}\!=\!25\,\mathrm{m}$, respectively, see Fig.\,\ref{ex_ring_rigid_in_undeform_plane}. We choose the Neo-Hookean material model with Young's modulus $210\,\mathrm{MPa}$ and Poisson's ratio $\nu=0.3$. Impenetrability is enforced by the penalty method using the penalty parameter ${\epsilon _{\rm{N}}}\!=\!{10^2}E/{L_0}$. The displacement is prescribed at the inner surface $\mathcal{S}^\mathrm{D}_0$ (see Fig.\,\ref{ex_ring_rigid_in_undeform_plane}), as
\begin{equation}
	\boldsymbol{u}\!\coloneqq\!\boldsymbol{x}-\boldsymbol{X}={\bar {\boldsymbol{u}}}\,\,\,\mathrm{on}\,\,\,\mathcal{S}^\mathrm{D}_0,
\end{equation}
where we choose the prescribed displacement vector $\bar{\boldsymbol{u}}=[-2\,\mathrm{m},0,0]^\mathrm{T}$, which means that the inner surface of the circular ring is rigidly translated in the negative $Y$-direction. For the beam formulation, the prescribed displacement boundary condition on the lateral surface $\mathcal{S}^\mathrm{D}_0$ can be enforced using a penalty method, and the detailed formulation can be found in Appendix \ref{app_contact_enforce_pre_trans}. The corresponding penalty parameter is chosen as ${\epsilon _{\rm{D}}} = {10^7}E/{L_0}$. For computational efficiency, we locally refine the mesh along the longitudinal direction in the domain quarter at the bottom (see Fig.\,\ref{ex1_cring_rflat_prob}). In the following, for brevity, we specify only the number of elements in the bottom parts for both beam and brick simulations. In the other parts, we use $n^\mathrm{L}_\mathrm{el}=10$ and $n_\mathrm{el}=10$ for brick and beam element solutions, respectively. In this example, we use uniform load increments with a total number of load steps $n_\mathrm{load}=5$, and $n_\mathrm{load}=40$ for beam and brick simulations, respectively. In both beam and brick simulations, we use $n_\mathrm{el}^\mathrm{sub}=10$ and $m_\mathrm{el}^\mathrm{sub}=20$.
Fig.\,\ref{ex1_cring_rflat_latdisp_compare} compares the displacement of the line $\overline{AB}$, defined in Fig.\,\ref{ex_ring_rigid_in_undeform}. This displacement contains a $Z$-component due to the Poisson effect. In case of $N=1$ (EAS), we have only one director along $\zeta^2$, see Fig.\,\ref{ex_ring_rigid_in_undeform}, and the homogeneous boundary condition on the $Z$-displacement at $\mathcal{S}^\mathrm{D}_0$ constrains the magnitude of the director $\boldsymbol{d}_2$. Therefore, the lateral displacement component $u_Z$ vanishes. As the order of approximation $N$ is increased, the lateral displacement of the beam formulation approaches that of the brick element solution, as Fig.\,\ref{ex1_cring_rflat_latdisp_compare} shows. Further we verify the contact pressure distribution at the bottom surface, i.e., the green-colored surface in Fig.\,\ref{ex_ring_rigid_in_undeform}. Due to the symmetry with respect to the $Y-Z$ plane, we plot the pressure only within the range ${0^ \circ } \le \theta  \le {30^ \circ }$. Fig.\,\ref{ex1_cring_rflat_prs_dist_surf} compares the contact pressure distributions of the beam and brick element models. Fig.\,\ref{ex1_cring_rflat_prs_d} compares the distribution of the contact force per unit initial arc-length, i.e., $r_\mathrm{N}$ of Eq.\,(\ref{def_cf_per_uinit_alen}). In case of $N=1$, the contact pressure is constant in transverse direction due to the first order approximation of the displacement field, which leads to an overestimation of the contact force. As the order of approximation $N$ increases, the agreement between the contact pressure distribution and the brick element solution improves significantly, as Figs.\,\ref{ex1_cring_rflat_prs_dist_surf} and \ref{ex1_cring_rflat_prs_d} show. Further, Fig.\,\ref{ex1_cring_rflat_prs_conv_rel_tot_f} shows that the total contact force of Eq.\,(\ref{tot_cf_intg}) converges to the corresponding brick element solution with increasing order $N$.
\subsection{Sliding contact between two initially straight beams}
\label{num_ex_slide_2beam}
\begin{figure}
	\centering
	\begin{subfigure}[b] {0.475\textwidth} \centering
		\includegraphics[width=\linewidth]{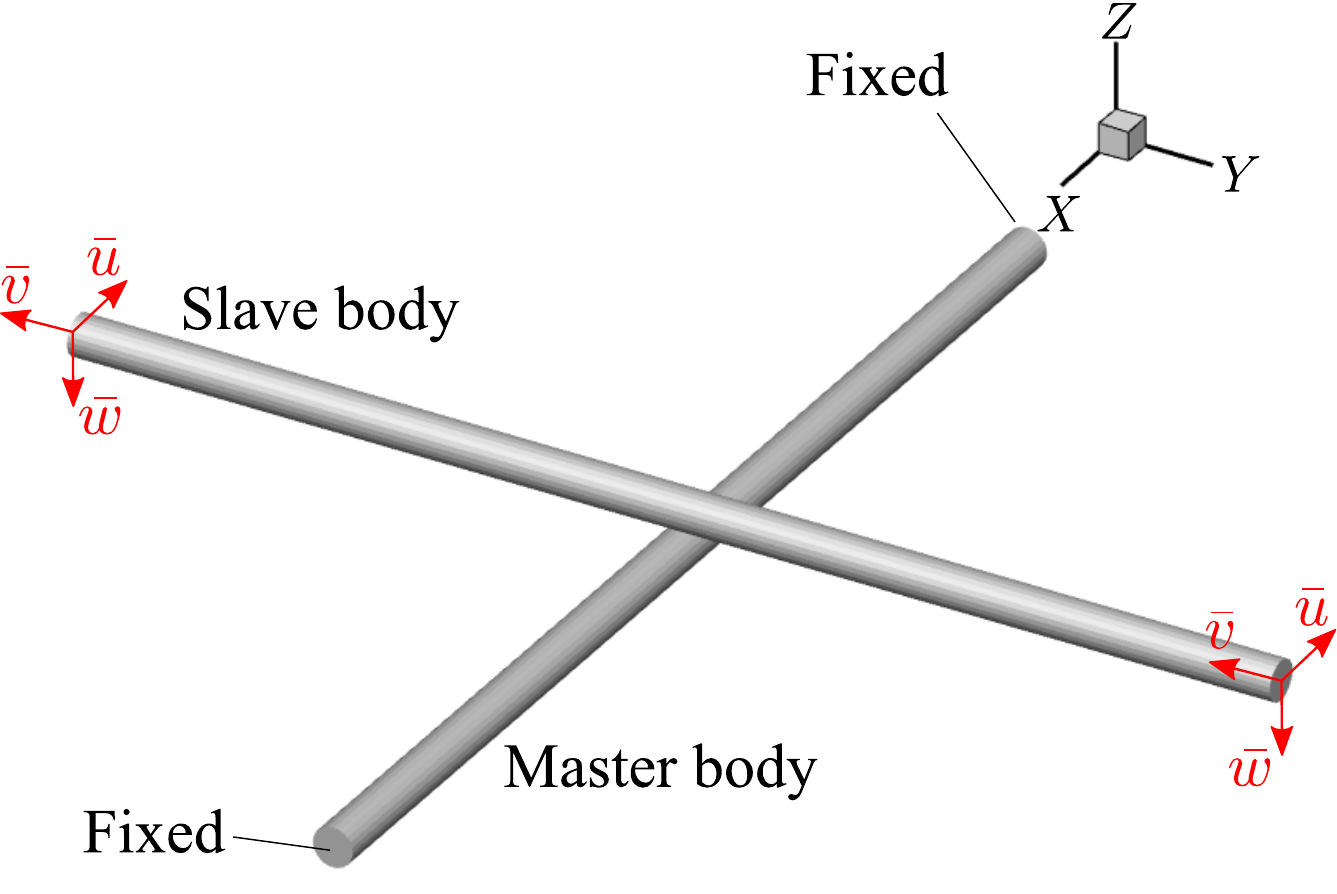}
		\caption{A perspective view}
		\label{ex_str2b_prob_desc_pers}
	\end{subfigure}
	\begin{subfigure}[b] {0.475\textwidth} \centering
	\includegraphics[width=\linewidth]{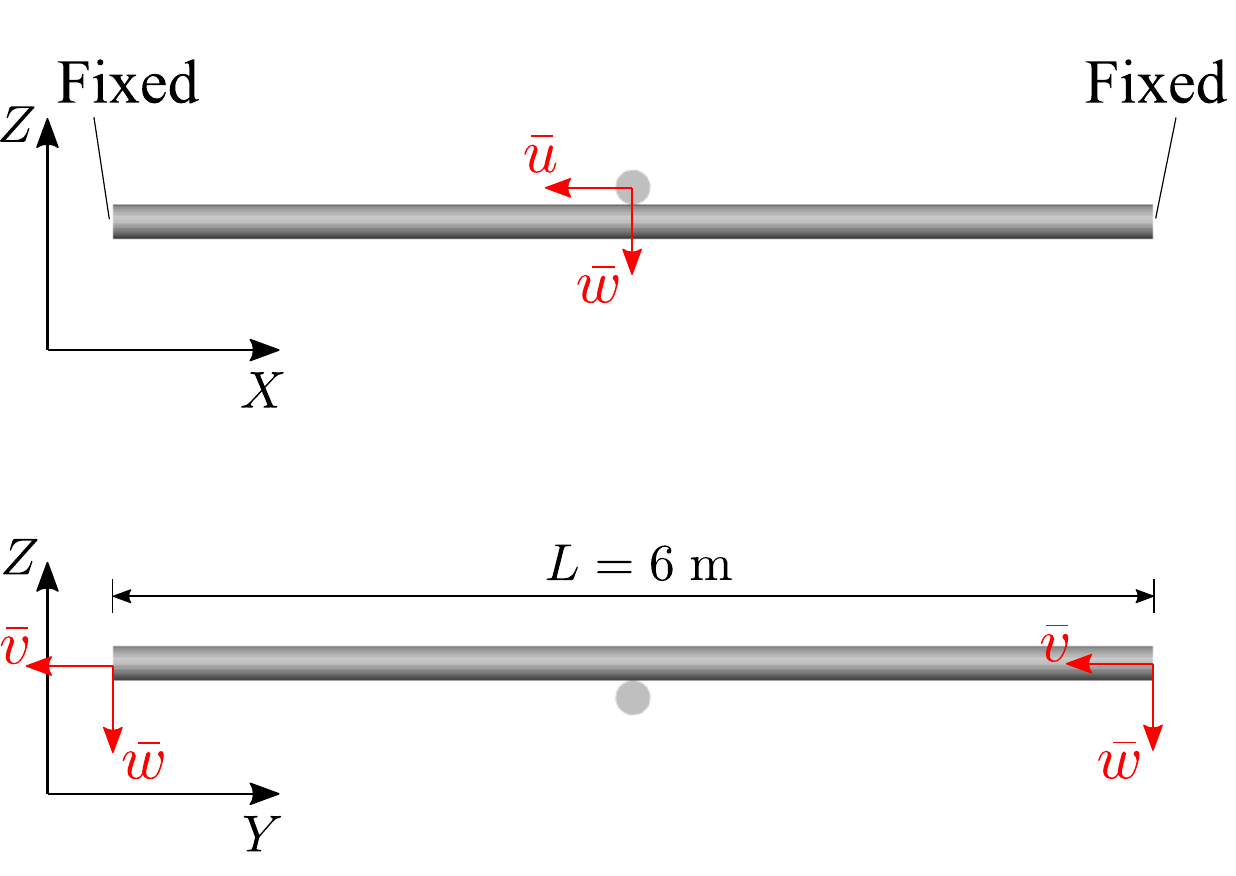}
	\caption{Planar views}
	\label{ex_str2b_prob_desc_planar}
	\end{subfigure}
	\caption{Sliding contact between two initially straight beams: Initial configuration and boundary conditions. The prescribed displacements in $X$-, $Y$-, and $Z$-directions at the ends of the upper beam are denoted by $\bar u$, $\bar v$, and $\bar w$, respectively, and chosen as ${\bar u}={\bar v}={\bar w}=-1.5\,\mathrm{m}$.}
	\label{ex_str2b_prob_desc}	
\end{figure}

We consider sliding contact between two initially straight beams. The two beams have the same length $L=6\,\mathrm{m}$ and a circular cross-section of radius $R=0.1\,\mathrm{m}$, and they are initially perpendicular to each other with vertical distance $d=10^{-3}\,\mathrm{m}$. We select the Neo-Hookean material model with Young's modulus $E=210\times10^9\,\mathrm{Pa}$ and Poisson's ratio $\nu=0.3$. In the contact formulation, we choose the lower beam as the master body, and the upper beam as the slave body. The lower beam is fixed at both ends, and the upper beam is under non-homogeneous displacement boundary conditions at both ends of the beam, see Fig.\,\ref{ex_str2b_prob_desc}. We investigate the following two cases of deformability.
\begin{itemize}
	\item Case 1: the slave body (upper beam) is rigid, and the master body (lower beam) is deformable,
	\item Case 2: both beams are deformable. It should be noted that the prescribed displacement boundary conditions apply to the whole end faces. Thus, the cross-sections at both ends are not deformable.
\end{itemize}
\noindent We choose the cutoff radius $r_\mathrm{c}=3R=0.3\,\mathrm{m}$, and $\varepsilon_\theta=0.4\pi$ in the global contact search.
\subsubsection{Case 1: Rigid slave beam}
\label{ex_btb_slide_cs1}
\begin{figure}
	\centering
	\begin{subfigure}[b] {0.465\textwidth} \centering
		\includegraphics[width=\linewidth]{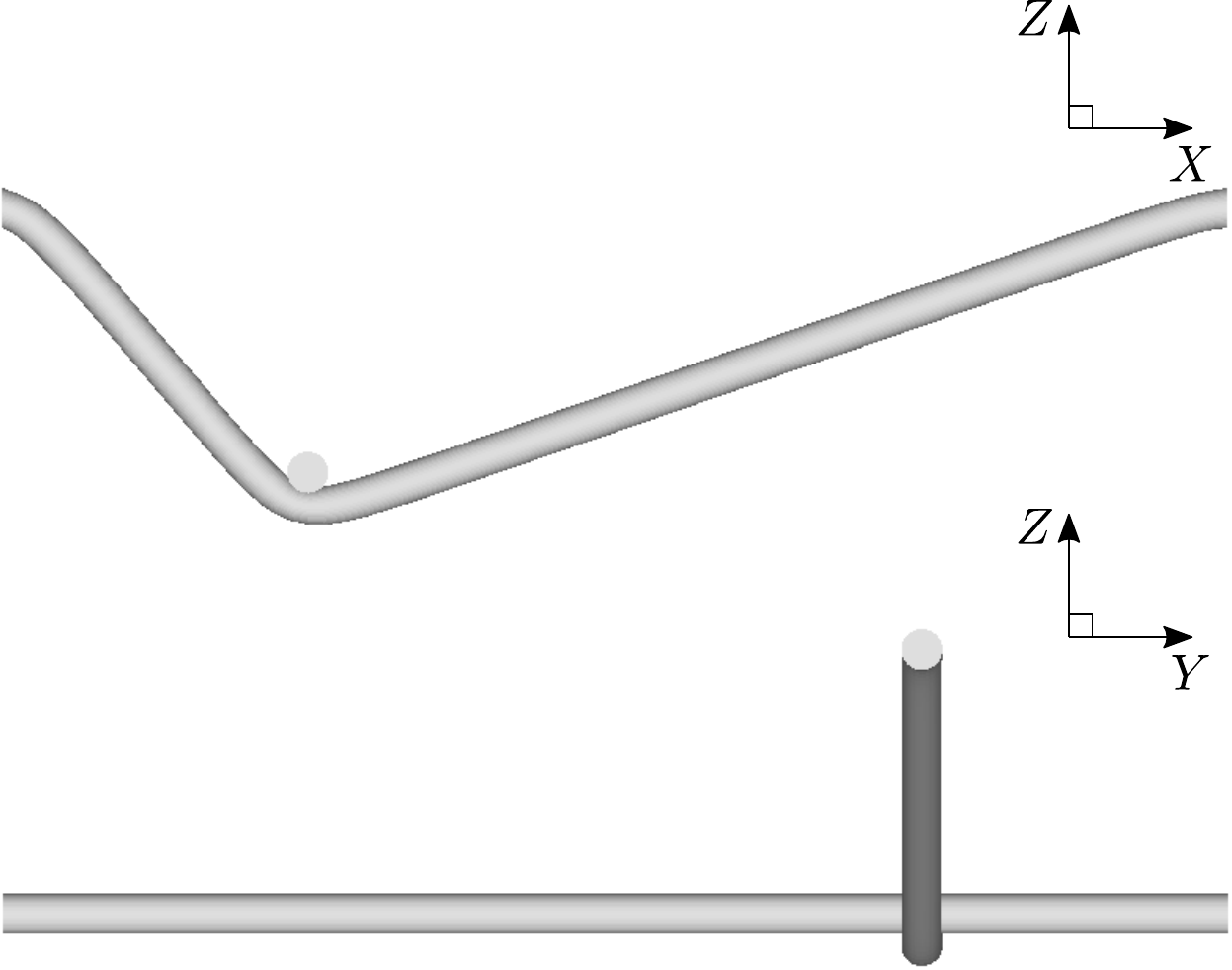}
	\end{subfigure}
	\caption{Sliding contact between two initially straight beams (case 1): Final deformed configuration in the case of rigid slave (upper) beam. For the discretization of the beam axis, we use B-spline basis functions with $p=3$, and $n_\mathrm{el}=80$ and $n_\mathrm{el}=160$ for the slave and master bodies, respectively. The chosen penalty parameter is $\epsilon_\mathrm{N}=20E/{L_0}$.}
	\label{ex_str2b_case1_deform}	
\end{figure}

\begin{figure}
	\centering
	\begin{subfigure}[b] {0.475\textwidth} \centering
		\includegraphics[width=\linewidth]{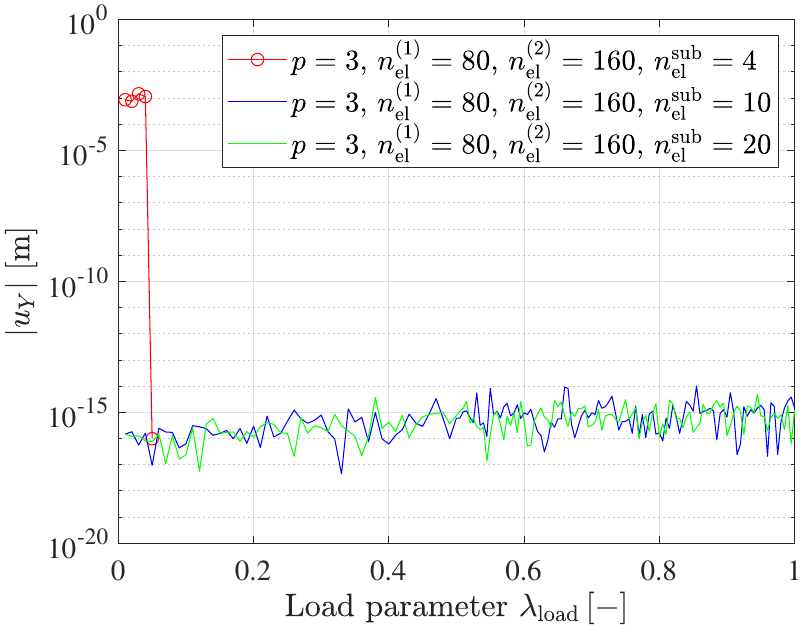}
	\end{subfigure}
	\caption{Contact between two initially straight beams (case 1): The magnitude of $Y$-displacement ($\left\lvert{u_Y}\right\rvert$) at the center of the beam's axis ($s=L/2$) in the master body (lower beam). The solution using $n^\mathrm{sub}_\mathrm{el}=4$ diverges at the 6th load step. In cases of $n^\mathrm{sub}_\mathrm{el}=10$ and $n^\mathrm{sub}_\mathrm{el}=20$, the $Y$-displacement vanishes to machine precision. The chosen penalty parameter is $\epsilon_\mathrm{N}=10E/{L_0}$. In all cases, $m_\mathrm{el}^\mathrm{sub}=100$.}
	\label{cont_slide_case1_cent_ydisp_master}	
\end{figure}

\begin{figure}
	\centering
	\begin{subfigure}[b] {0.475\textwidth} \centering
		\includegraphics[width=\linewidth]{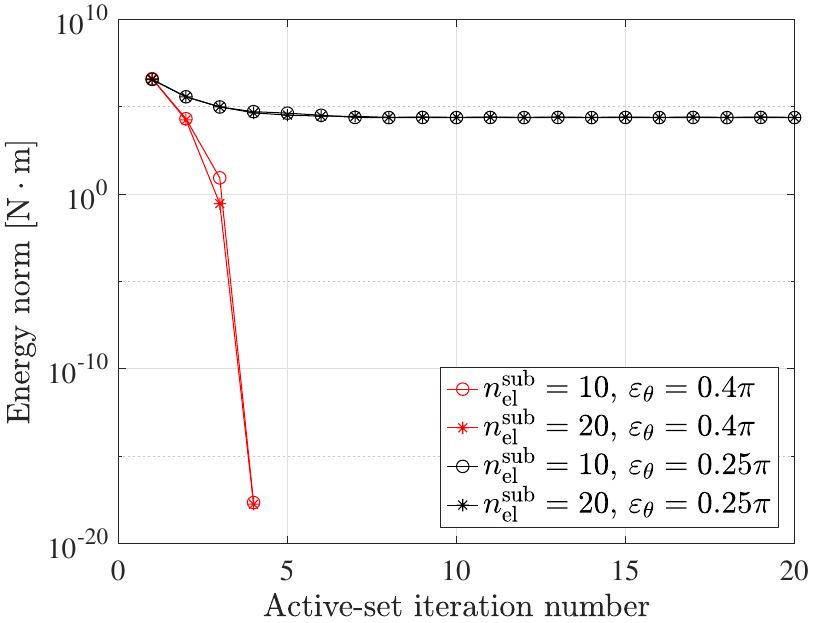}
		\caption{History of energy norm}
		\label{ex_str2b_case1_hist_enorm}		
	\end{subfigure}
	\begin{subfigure}[b] {0.475\textwidth} \centering
		\includegraphics[width=\linewidth]{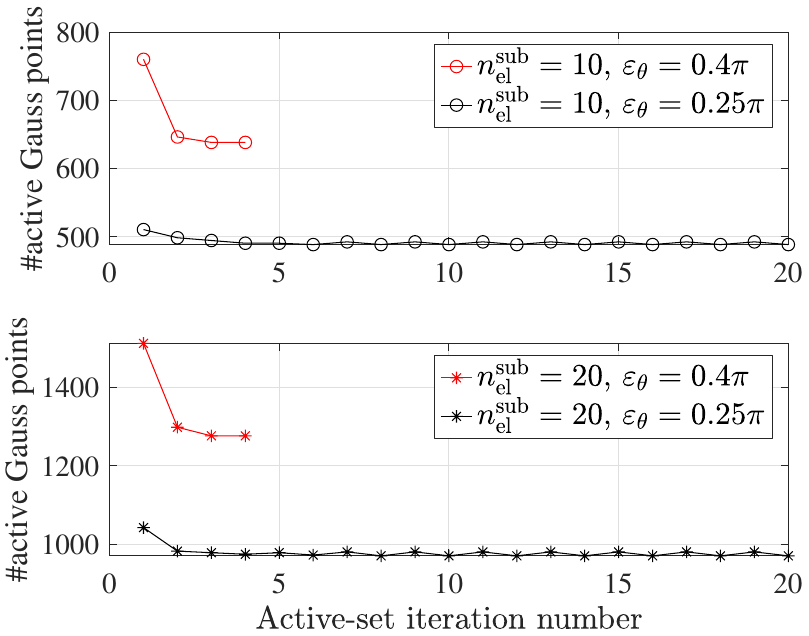}
		\caption{History of the number of active Gauss points}
		\label{ex_str2b_case1_hist_ngp}
	\end{subfigure}
	\caption{Sliding contact between two initially straight beams (case 1): Convergence history of the energy norm and the total number of active surface Gauss integration points in the slave body during the active set iteration at the last load step for two different numbers of Gauss integration points along the axis of slave body in each case of $\varepsilon_\theta$. The chosen penalty parameter is $\epsilon_\mathrm{N}=10E/{L_0}$.}
	\label{ex_str2b_case1_hist}	
\end{figure}

\begin{figure}
	\centering
	\begin{subfigure}[b] {0.475\textwidth} \centering
		\includegraphics[width=\linewidth]{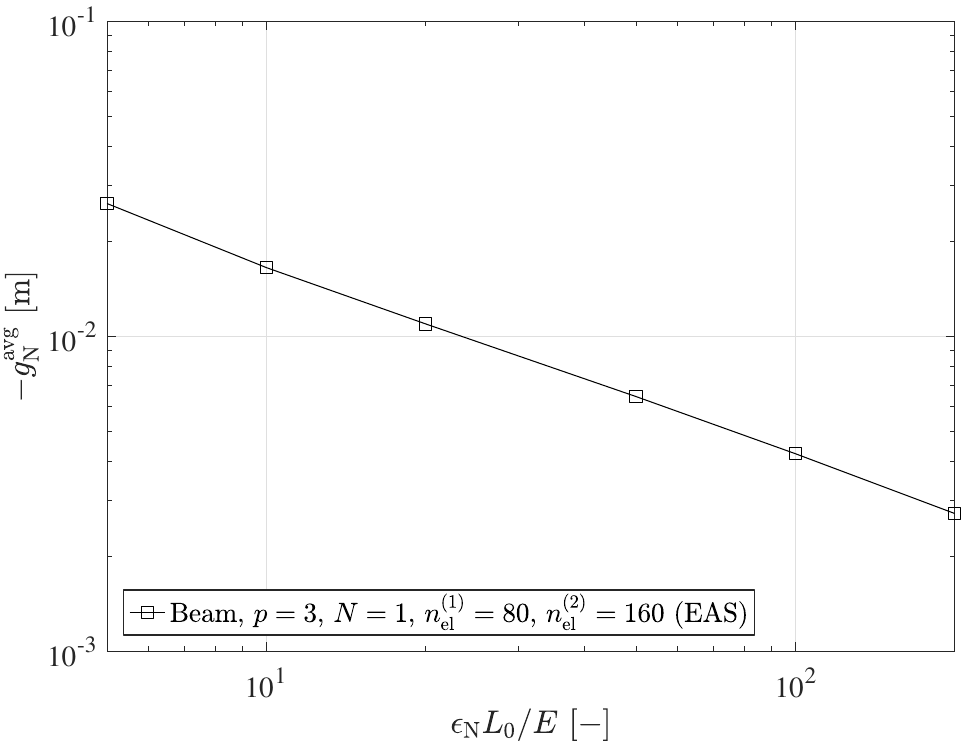}
		\caption{Average normal gap}
		\label{btb_slide_case1_conv_penet}	
	\end{subfigure}
	\begin{subfigure}[b] {0.475\textwidth} \centering
	\includegraphics[width=\linewidth]{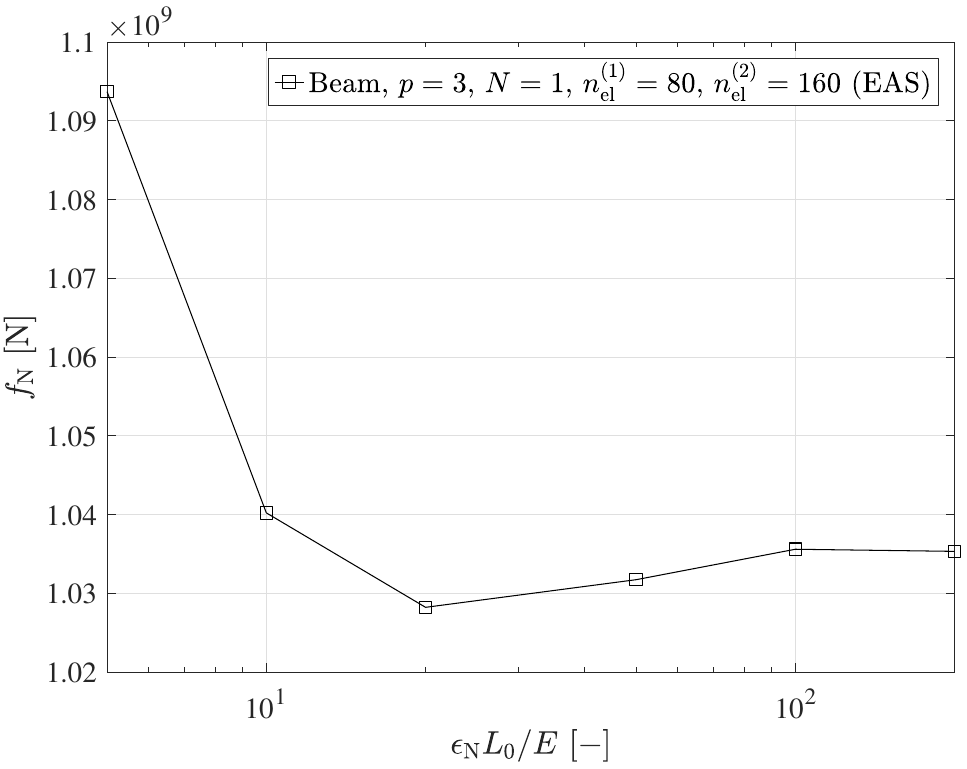}
	\caption{Total contact force}
	\label{btb_slide_case1_conv_cforce}	
	\end{subfigure}
	\caption{Contact between two initially straight beams (case 1): Convergence of (a) the average normal gap from Eq.\,(\ref{pent_per_unit_undef_a}) (b) and the total contact force from Eq.\,(\ref{tot_cf_intg}) versus the penalty parameter $\epsilon_\mathrm{N}$. The dashed line represents a linear rate of convergence.}
	\label{btb_slide_case1_conv_penet_cforce}
\end{figure}

The first case considers a rigid slave body, and we implement the \textit{rigidity} by constraining all DOFs in the finite element discretization of the beam, and the closest point projection employs the numerical scheme presented in Sections \ref{glob_csearch} and \ref{loc_csearch}. Fig.\,\ref{ex_str2b_case1_deform} shows the deformed configuration. During the sliding contact, as we consider frictionless contact, no $Y$-directional contact force should act on the lower beam. However, if the number of Gauss integration points along the axis of the slave body is not sufficient like the case $n^\mathrm{sub}_\mathrm{el}=4$, it is seen in Fig.\,\ref{cont_slide_case1_cent_ydisp_master} that an unphysical $Y$-displacement occurs. For both cases $n^\mathrm{sub}_\mathrm{el}=10$ and $20$, it is shown that the $Y$-displacement vanishes to machine precision. However, it is observed that, if $\varepsilon_\theta$ is too small (e.g., $\varepsilon_\theta=0.25\pi$), the active-set iteration does not converge but oscillates between two different contact states in several load steps, e.g., the last one, see Fig.\,\ref{ex_str2b_case1_hist}. With $\varepsilon_\theta=0.4\pi$, the active set iteration converges in all load steps. We calculate the \textit{average normal gap} over the whole contact area by
\begin{align}
	\label{pent_per_unit_undef_a}
	{g^\mathrm{avg}_{\rm{N}}} &\coloneqq  \frac{1}{\int_{{\mathcal{R}_0}} {{\rm{d}}{\mathcal{R}_0}}}\int_{{\mathcal{R}_0}} {{g_{\rm{N}}}\,{\rm{d}}{\mathcal{R}_0}}.
\end{align}
Fig.\,\ref{btb_slide_case1_conv_penet} shows that the average normal gap decreases, and the total contact force converges, with increasing penalty parameter. Table\,\ref{app_2beam_rigid_flex_ngp_info} shows the selected number of sub-elements in axial and circumferential directions for the contact integral, and Table\,\ref{app_2beam_rigid_flex_nload} shows the selected load increment sizes for each case of the penalty parameter. It is typically required to increase the number of surface Gauss integration points for the contact integral and reduce the load increment size in order to achieve convergence in the solution process using a larger penalty parameter. Too many Gauss integration points or high penalty parameters may lead to an overconstrained system, a case also called \textit{contact locking}, especially for low order finite elements. In this paper, we basically use higher order basis functions, so that we do not observe such locking. Further steps to alleviate this locking, e.g., using a mortar-type discretization, remains future work. One can also develop a scheme to deactivate superfluous Gauss integration points in the active set (outer) loop, or adaptively control the penalty parameter.
\subsubsection{Case 2: Two deformable beams}
\label{num_ex_sliding_2b_c2_def}
\begin{figure}
	\centering
	\begin{subfigure}[b] {0.465\textwidth} \centering
		\includegraphics[width=\linewidth]{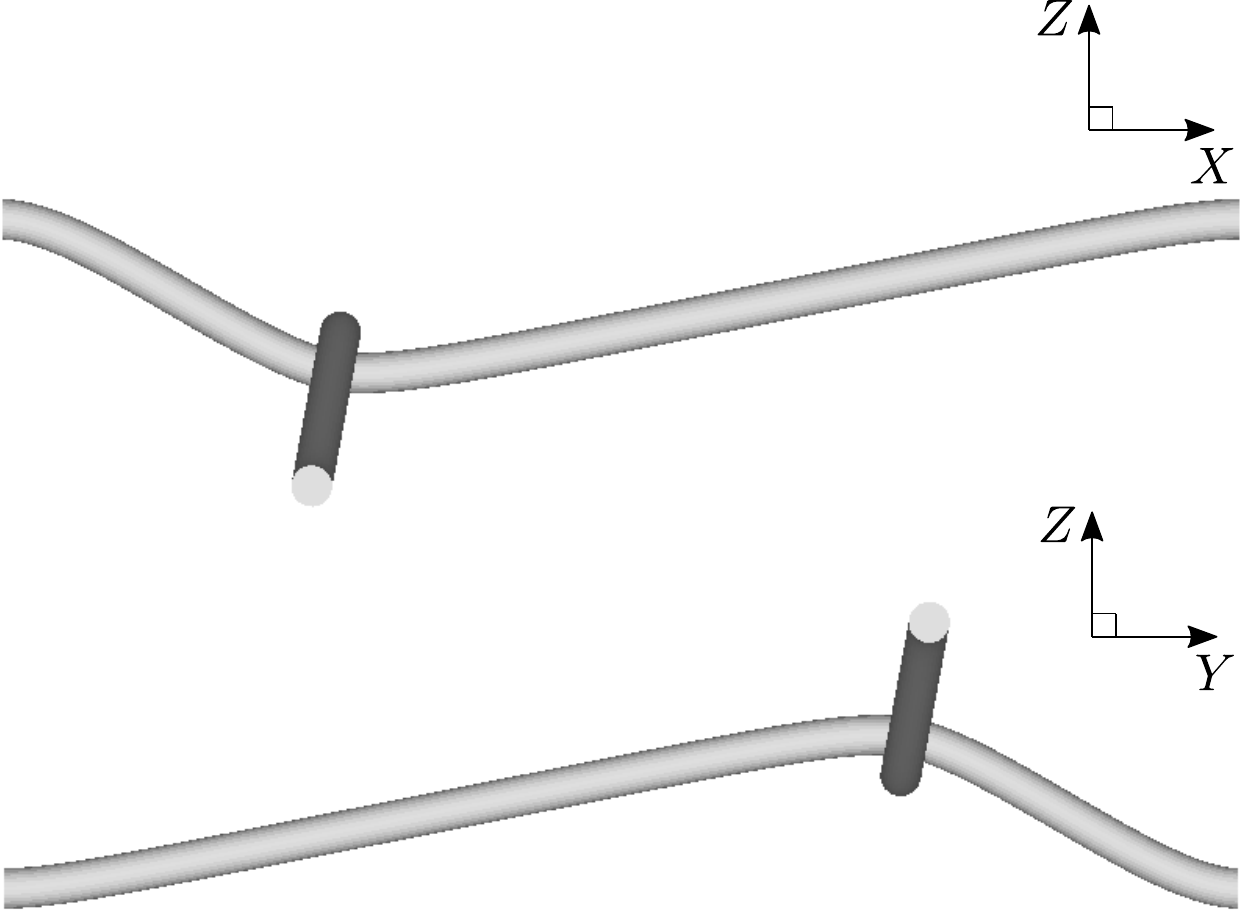}
	\end{subfigure}	
	\caption{Sliding contact between two initially straight beams (case 2): For the discretization of the beam axis, we use B-spline basis functions with $p=3$, and $n_\mathrm{el}=160$ for both the slave and master bodies. The chosen penalty parameter is $\epsilon_\mathrm{N}={10^2}E/{L_0}$.}
	\label{ex_str2b_case2_deform}
\end{figure}

\begin{figure}
	\centering
	\begin{subfigure}[b] {0.475\textwidth} \centering
		\includegraphics[width=\linewidth]{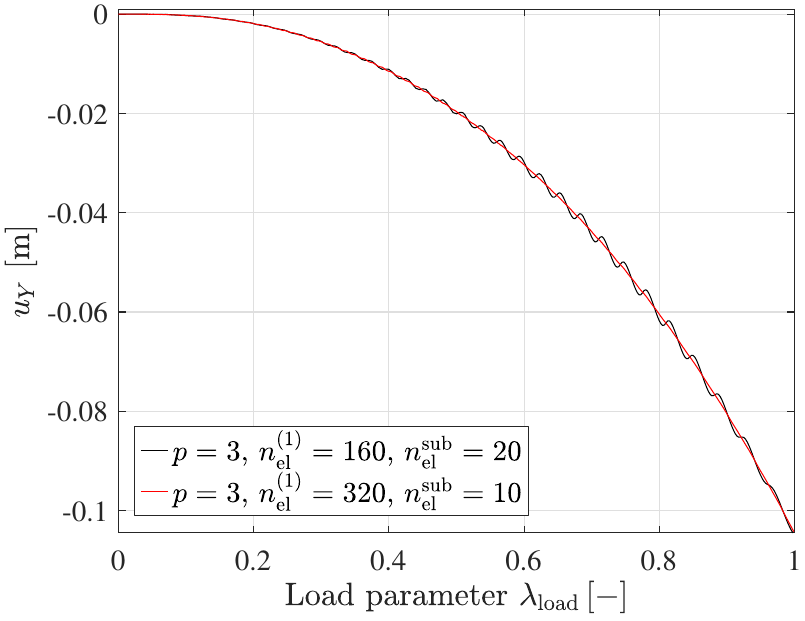}
	\end{subfigure}\\
	\caption{Contact between two initially straight beams (case 2): $Y$-displacement at the center of axis ($s=L/2$) in the master beam for two different discretization. It should be noted that those two cases have the same total number of Gauss integration points for the evaluation of the contact integral on the slave body. In all cases, we use $m^\mathrm{sub}_\mathrm{el}=100$.}
	\label{cont_slide_cent_ydisp_master}	
\end{figure}

\begin{figure}
	\centering
	\begin{subfigure}[b] {0.475\textwidth} \centering
		\includegraphics[width=\linewidth]{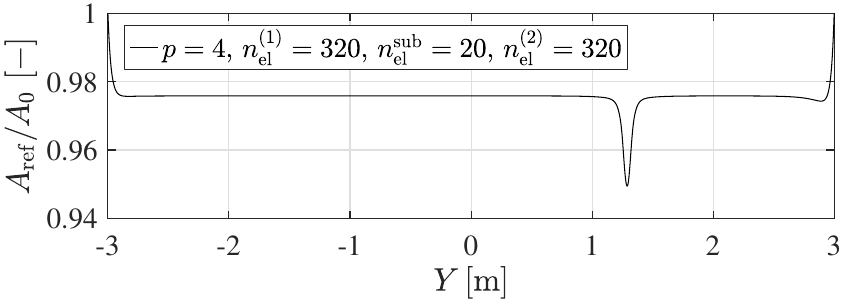}
		\caption{Ratio between the cross-section area of the reference solution and the initial solution}
		\label{cont_slide_carea_ref_sol}			
	\end{subfigure}
	\begin{subfigure}[b] {0.475\textwidth} \centering
		\includegraphics[width=\linewidth]{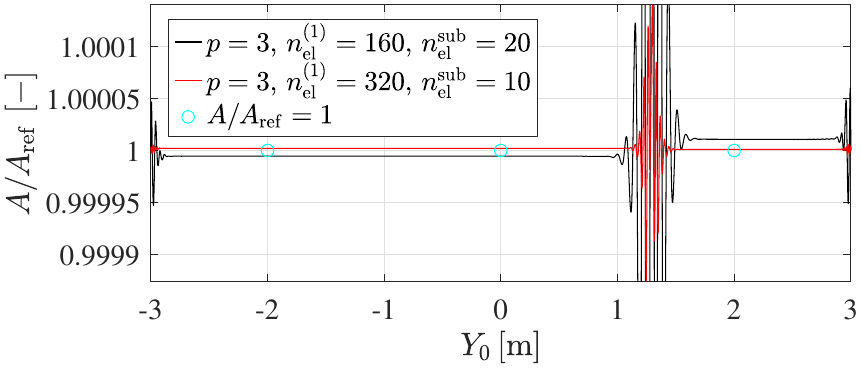}
		\caption{Ratio between the cross-section area of the beam solution and the reference solution}
		\label{cont_slide_carea_ratio_mag_view}		
	\end{subfigure}
	\caption{Contact between two initially straight beams (case 2): Distribution of the cross-sectional area along the axis for (a) the reference solution and (b) the beam element solution. In (b), the black curve is cut off by the minimum and maximum values of the red curve in the vertical axis for a clear visualization. In all cases, we use $m_\mathrm{el}^\mathrm{sub}=100$. The original graph can be found in Fig.\,\ref{cont_slide_carea_org_data}.}
	\label{cont_slide_flex_carea_ratio}	
\end{figure}

Next, we consider both beams deformable. Fig.\,\ref{ex_str2b_case2_deform} shows the final deformed configuration. In Fig.\,\ref{cont_slide_cent_ydisp_master}, we compare the lateral ($Y$-directional) displacement at the center of the lower beam's axis during the deformation for two different numbers of DOFs in the upper beam (slave body): $n^{(1)}_\mathrm{el}=160$ and $n^{(1)}_\mathrm{el}=320$ with the number of sub-elements $n^\mathrm{sub}_\mathrm{el}=20$, and $n^\mathrm{sub}_\mathrm{el}=10$ for the evaluation of the contact integral, respectively, such that both cases have the same total number of Gauss integration points. In the former case with less DOFs in the slave body, the lateral displacement exhibits spurious oscillations. This is associated with an oscillatory change of cross-sectional area along the axis in the deformed configuration of the slave body due to \textit{curvature thickness locking}, which means an artificial coupling between the cross-sectional stretching and the bending deformation. In Fig.\,\ref{cont_slide_flex_carea_ratio}, we compare the distribution of the cross-sectional area along the axis at the final deformed configuration for the two cases of the number of DOFs in the slave body. In Fig.\,\ref{cont_slide_flex_carea_ratio}, the reference solution of the cross-sectional area ($A_\mathrm{ref}$) is obtained by using B-spline basis functions with $p=4$, $n^{(1)}_\mathrm{el}=n^{(2)}_\mathrm{el}=320$, $n^\mathrm{sub}_\mathrm{el}=20$, and $m^\mathrm{sub}_\mathrm{el}=100$. Fig.\,\ref{cont_slide_flex_carea_ratio} shows that the cross-sectional area decreases in the whole domain due to the Poisson effect caused by the axial stretching. It leads to large curvature in the lateral surface around the loaded area as well as the fixed boundary, which eventually leads to the oscillatory lateral displacement during the sliding contact. As shown in \citet{choi2021isogeometric}, this locking can be alleviated by mesh refinement, see Fig.\,\ref{cont_slide_carea_ratio_mag_view}. Thus, in Fig.\,\ref{cont_slide_cent_ydisp_master}, it is seen that as we increase the number of DOFs in the slave body while maintaining the same total number of Gauss integration points for the evaluation of the contact integral, the amplitude of oscillation significantly decreases. The selected load increment sizes for each case of the results in Figs.\,\ref{cont_slide_cent_ydisp_master} and \ref{cont_slide_flex_carea_ratio} can be found in Table \ref{app_2beam_flex_flex_nload}.
\subsection{Twisting of a wire strand}
\label{num_ex_twist}
\begin{figure}	
	\centering
	\begin{subfigure}[b] {0.4875\textwidth} \centering
		\includegraphics[width=\linewidth]{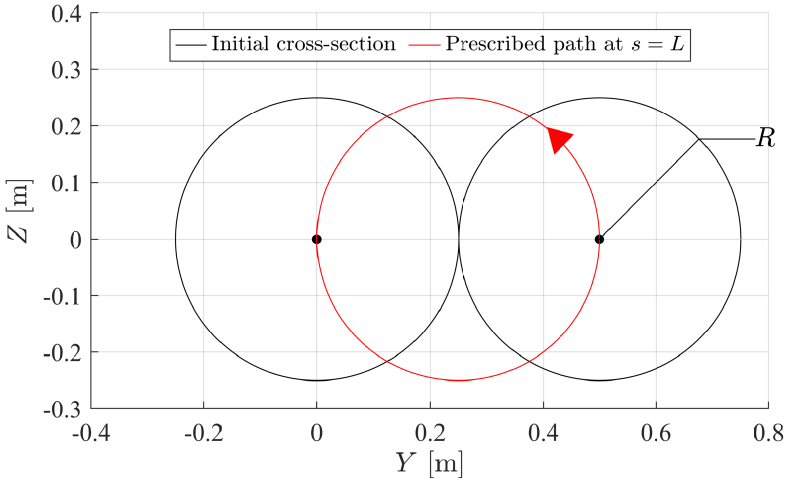}
	\end{subfigure}
	\caption{Twisting of a wire strand (two wires): Initial configuration of the wires in $Y$-$Z$ plane (cross-sectional) view, and the same circular path of radius $R$ is prescribed to the ends of the axes ($s=L$) of both beams. This is shown by the red curve with the arrow indicating the direction of prescribed rotation. Note that we choose the center of the prescribed circular path $c_2=R$ and $c_3 = 0$. We choose the left and right beams as slave and master bodies, respectively.}
	\label{twist_2b_init_yz}	
\end{figure}

A wire rope usually consists of twisted strands, where each strand is also composed of several twisted wires. Here we show a simulation of the twisting process of strands made of initially straight wires, which has been commonly employed to verify the applicability of the developed beam-to-beam contact formulations to cases with significantly small intersection angles, for example, see the relevant examples in \citet{meier2016finite}, \citet{konyukhov2018consistent}, and \citet{durville2010simulation,durville2012contact}. We consider two different cases of strands with either two or seven initially straight wires aligned with the $X$-axis. In both cases, we consider initially straight beams with $L\!=\!10\,\mathrm{m}$, and initially circular cross-sections of radius $R$. Two initial directors of the cross-section are chosen as ${{\boldsymbol{D}}_1} = {\boldsymbol{e}_2}$ and ${{\boldsymbol{D}}_2} = {\boldsymbol{e}_3}$. A Neo-Hookean material model with Young's modulus $E=210\,\mathrm{GPa}$, and Poisson's ratio $\nu=0.3$ is considered. We consider the following boundary conditions.
\begin{itemize}
	\item First, we constrain the axis displacement $\Delta \boldsymbol{\varphi}$ at the end $s=0$, i.e, 
	\begin{equation}
		\Delta \boldsymbol{\varphi}=\boldsymbol{0}\,\,\,\,\mathrm{at}\,\,s=0.	
	\end{equation}
	\item Second, an axial tension is imposed as a \textit{pre-deformation} by prescribing the $X$-displacement of the axis, as
	\begin{equation}
		\label{ex_twist_bdc_axial_tens}
		\Delta \boldsymbol{\varphi} = {\bar u_1}{\boldsymbol{e}_1}\,\,\,\,\mathrm{at}\,\,s=L,
	\end{equation}
	where we choose $\bar u_1=1\,\mathrm{m}$. This leads to a loss of contact due to the Poisson effect in the early phase of the twisting motion.
	\item Third, the end position of the axis at $s=L$ are prescribed such that they follow a circular path (see Fig.\,\ref{twist_2b_init_yz}), i.e.,
	\begin{subequations}
		\label{twist_cir_path}	
	\begin{equation}
		\left\{ {\begin{array}{*{20}{c}}
				{{\boldsymbol{\varphi}\cdot \boldsymbol{e}_2}-{c_2}}\\
				{{\boldsymbol{\varphi}\cdot \boldsymbol{e}_3}-{c_3}}
		\end{array}} \right\} = {\bar{\boldsymbol{\Lambda}}}\left\{ {\begin{array}{*{20}{c}}
				{{{\boldsymbol{\varphi}_0\cdot \boldsymbol{e}_2}} - {c_2}}\\
				{{{\boldsymbol{\varphi}_0\cdot \boldsymbol{e}_3}} - {c_3}}
		\end{array}} \right\},
	\end{equation}
	with
	\begin{equation}
		{\bar{\boldsymbol{\Lambda}}}\coloneqq\left[ {\begin{array}{*{20}{c}}
				{\cos {\bar \theta}}&{ - \sin {\bar \theta}}\\
				{\sin {\bar \theta}}&{\cos {\bar \theta}}
		\end{array}} \right],\,\,0 \le {\bar \theta} \le 2\pi,
	\end{equation}	
	\end{subequations}
	where $c_2$ and $c_3$ denote the fixed $Y$- and $Z$-coordinates of the selected rotation center, and ${\bar \theta}=2\pi\lambda_\mathrm{load}$ denotes the prescribed rotation angle with the load parameter $0\le\lambda_\mathrm{load}\le1$. Note that this prescribed circular path is imposed after applying the axial tension of Eq.\,(\ref{ex_twist_bdc_axial_tens}).
	\item Fourth, in order to prevent the rigid body rotation around the axis, we constrain the displacement components of directors $\boldsymbol{d}_1$ and $\boldsymbol{d}_2$ at the end $s=0$, as
	\begin{equation}
		\Delta {\boldsymbol{d}_{1}}\cdot{\boldsymbol{e}_3} = \Delta {\boldsymbol{d}_{2}}\cdot{\boldsymbol{e}_2} = 0\,\,\,\,\mathrm{at}\,\,s=0.
	\end{equation}
\end{itemize}
\noindent We choose the cutoff radius $r_\mathrm{c}=3R$, and $\varepsilon_\theta=0.25\pi$ in the global contact search.
\subsubsection{A strand of two wires}
\begin{figure}
	\centering
	\begin{subfigure}[b] {0.475\textwidth} \centering
		\includegraphics[width=\linewidth]{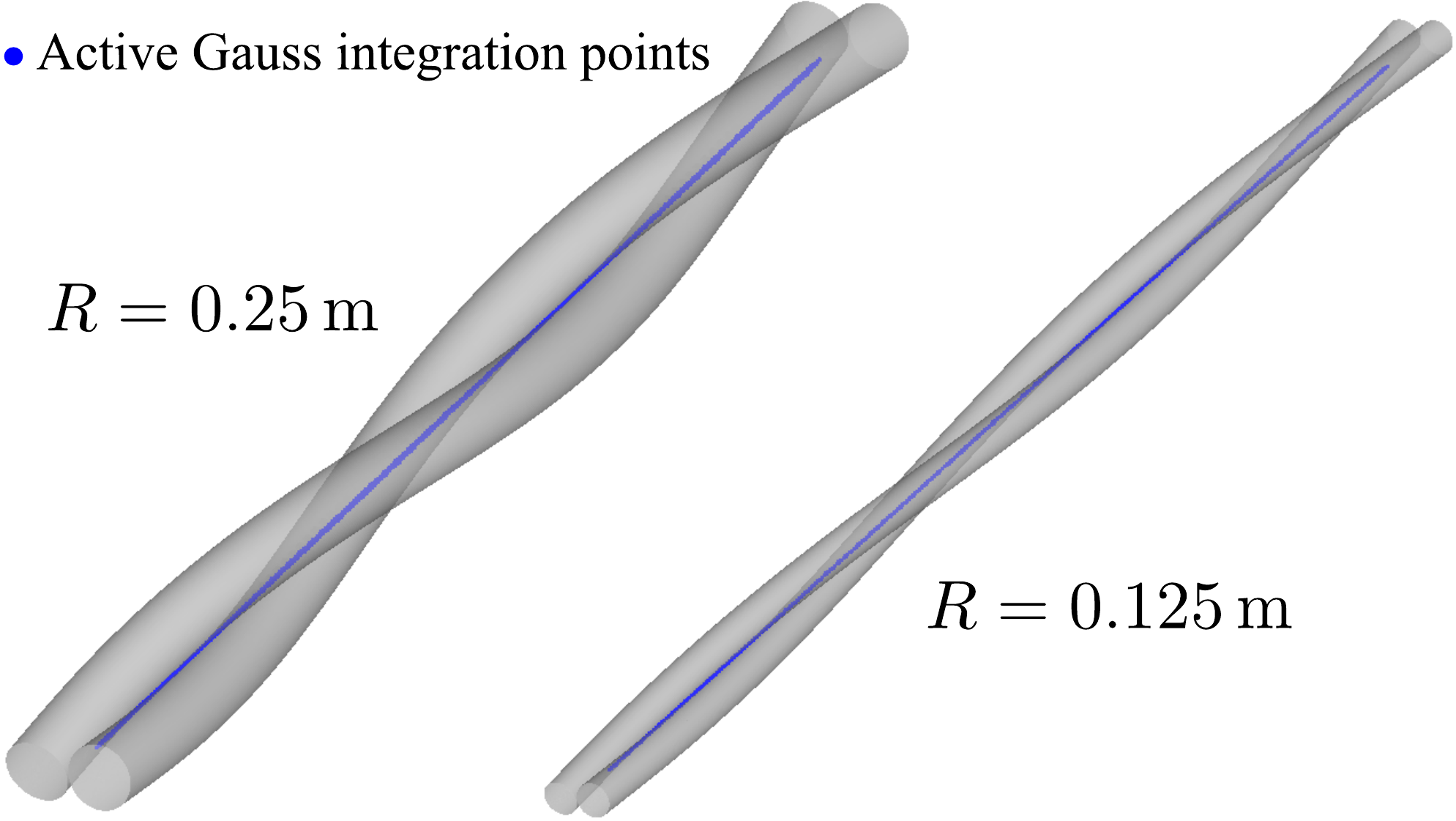}
	\end{subfigure}
	\caption{Twisting of a wire strand (two wires): The final deformed configurations for two cases of the initial cross-section's radius: $R=0.25\,\mathrm{m}$, and $R=0.125\,\mathrm{m}$. The blue dots indicate the active surface Gauss integration points. In both cases, we use $p=3$, $n_\mathrm{el}=40$, and $\epsilon_\mathrm{N}=10^2E/{L_0}$.}
	\label{twist_2beam_deform_len10}	
\end{figure}

\begin{figure}
	\centering
	\begin{subfigure}[b] {0.475\textwidth} \centering
		\includegraphics[width=\linewidth]{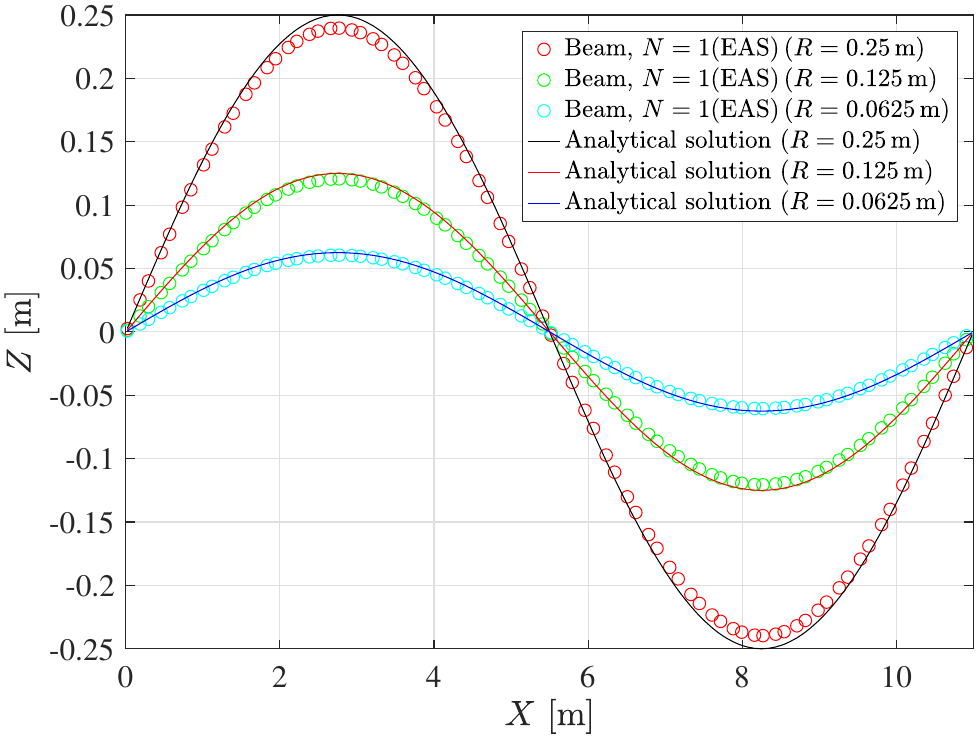}
	\end{subfigure}
	\caption{Twisting of a wire strand (two wires): Comparison of the final deformed configuration of the axis of the master body in the $X$-$Z$ plane for three different values of the initial cross-section radius: $R=0.25\,\mathrm{m},\,0.125\,\mathrm{m}$, and $0.0625\,\mathrm{m}$. Solid lines indicate the analytical solutions of Eq.\,(\ref{ex_str2b_helix_analytic}). In all cases, we use $p=3$, $n_\mathrm{el}=40$, and $\epsilon_\mathrm{N}=10^2E/{L_0}$.}
	\label{twist_2beam_caxis_config_plane}	
\end{figure}

\begin{figure}
	\centering
	\begin{subfigure}[b] {0.475\textwidth} \centering
		\includegraphics[width=\linewidth]{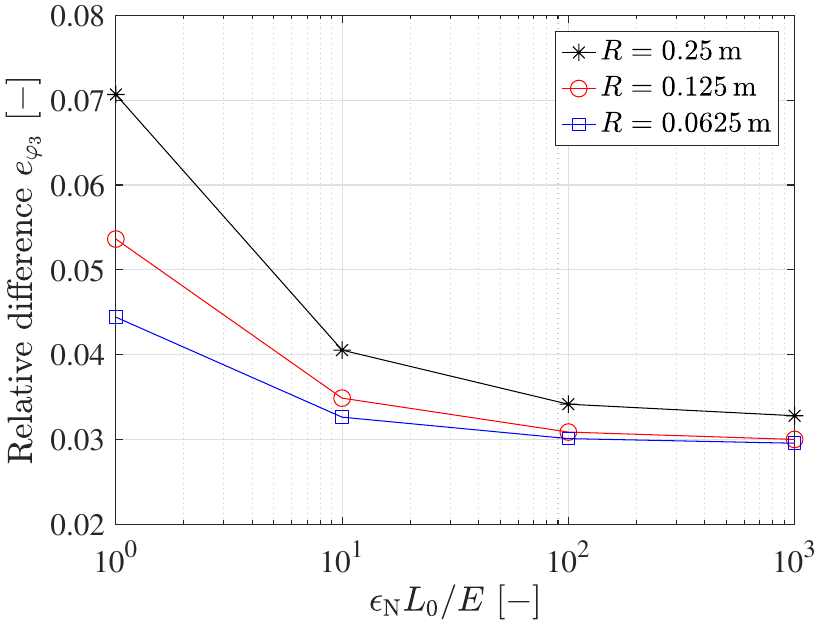}
	\end{subfigure}
	\caption{Twisting of a wire strand (two wires): The relative $L^2$-norm of Eq.\,(\ref{twist2b_rel_l2_norm_zdisp}) versus the normal contact penalty parameter. We use $N=1$ combined with the EAS method, and $p=3$ and $n_\mathrm{el}=40$.}
	\label{twist_2beam_caxis_diff_ref}	
\end{figure}

\begin{figure}
	\centering
	\begin{subfigure}[b] {0.475\textwidth} \centering
		\includegraphics[width=\linewidth]{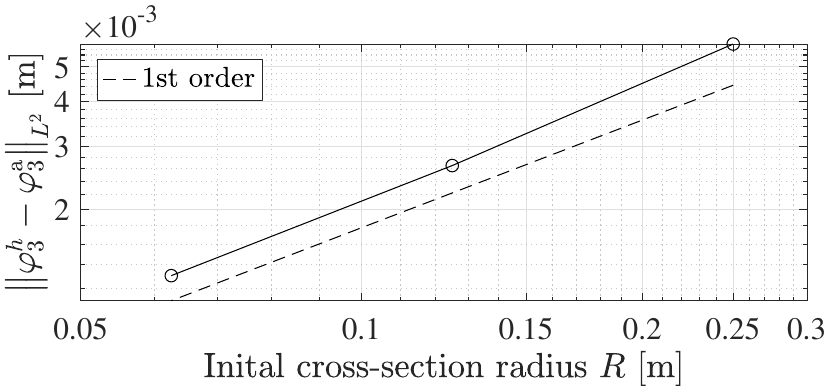}
		\caption{$L^2$-norm of the difference}
		\label{twist_2b_abs_diff_analytic}		
	\end{subfigure}
	\centering
	\begin{subfigure}[b] {0.475\textwidth} \centering
		\includegraphics[width=\linewidth]{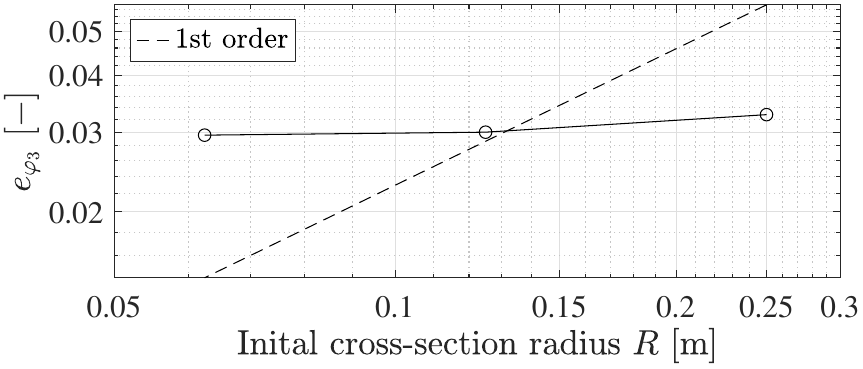}
		\caption{Relative $L^2$-norm of the difference}
		\label{twist_2b_rel_diff_analytic}		
	\end{subfigure}
	\caption{Twisting of wire strands (two wires): The convergence of the difference in the axial $Z$-coordinates of the final deformed configuration of the master body between the beam solution and the analytical solution of Eq.\,(\ref{ex_str2b_helix_analytic}), for the penalty parameter $\epsilon_\mathrm{N}=10^3E/{L_0}$. We use the beam formulation with $N=1$ combined with the EAS method, and $p=3$ and $n_\mathrm{el}=40$.}
	\label{twist_2beam_caxis_diff_ref_conv_r}	
\end{figure}

We first consider a strand of two wires initially parallel to the $X$-axis, and a circular path is prescribed by the angle ${\bar \theta}=2\pi\lambda_\mathrm{load}\,[\mathrm{rad}]$ at $s=L$ using Eq.\,(\ref{twist_cir_path}), see Fig.\,\ref{twist_2b_init_yz}. Fig.\,\ref{twist_2beam_deform_len10} shows the deformed configurations for two different values of the initial cross-section radius $R$. In Fig.\,\ref{twist_2beam_caxis_config_plane}, we compare the $Z$-coordinate in the axis of the final deformed configuration of the wires with that of an analytical circular helix of radius $R$, given by 
\begin{equation}
	\varphi^\mathrm{a}_3 = R\sin \left( {2\pi {\varphi^h_1}/\ell} \right),\,\,{\varphi^h_1} \in \left[ {0,\ell} \right],
	\label{ex_str2b_helix_analytic}
\end{equation}
where $\ell$ denotes the deformed length of the beam's axis after the pre-stretch, and ${\varphi^h_i}\coloneqq{\boldsymbol{\varphi}^h}\cdot\boldsymbol{e}_i$, $i\in\left\{1,2,3\right\}$. Further we define the \textit{relative} $L^2$-norm of the difference by
\begin{equation}
	\label{twist2b_rel_l2_norm_zdisp}
	{e_{{\varphi _3}}} \coloneqq \frac{{{{\left\| {\varphi _3^h - \varphi _3^{\rm{a}}} \right\|}_{{L^2}}}}}{{{{\left\| {\varphi _3^{\rm{a}}} \right\|}_{{L^2}}}}},
\end{equation}
where the $L^2$-norm for $u=u(s)$ in the domain $(0,L)\ni{s}$ is defined as
\begin{equation}
	{\left\| u \right\|_{{L^2}}} \coloneqq \sqrt {\int_0^L {{u^2}{\rm{d}}s} }.	
\end{equation}
It is shown in Fig.\,\ref{twist_2beam_caxis_config_plane} that the deformed axis of the beam slightly deviates from the analytical helix curve, and the difference decreases as the initial cross-sectional radius decreases. The difference from the analytical solution is mainly attributed to the following reasons, which are not considered in the analytical solution:
\begin{itemize}
	\item Cross-sectional contraction due to the pre-stretch,
	\item Cross-sectional deformations due to the contact interactions,
	\item A slight penetration allowed in the beam contact formulation using the penalty method.	
\end{itemize}
Fig.\,\ref{twist_2beam_caxis_diff_ref} shows that the difference decreases as the amount of penetration decreases due to increasing the penalty parameter. It is also seen in Fig.\,\ref{twist_2b_abs_diff_analytic} that the difference linearly decreases, as the slenderness ratio increases, since the $Z$-displacement is linearly proportional to the cross-sectional radius. However, the relative difference does not completely vanish but converges to a value around $e_{\varphi_3}=0.03$, which is attributed to the nominal transverse normal strain in the cross-section. Tables\,\ref{app_twist2b_nlstep_info_r025}\,-\,\ref{app_twist2b_nlstep_info_r00625} show the selected load increment sizes used in the results of Fig. \ref{twist_2beam_caxis_diff_ref} for initial cross-sectional radii $R=0.25\,\mathrm{m}$, $0.125\,\mathrm{m}$, and $0.0625\,\mathrm{m}$, respectively. Table \ref{app_twist2b_nlstep_info_r025_r0125_r00625} shows the selected number of sub-elements for the contact integral.
\subsubsection{A strand of seven wires}
\begin{figure}	
	\centering
	\begin{subfigure}[b] {0.425\textwidth} \centering
		\includegraphics[width=\linewidth]{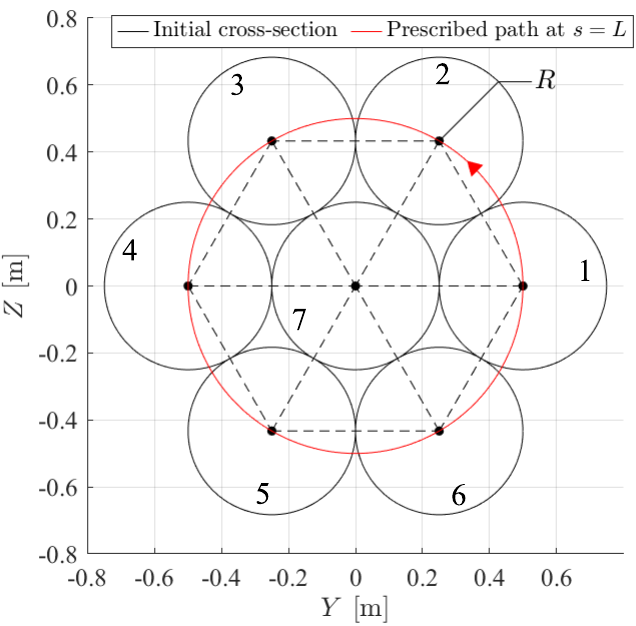}
	\end{subfigure}
	\caption{Twisting of a wire strand (seven wires): Initial configurations of strands in the cross-sectional $Y$-$Z$ plane, and the prescribed circular paths at the ends ($s=L$), indicated by the red circle.}
	\label{twist_7b_init_yz}	
\end{figure}

\begin{figure}	
	\centering
	\begin{subfigure}[b] {0.425\textwidth} \centering
		\includegraphics[width=\linewidth]{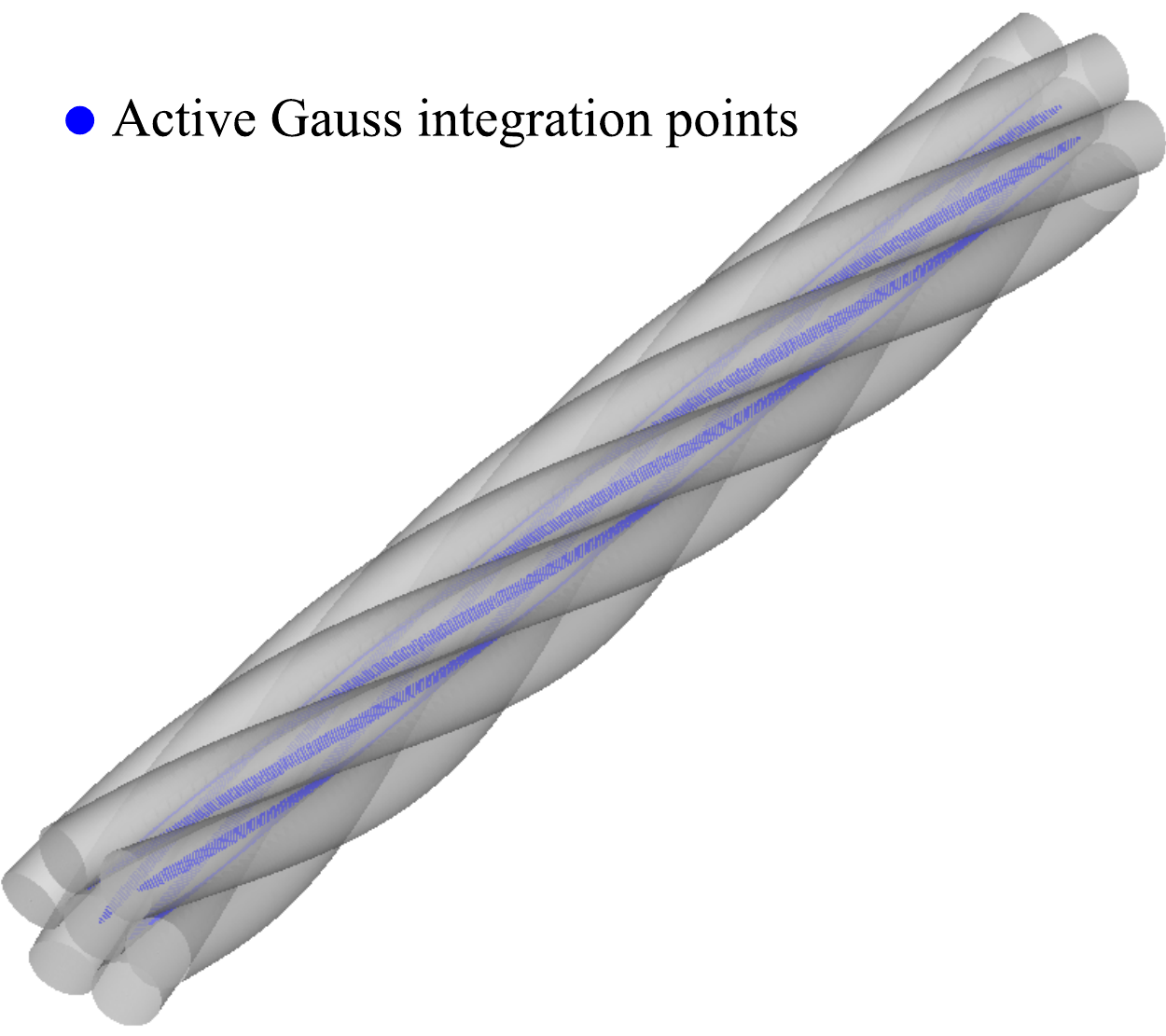}
	\end{subfigure}
	\caption{Twisting of wire strands (seven wires): The final deformed configuration. We use $p=3$, $n_\mathrm{el}=20$, and $\epsilon_\mathrm{N}=10E/{L_0}$, and select the outer beams as the slave bodies, and the inner beam as the master body (M1S6). The blue dots indicate the active surface Gauss points.}
	\label{twist_7b_final_deformed_m1s6}	
\end{figure}

\begin{figure}	
	\centering
	\begin{subfigure}[b] {0.4875\textwidth} \centering
	\includegraphics[width=\linewidth]{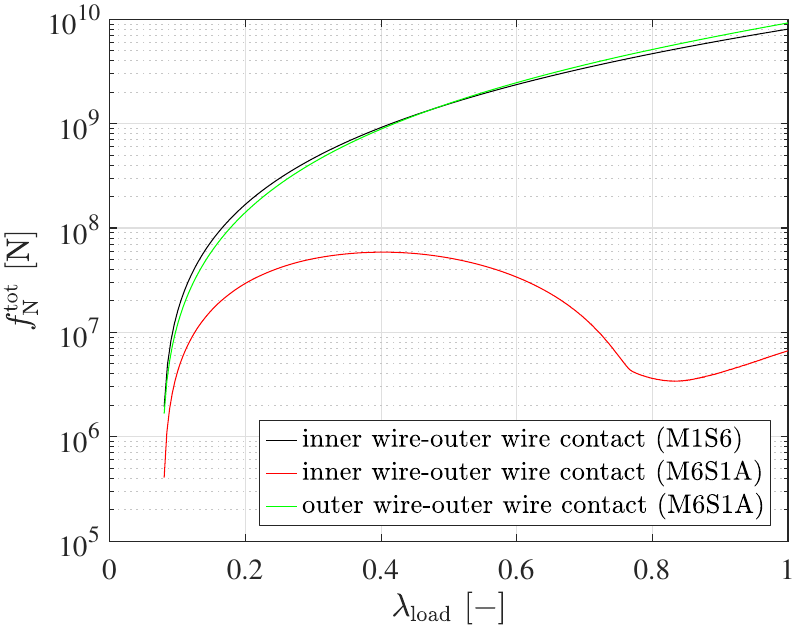}
	\end{subfigure}
	\caption{Twisting of a wire strand (seven wires): Change of the average of the total contact force from Eq.\,(\ref{tot_cf_intg}) between corresponding contact pairs. We use $p=3$, $n_\mathrm{el}=20$, and $\epsilon_\mathrm{N}=10E/{L_0}$.}
	\label{twist_7b_cf_ng_lparam}	
\end{figure}

We consider twisting of a strand with seven wires. The beams are initially aligned with the $X$-axis, and have the same length $L\!=\!10\,\mathrm{m}$, and a circular cross-section of radius 
$R=0.25\,\mathrm{m}$. Fig.\,\ref{twist_7b_init_yz} shows the initial arrangement of the wires and the prescribed rotation ${\bar \theta}=2\pi\lambda_\mathrm{load}\,[\mathrm{rad}]$ at the end $s=L$ of the outer wires. We consider the following two cases of selecting contact pairs:
\begin{itemize}
	\item M1S6: the inner beam is selected as master body, and the other six outer beams are selected as slave bodies. Thus, we have a total of six contact pairs.
	\item M6S1A: the inner beam is selected as slave body, and the other six outer beams are selected as master bodies. Additionally, the interaction between outer wires is considered. Thus, we have six \textit{inner wire-outer wire} contact pairs, and also six \textit{outer wire-outer wire} contact pairs, see Table \ref{app_twist7b_pair_info_m6s1a} for the chosen contact pair information.
\end{itemize}
Fig.\,\ref{twist_7b_final_deformed_m1s6} shows the final deformed configuration in case 1. In Fig.\,\ref{twist_7b_cf_ng_lparam}, the averge of the total contact forces between the contact pairs in each case is plotted. In the result for M6S1A, the contact force between inner and outer wires, and outer wires are plotted separately. Interestingly, the contact force between inner and outer wires in those two cases are significantly different, if we additionally consider the interaction between outer wires. Table \ref{app_twist7b_nlstep_info_2case} shows the selected load increment sizes in each case of M1S6 and M6S1A. In both cases, we use $n_\mathrm{el}^\mathrm{sub}=20$, and $m_\mathrm{el}^\mathrm{sub}=300$ for the contact integral.
\section{Conclusions}\label{sec_conclusion} 
In this paper, we present an isogeometric finite element formulation for beams and beam-to-beam frictionless contact, based on the kinematics of Cosserat rod with unconstrained directors. The beam cross-sectional deformation is efficiently and accurately described by unconstrained directors of an arbitrary order. The beam contact formulation is based on a Gauss point-to-surface contact algorithm, where the impenetrability constraint is enforced by a penalty method combined with an active set method. Further we present a geometrical approach to efficiently determine an initial guess in the local Newton-Raphson iteration of the closest point projection. In numerical examples, we show that the beam contact formulation can provide contact pressure distributions which agree well with brick element solutions while requiring much less DOFs. The following areas could be interesting future research directions.
\begin{itemize}
\item Alleviation of locking: In Section \ref{num_ex_sliding_2b_c2_def}, it is shown that the curvature-thickness locking may lead to an unphysical oscillation of displacements in sliding contact. An extension of the assumed natural strain method in \citet{betsch1995assumed} to a higher order basis functions would be interesting. Further, an alleviation of contact locking (overconstrained system) due to the contact constraints will be investigated further in the future, e.g., using a mortar-type discretization method, a deactivation of superfluous Gauss integration points in the active set (outer) loop, or an adaptive control of the penalty parameter.

\item A mixed-variational formulation could alleviate locking and improve robustness in larger load steps, see \citet{wackerfuss2009mixed}. The static condensation of the strains associated with higher order directors could also improve the efficiency of the beam formulation.
%
\item It is observed that an increase of the penalty parameter typically requires larger number of surface Gauss integration points for the contact integral and larger number of load steps, which eventually makes the solution process less efficient. One can consider an augmented Lagrangian method or Nitsche's method in order to exactly satisfy the impenetrability constraint using low or moderate penalty parameters.
\item In the surface-to-surface beam contact formulation, we consider contact between smooth lateral surface domains only. If a non-smooth edge exists in the lateral surface, it is required to additionally consider a contact integral along the edge, for example, see \cite{litewka2002contact}. Further, the cross-sections at the ends of axis may also contribute to contact interactions. For example, their boundary edges may undergo \textit{edge-to-surface} or \textit{edge-to-edge} contact.
\item The present beam contact formulation can be extended to incorporate tangential contact conditions including sticking and sliding friction. In sliding friction, an additional history variable is required for the amount of irreversible tangential slip, and it typically requires much smaller load increment sizes for accuracy, compared with frictionless or sticking contact formulations.
\end{itemize}
\section*{Acknowledgement}
M.-J Choi would like to gratefully acknowledge the financial support of a postdoctoral research fellowship from the Alexander von Humboldt Foundation in Germany.


\begin{appendices}
	\section{}
	\subsection{Gauss point-to-axis closest point projection}
	\label{app_glob_contact_search}
	We first define the distance between the material points on the axes of slave and master bodies, as
	\begin{align}
		\label{dist_f_axis_cpp}
		d_{\varphi}\big(\boldsymbol{\varphi}^{(1)},\boldsymbol{\varphi}^{(2)}\big({{\xi^1_{(2)}}}\big)\big)\coloneqq\left\| {{{\boldsymbol{\varphi }}^{(1)}} - {{\boldsymbol{\varphi }}^{(2)}}\big({{\xi^1_{(2)}}}\big)} \right\|.
	\end{align}	
	The convective axial coordinate $\widetilde{\xi^1}$ of the closest point in the master body to a given point $\boldsymbol{\varphi}^{(1)}$ along the axis of the slave body is determined as a solution of the following unilateral minimal distance problem:
	\begin{equation}\label{ul_dist_prob}
		{{\widetilde {\xi^1}} } \coloneqq \arg \mathop {\min }\limits_{{\xi^1_{(2)}}} d_{\varphi}\big({\boldsymbol{\varphi}^{(1)}},\boldsymbol{\varphi}^{(2)}\big({{\xi^1_{(2)}}}\big)\big),
	\end{equation}
	where
	\begin{equation}\label{ul_dist_prob_scnst}
	{{\xi^1_\mathrm{min}}\lt{{\xi^1_{(2)}}\lt{\xi^1_\mathrm{max}}}}.
	\end{equation}
	Hereafter, for brevity, we often omit the argument $\xi^1_{(2)}$. The solution of Eq.\,(\ref{ul_dist_prob}) satisfies the first order necessary condition
	\begin{equation}\label{ul_dist_prob_nec_cond}
		{f_\varphi} \coloneqq {\boldsymbol{\varphi }}^{(2)}_{,1}\cdot\left( {{{\boldsymbol{\varphi }}^{(1)}} - {{\boldsymbol{\varphi }}^{(2)}}} \right)= 0.
	\end{equation}
	This is a single nonlinear equation with respect to ${{\xi^1_{(2)}}}$, which can be iteratively solved using a Newton-Raphson iteration. The linearization of Eq.\,(\ref{ul_dist_prob_nec_cond}) leads to
	\begin{equation}\label{ul_dist_prob_lin}
		f_{\varphi}^ * \Delta {{\xi^1_{(2)}}} =  - {f_{\varphi}},
	\end{equation}
	such that the convective axial coordinate is updated by 
	\begin{equation}
		\prescript{(i)}{}\xi^1_{(2)}=\prescript{(i-1)}{}\xi^1_{(2)} + \Delta\xi^1_{(2)},\,\,i=1,2,...,
	\end{equation}
	until a convergence criterion $\left\lvert {{f_\varphi }} \right\rvert < \varepsilon _{{\rm{cpp}}}^\varphi$ is satisfied, where $\varepsilon _{{\rm{cpp}}}^\varphi>0$ is a chosen tolerance, and 
	\begin{align}\label{ul_def_f_ul_st}
		f_{\varphi}^* &\coloneqq {\partial{f_\varphi}/\partial{{\xi^1_{(2)}}}}\nonumber\\
		&=\left( {{{\boldsymbol{\varphi }}^{(1)}} - {{\boldsymbol{\varphi }}^{(2)}}} \right) \cdot {\boldsymbol{\varphi }}_{,{11}}^{(2)} - {\boldsymbol{\varphi }}_{,1}^{(2)}\cdot{\boldsymbol{\varphi }}_{,1}^{(2)},
	\end{align}
	and $\prescript{(0)}{}\xi^1_{(2)}$ is a chosen initial guess. Assuming ${f_\varphi^ *}\ne0$, from Eq.\,(\ref{ul_dist_prob_lin}), we obtain
	\begin{equation}\label{lcontact_xi2_inc}
		\Delta {\xi ^1_{(2)}} =  - {f_{\varphi}}/f_{\varphi}^{*}.
	\end{equation}
	
	
	\subsection{Contact variational form}
	\label{app_lin_contact_var_form}
	\subsubsection{Linearization}
	Taking the directional derivative of Eq.\,(\ref{var_normal_gap_prv}) leads to \citep{wriggers2006}
	\begin{align}
		\label{app_compact_del_del_gN}
		\Delta \delta {g_{\rm{N}}} &= {\left\{ {\begin{array}{*{20}{c}}
					{\delta {{{\boldsymbol{\bar u}}}_{,1}}}\\
					{\delta {{{\boldsymbol{\bar u}}}_{,2}}}
			\end{array}} \right\}^{\mathrm{T}}}{{\boldsymbol{k}}_{{\mathrm{uu}}}}\left\{ {\begin{array}{*{20}{c}}
				{\Delta {{{\boldsymbol{\bar u}}}_{,1}}}\\
				{\Delta {{{\boldsymbol{\bar u}}}_{,2}}}
		\end{array}} \right\} \nonumber\\
	&+ {\left\{ {\begin{array}{*{20}{c}}
					{\delta {{{\boldsymbol{\bar u}}}_{,1}}}\\
					{\delta {{{\boldsymbol{\bar u}}}_{,2}}}
			\end{array}} \right\}^\mathrm{T}}{{\boldsymbol{k}}_{{\mathrm{ux}}}}{\Delta \bar{\boldsymbol{\xi }}}+ {\delta \bar{\boldsymbol{\xi }}^\mathrm{T}}{\boldsymbol{k}}_{{\mathrm{ux}}}^{\mathrm{T}}\left\{ {\begin{array}{*{20}{c}}
				{\Delta {{{\boldsymbol{\bar u}}}_{,1}}}\\
				{\Delta {{{\boldsymbol{\bar u}}}_{,2}}}
		\end{array}} \right\}\nonumber\\
	&+ {{\delta \bar{\boldsymbol{\xi }}^\mathrm{T}}}{{\boldsymbol{k}}_{{\mathrm{xx}}}}{\Delta \bar{\boldsymbol{\xi }}}
	\end{align}
	with the following matrix operators
	\begin{align}
		\label{app_kuu_def}
		{{\boldsymbol{k}}_{{\mathrm{uu}}}} \coloneqq \left[ {\begin{array}{*{20}{c}}
				{{{\bar a}^{11}}{{{\boldsymbol{\bar \nu }}_t}} \otimes {{{\boldsymbol{\bar \nu }}_t}}}&{{{\bar a}^{12}}{{{\boldsymbol{\bar \nu }}_t}} \otimes {{{\boldsymbol{\bar \nu }}_t}}}\\
				{{\mathrm{sym}}.}&{{{\bar a}^{22}}{{{\boldsymbol{\bar \nu }}_t}} \otimes {{{\boldsymbol{\bar \nu}}_t}}}
		\end{array}} \right]{g_{\mathrm{N}}},
	\end{align}
	\begin{align}
	{{\boldsymbol{k}}_{{\mathrm{ux}}}} \coloneqq \left[ {\begin{array}{*{20}{c}}
			{\big({g_\mathrm{N}}{{{\bar a}^{1\gamma }}{{\bar b}_{\gamma 1}}}-1\big){{{\boldsymbol{\bar \nu }}}}_t}&{{g_\mathrm{N}}{{\bar a}^{1\gamma }}{{\bar b}_{\gamma 2}}{{{\boldsymbol{\bar \nu }}_t}}}\\
			{{g_\mathrm{N}}{{\bar a}^{2\gamma }}{{\bar b}_{\gamma 1}}{{{\boldsymbol{\bar \nu }}}}_t}&{\big({g_\mathrm{N}}{{{\bar a}^{2\gamma }}{{\bar b}_{\gamma 2}}}-1\big){{{\boldsymbol{\bar \nu }}_t}}}
	\end{array}} \right],
	\end{align}
	and
	\begin{align}
		{{\boldsymbol{k}}_{{\mathrm{xx}}}} &\coloneqq \left[ {\begin{array}{*{20}{c}}
				{{{\bar a}^{\alpha \gamma }}{{\bar b}_{\alpha 1}}{{\bar b}_{\gamma 1}}}&{{{\bar a}^{\alpha \gamma }}{{\bar b}_{\alpha 1}}{{\bar b}_{\gamma 2}}}\\
				{{\mathrm{sym}}.}&{{{\bar a}^{\alpha \gamma }}{{\bar b}_{\alpha 2}}{{\bar b}_{\gamma 2}}}
		\end{array}} \right]{g_{\mathrm{N}}}\nonumber\\
	&-\left[ {\begin{array}{*{20}{c}}
				{{{\bar b}_{11}}}&{{\bar b}_{12}}\\
				{{\mathrm{sym}}.}&{{\bar b}_{22}}
		\end{array}} \right].
	\end{align}
	$\delta{\bar {\boldsymbol{\xi}}}$ in Eq.\,(\ref{app_compact_del_del_gN}) can be expressed in terms of $\delta{\boldsymbol{\bar u}}$. Taking the first order variation of Eq.\,(\ref{cpp_uni_lat}), we have
	\begin{align}
		\label{express_del_xi_del_x_u}
		{\bar{{\boldsymbol{f}}^{*}}}\delta \bar{\boldsymbol{\xi }} = {{\boldsymbol{\Xi }}_{\mathrm{f}}}\left\{ {\begin{array}{*{20}{c}}
				{\delta {{\boldsymbol{q}}^{(1)}}}\\
				{\delta \widetilde{\boldsymbol{q}}}
		\end{array}} \right\},
	\end{align}
	where we obtain \citep{wriggers2006}
	\begin{align}
		\label{mat_f_u_at_xi_bar}
		{\bar{{\boldsymbol{f}}^{*}}}&\coloneqq{{\boldsymbol{f}}^{*}}({\bar{\boldsymbol{\xi}}}) \nonumber\\
		&= {g_{\mathrm{N}}}\left[ {\begin{array}{*{20}{c}}
				{{{\bar b}_{11}}}&{{{\bar b}_{12}}}\\
				{{\mathrm{sym}}.}&{{{\bar b}_{22}}}
		\end{array}} \right] - \left[ {\begin{array}{*{20}{c}}
				{{{\bar a}_{11}}}&{{{\bar a}_{12}}}\\
				{{\mathrm{sym}}.}&{{{\bar a}_{22}}}
			\end{array}} \right]
	\end{align}
	by evaluating Eq.\,(\ref{f_u_st_1}) at ${{\boldsymbol{\xi }}^{(2)}} = {\boldsymbol{\bar \xi }}$, and we define the operator
	\begin{subequations}
		\begin{equation}
			{\boldsymbol{\Xi }_\mathrm{f}} \coloneqq \left[ {\setlength{\arraycolsep}{1pt}
				\renewcommand{\arraystretch}{1.3}\begin{array}{*{20}{c}}
					{{\boldsymbol{\Xi}}_\mathrm{f}^{11}}&{{\boldsymbol{\Xi}}_\mathrm{f}^{12}}\\
					{{\boldsymbol{\Xi}}_\mathrm{f}^{21}}&{{\boldsymbol{\Xi}}_\mathrm{f}^{22}}
			\end{array}} \right],
		\end{equation}		
		where
		\begin{align}
			{\boldsymbol{\Xi}}_\mathrm{f}^{11} &\coloneqq {-{{\boldsymbol{\bar a}}}_1^\mathrm{T}}{\boldsymbol{\Pi}_{(1)}^\mathrm{T}},\\
			{\boldsymbol{\Xi}}_\mathrm{f}^{12} &\coloneqq {\bar{\boldsymbol{a}}}_1^\mathrm{T}{\bar{\boldsymbol{\Pi}}^\mathrm{T}}\nonumber\\
			&-{g_{\mathrm{N}}}{{\bar{\boldsymbol{\nu }}}}_t^\mathrm{T}\left({\bar{\boldsymbol{\Pi}}^\mathrm{T}}{{\left(  \bullet\right)}_{,1}}+{{\bar{\boldsymbol{\Pi}}}_{,\zeta^\alpha}^\mathrm{T}}{{\bar \zeta}^\alpha_{,1}}\right),\\
			{\boldsymbol{\Xi}}_\mathrm{f}^{21} &\coloneqq {-{\bar{\boldsymbol{a}}}_2^\mathrm{T}}{\boldsymbol{\Pi}_{(1)}^\mathrm{T}},\\
			{\boldsymbol{\Xi}}_\mathrm{f}^{22} &\coloneqq {\bar{\boldsymbol{a}}}_2^\mathrm{T}{{\bar{\boldsymbol{\Pi}}}^\mathrm{T}}-{{g_{\mathrm{N}}}{{\bar{\boldsymbol{\nu}}}}_t^\mathrm{T}{{\bar{\boldsymbol{\Pi}}}_{,\zeta^\alpha}^\mathrm{T}}{{\bar \zeta}^\alpha_{,2}}}.
		\end{align}
	\end{subequations}
	Here we assume the matrix $\bar{{\boldsymbol{f}}^{*}}$ is invertible, then we have
	\begin{equation}
		\label{app_del_xi_bar}
		\delta \bar{\boldsymbol{\xi }} = {\bar{{\boldsymbol{f}}^{*}}}^{ - 1}{{\boldsymbol{\Xi }}_{\mathrm{f}}}\left\{ {\begin{array}{*{20}{c}}
				{\delta {{\boldsymbol{q}}^{(1)}}}\\
				{\delta \widetilde{\boldsymbol{q}}}
		\end{array}} \right\}.
	\end{equation}
	We further define a matrix operator ${\boldsymbol{\Xi }_\mathrm{u}}$ such that
	\begin{subequations}
	\begin{equation}
		\label{app_mat_op_del_u_deriv}
		\left\{ {\begin{array}{*{20}{c}}
				{\delta {{{\boldsymbol{\bar u}}}_{,1}}}\\
				{\delta {{{\boldsymbol{\bar u}}}_{,2}}}
		\end{array}} \right\} = {\boldsymbol{\Xi }_\mathrm{u}}\left\{ {\begin{array}{*{20}{c}}
				{\delta {{\boldsymbol{q}}^{(1)}}}\\
				{\delta {\boldsymbol{\widetilde q}}}
		\end{array}} \right\}
	\end{equation}
	with
	\begin{equation}
		{\boldsymbol{\Xi }_\mathrm{u}} \coloneqq \left[ {\begin{array}{*{20}{c}}
				{{{\boldsymbol{0}}_{3 \times {n_\mathrm{cs}}}}}&{{{\bar{\boldsymbol{\Pi}}}_{,\zeta^\alpha}^\mathrm{T}}{{\bar \zeta}^\alpha_{,1}}+{{\bar{\boldsymbol{\Pi}}}^\mathrm{T}}{(\bullet)_{,1}}}\\
				{{{\boldsymbol{0}}_{3 \times {n_\mathrm{cs}}}}}&{{\bar {\boldsymbol{\Pi}}}_{,\zeta^{\alpha}}^\mathrm{T}}{{\bar \zeta}^{\alpha}_{,2}}
		\end{array}} \right].
	\end{equation}
	\end{subequations}
	Then, substituting Eqs.\,(\ref{app_del_xi_bar}) and (\ref{app_mat_op_del_u_deriv}) into Eq.\,(\ref{app_compact_del_del_gN}), we finally obtain 
	\begin{align}
		\label{dir_deriv_weak_form_3d_gen}
		\Delta \delta {g_{\mathrm{N}}} = {\left\{ {\begin{array}{*{20}{c}}
					{\delta {{\boldsymbol{q}}^{(1)}}}\\
					{\delta {\boldsymbol{\widetilde q}}}
			\end{array}} \right\}^{\mathrm{T}}}{{\boldsymbol{k}}}^\mathrm{G}_\mathrm{N}\left\{ {\begin{array}{*{20}{c}}
				{\Delta {{\boldsymbol{q}}^{(1)}}}\\
				{\Delta {\boldsymbol{\widetilde q}}}
		\end{array}} \right\},
	\end{align}
	where
	\begin{align}
		\label{deriv_kg_n}
		{{\boldsymbol{k}}}^\mathrm{G}_\mathrm{N}\coloneqq{{{\tilde {\boldsymbol{k}}}_{{\mathrm{uu}}}} + {{\tilde {\boldsymbol{k}}}_{{\mathrm{ux}}}} + {\tilde {\boldsymbol{k}}}_{{\mathrm{ux}}}^{\mathrm{T}} + {{\tilde {\boldsymbol{k}}}_{{\mathrm{xx}}}}},
	\end{align}
	with
	\begin{subequations}
		\label{lin_con_var_mat_def}		
		\begin{align}
			{{\tilde{\boldsymbol{k}}}_{{\mathrm{uu}}}} &\coloneqq {\boldsymbol{\Xi }}_{\mathrm{u}}^{\mathrm{T}}{{\boldsymbol{k}}_{{\mathrm{uu}}}}{{\boldsymbol{\Xi }}_{\mathrm{u}}},\\
			{{\tilde{\boldsymbol{k}}}_{{\mathrm{ux}}}} &\coloneqq {\boldsymbol{\Xi }}_{\mathrm{u}}^{\mathrm{T}}{{\boldsymbol{k}}_{{\mathrm{ux}}}}{\bar{\boldsymbol{f}^*}^{ - 1}}{{\boldsymbol{\Xi }_\mathrm{f}}},\\
			{\tilde{{\boldsymbol{k}}}_{{\mathrm{xx}}}} &\coloneqq {\boldsymbol{\Xi }}_{\mathrm{f}}^{\mathrm{T}}{\bar{\boldsymbol{f}^*}^{ - \mathrm{T}}}{{\boldsymbol{k}}_{{\mathrm{xx}}}}{\bar{\boldsymbol{f}^*}^{ - 1}}{{\boldsymbol{\Xi }}_{\mathrm{f}}}.
		\end{align}
	\end{subequations}
	\subsubsection{Spatial discretization}
	In the material and geometric part of the tangent stiffness matrices, we define
	\begin{subequations}
		\label{disc_km_n}
		\begin{align}
			{{{\Bbbk}}^\mathrm{M}_\mathrm{N}}=\left[\renewcommand{\arraystretch}{1.25}{\begin{array}{*{20}{c}}
					\big({{{\Bbbk}}^\mathrm{M}_\mathrm{N}}\big)_{11}&\big({{{\Bbbk}}^\mathrm{M}_\mathrm{N}}\big)_{12}\\
					{\mathrm{sym.}}&\big({{{\Bbbk}}^\mathrm{M}_\mathrm{N}}\big)_{22}
			\end{array}} \right],
		\end{align}
		with
		\begin{align}
			\big({{{\Bbbk}}^\mathrm{M}_\mathrm{N}}\big)_{11}&\coloneqq{{\Bbb{N}_{e}^{(1)\mathrm{T}}}{{\boldsymbol{\Pi }}^{(1)}}{{{\boldsymbol{\bar \nu }}}_t} \otimes {{\Bbb{N}_{e}^{(1)\mathrm{T}}}{\boldsymbol{\Pi }}^{(1)}}{{{\boldsymbol{\bar \nu }}}_t}},\\
			\big({{{\Bbbk}}^\mathrm{M}_\mathrm{N}}\big)_{12}&\coloneqq{ - {\Bbb{N}_{e}^{(1)\mathrm{T}}}{{\boldsymbol{\Pi }}^{(1)}}{{{\boldsymbol{\bar \nu }}}_t} \otimes {\bar{\Bbb{N}}_{\bar e}^{\mathrm{T}}}{\boldsymbol{\bar \Pi }}\,{{{\boldsymbol{\bar \nu }}}_t}},\\			
			\big({{{\Bbbk}}^\mathrm{M}_\mathrm{N}}\big)_{22}&\coloneqq{\bar{\Bbb{N}}_{\bar e}^{\mathrm{T}}}{{\boldsymbol{\bar \Pi }}\,{{{\boldsymbol{\bar \nu }}}_t} \otimes {\bar{\Bbb{N}}_{\bar e}^{\mathrm{T}}}{\boldsymbol{\bar \Pi }}\,{{{\boldsymbol{\bar \nu }}}_t}},			
		\end{align}
	\end{subequations}
	and
	\begin{align}
		\label{disc_kg_n}
		{{\Bbbk}}^\mathrm{G}_\mathrm{N}\coloneqq{{{\tilde {\Bbbk}}^e_{{\mathrm{uu}}}} + {{\tilde {\Bbbk}}^e_{{\mathrm{ux}}}} + {\tilde {\Bbbk}}_{{\mathrm{ux}}}^{e\mathrm{T}} + {{\tilde {\Bbbk}}^e_{{\mathrm{xx}}}}},
	\end{align}
	where
	\begin{subequations}
	\label{c_disc_matxi_f}	
		\begin{align}
			{{\tilde{\Bbbk}}^e_{{\mathrm{uu}}}} &\coloneqq {{\boldsymbol{\Xi }}_{\mathrm{u}}^{{e\mathrm{T}}}}{{\boldsymbol{k}}_{{\mathrm{uu}}}}{{\boldsymbol{\Xi }}^e_{\mathrm{u}}},\\
			{{\tilde{\Bbbk}}^e_{{\mathrm{ux}}}} &\coloneqq {{\boldsymbol{\Xi }}_{\mathrm{u}}^{e\mathrm{T}}}{{\boldsymbol{k}}_{{\mathrm{ux}}}}{\bar{\boldsymbol{f}^*}^{ - 1}}{{\boldsymbol{\Xi }^{e}_\mathrm{f}}},\\
			{\tilde{{\Bbbk}}^e_{{\mathrm{xx}}}} &\coloneqq {{\boldsymbol{\Xi }}_{\mathrm{f}}^{{e\mathrm{T}}}}{\bar{\boldsymbol{f}^*}^{ - \mathrm{T}}}{{\boldsymbol{k}}_{{\mathrm{xx}}}}{\bar{\boldsymbol{f}^*}^{ - 1}}{{\boldsymbol{\Xi }}^e_{\mathrm{f}}}.
		\end{align}
	\end{subequations}
	Those matrices ${\boldsymbol{\Xi }}_{\mathrm{u}}^e$ and ${{\boldsymbol{\Xi }}_{\mathrm{f}}^e}$ are defined for the $e$-th element, as
	\begin{equation}
		\left.\begin{array}{*{20}{c}}
			{\boldsymbol{\Xi }}_{\rm{u}}^e &\coloneqq {\left[ {\begin{array}{*{20}{c}}
						{\tilde {\boldsymbol{\Xi }}_{\rm{u}}^1}&{\tilde {\boldsymbol{\Xi }}_{\rm{u}}^2}& \cdots &{\tilde {\boldsymbol{\Xi }}_{\rm{u}}^{{n_e}}}
				\end{array}} \right]_{6 \times 2{n_e}n_{{\rm{cs}}}}}\\
			{\boldsymbol{\Xi }}_{\rm{f}}^e &\coloneqq {\left[ {\begin{array}{*{20}{c}}
						{\tilde {\boldsymbol{\Xi }}_{\rm{f}}^1}&{\tilde {\boldsymbol{\Xi }}_{\rm{f}}^2}& \cdots &{\tilde {\boldsymbol{\Xi }}_{\rm{f}}^{{n_e}}}
				\end{array}} \right]_{2 \times 2{n_e}n_{{\rm{cs}}}}}
		\end{array}\right\},
	\end{equation}
	where
	\begin{equation}
		\label{c_disc_matxi_u}
		{\tilde {\boldsymbol{\Xi }}^I_{\mathrm{u}}} \coloneqq \left[ {\begin{array}{*{20}{c}}
				{{{\boldsymbol{0}}_{3 \times {n_\mathrm{cs}}}}}&{{{{\bar{\boldsymbol{\Pi}}}_{,\zeta^\alpha}^\mathrm{T}}{{\bar \zeta}^\alpha_{,1}}{{\bar N}_{I}}+{{\bar{\boldsymbol{\Pi}}}^\mathrm{T}}{{\bar N}_{I,1}}}}\\
				{{{\boldsymbol{0}}_{3 \times {n_\mathrm{cs}}}}}&{{\bar{\boldsymbol{\Pi}}}_{,\zeta^{\alpha}}^\mathrm{T}}{{\bar \zeta}^{\alpha}_{,2}}{{\bar N}_I}
		\end{array}} \right],
	\end{equation}
	and
	\begin{subequations}
	\begin{equation}
		{\tilde {\boldsymbol{\Xi }}^I_{\mathrm{f}}} \coloneqq \left[ {\begin{array}{*{20}{c}}
				{({\tilde {\boldsymbol{\Xi }}^I_{\mathrm{f}}})_{11}}&{({\tilde {\boldsymbol{\Xi }}^I_{\mathrm{f}}})_{12}}\\
				{({\tilde {\boldsymbol{\Xi }}^I_{\mathrm{f}}})_{21}}&{({\tilde {\boldsymbol{\Xi }}^I_{\mathrm{f}}})_{22}}
		\end{array}} \right]_{2 \times 2{n_{{\rm{cs}}}}},
	\end{equation}
	with
	\begin{align}
		({\tilde {\boldsymbol{\Xi }}^I_{\mathrm{f}}})_{11}&\coloneqq{ - {{{\boldsymbol{\bar a}}}_1^\mathrm{T}}}{\boldsymbol{\Pi}_{(1)}^\mathrm{T}}{N_I},\\
		({\tilde {\boldsymbol{\Xi }}^I_{\mathrm{f}}})_{12}&\coloneqq{{\bar{\boldsymbol{a}}}_1^\mathrm{T}}{\bar{\boldsymbol{\Pi}}^\mathrm{T}}{{\bar N}_I}\nonumber\\
			&-{g_{\mathrm{N}}}{{\bar{\boldsymbol{\nu }}}}_t^\mathrm{T}({\bar{\boldsymbol{\Pi}}^\mathrm{T}}{{{\bar N}_{I,1}}}+{{\bar{\boldsymbol{\Pi}}}_{,\zeta^\alpha}^\mathrm{T}}{{\bar \zeta}^\alpha_{,1}}{{\bar N}_I}),\\
		({\tilde {\boldsymbol{\Xi }}^I_{\mathrm{f}}})_{21}&\coloneqq{ - {{\bar{\boldsymbol{a}}}_2^\mathrm{T}}}{\boldsymbol{\Pi}_{(1)}^\mathrm{T}}{N_I},\\
		({\tilde {\boldsymbol{\Xi}}^I_{\mathrm{f}}})_{22}&\coloneqq({{{\bar{\boldsymbol{a}}}_2^\mathrm{T}}{{\bar{\boldsymbol{\Pi}}}^\mathrm{T}}-{g_{\mathrm{N}}}{{\bar{\boldsymbol{\nu}}}}_t^\mathrm{T}{{\bar{\boldsymbol{\Pi}}}_{,\zeta^\alpha}^\mathrm{T}}{{\bar \zeta}^\alpha_{,2}}}){{\bar N}_I}.
	\end{align}
\end{subequations}
\subsection{Weak enforcement of displacement boundary conditions}
\label{app_contact_enforce_pre_trans}
We present a weak enforcement of the displacement boundary condition on the lateral boundary surface, which is utilized in the numerical example of Section \ref{contact_num_ex_ring_flat}. Let $\mathcal{S}^\mathrm{D}_0\subset\mathcal{S}^\mathrm{L}_0$ be a region of the initial boundary surface where the displacement $\boldsymbol{u}$ is prescribed, i.e., $\boldsymbol{u}={\bar {\boldsymbol{u}}}$ on $\mathcal{S}^\mathrm{D}_0$, ${\bar {\boldsymbol{u}}}\in\Bbb{R}^3$. 
\subsubsection{Variational formulation}
We employ a penalty method, and the penalty functional can be expressed, using Eq.\,(\ref{inf_area_init_config}), as
\begin{align}
	\label{app_pre_disp_pen_func}
	{\Pi _{\mathrm{D}}} &\coloneqq \frac{1}{2}\int_{\mathcal{S}_0^{\rm{D}}} {{\epsilon_{\rm{D}}}\,{{\left\|{\boldsymbol{u} - {\bar{\boldsymbol{u}}}} \right\|}^2}\,{\rm{d}}\mathcal{S}_0^{\rm{D}}}\nonumber\\
	&= \frac{1}{2}{\int_{\mathcal{S}_0^{\rm{L}}} {{\epsilon_{\rm{D}}}\,{\bar \omega}\,{{\left\|{\boldsymbol{u} - {{\bar{\boldsymbol{u}}}}} \right\|}^2}{\mathrm{d}}{\mathcal{S}_0^{\rm{L}}}}},
\end{align}
where $\epsilon_\mathrm{D}>0$ denotes the chosen penalty parameter, and ${\bar\omega}={\bar\omega}(\boldsymbol{X})$ denotes a Heaviside function, defined by 
\begin{equation}\label{gn_bracket}
	{{{\bar\omega}}}\coloneqq \left\{ {\begin{array}{*{20}{c}}
			\begin{array}{l}
				1\\
				0
			\end{array}&\begin{array}{l}
				{\rm{if}}\,\,\,{\boldsymbol{X}} \in {\mathcal{S}_0^\mathrm{D}},\\
				{\rm{if}}\,\,\,{\boldsymbol{X}} \in {\mathcal{S}_0^\mathrm{L}}\setminus{\mathcal{S}_0^\mathrm{D}}.
			\end{array}
	\end{array}} \right.
\end{equation}
Taking the first variation of Eq.\,(\ref{app_pre_disp_pen_func}), and using Eq.\,(\ref{inf_area_init_config}), we have
\begin{equation}
	\label{var_penalty_energy_pre_d}
	\delta {\Pi_\mathrm{D}} = \int_0^L {{\delta {{\boldsymbol{q}}^{\rm{T}}}}{{\bar {\boldsymbol{R}}}_\mathrm{D}}\,{\rm{d}}s},
\end{equation}
where 
\begin{align}
	{\bar {\boldsymbol{R}}}_\mathrm{D} =\int_{{\varXi}^{2}} \frac{\tilde J}{\tilde j}\,{\boldsymbol{\Pi }}{\boldsymbol{f}_\mathrm{D}\,{\rm{d}}{\xi^2}},
\end{align}
with $\boldsymbol{f}_\mathrm{D} \coloneqq {\epsilon_{\rm{D}}}\,{\bar \omega}\,(\boldsymbol{u} - {\bar{\boldsymbol{u}}})$. Taking the directional derivative of Eq.\,(\ref{var_penalty_energy_pre_d}), we obtain
\begin{equation}
	\label{pre_d_direc_deriv_del_pi}
	\Delta \delta {\Pi_\mathrm{D}} = \int_0^L {\delta {{\boldsymbol{q}}^{\rm{T}}}{{\boldsymbol{k}}}_{\rm{D}}\Delta {\boldsymbol{q}}\,{\rm{d}}s},
\end{equation}
where
\begin{align}
	\label{pre_d_k_tild_d_def}
	{{\boldsymbol{k}}}_{\rm{D}} = \int_{\varXi^{2}} \frac{\tilde J}{\tilde j}\,{{\epsilon_{\rm{D}}}\,{\bar \omega}\,{\boldsymbol{\Pi }}{\boldsymbol{\Pi }}^\mathrm{T}{\rm{d}}{{\xi^2}}}.
\end{align}
\subsubsection{Spatial discretization}
Substituting Eq.\,(\ref{nurbs_disc_general_dir_var}) into Eq.\,(\ref{var_penalty_energy_pre_d}) leads to
\begin{equation}
	\delta \Pi _{\rm{D}}^h = \delta {{\bf{q}}^{\rm{T}}}{{\bf{F}}_{\rm{D}}}\,\,\mathrm{with}\,\,{{\bf{F}}_{\rm{D}}} \coloneqq \mathop {\bf{A}}\limits_{e = 1}^{{n_{{\rm{el}}}}} {\bf{F}}_{\rm{D}}^e,
\end{equation}
where the element load vector is defined as
\begin{equation}
	{\bf{F}}_{\rm{D}}^e \coloneqq \int_{\varXi _e} {{{\Bbb{N}}_e^{\rm{T}}}{{\bar{\boldsymbol{R}}}_{\rm{D}}}\,{\tilde j}\,{\rm{d}}\xi^1}.
\end{equation}
Similarly, using (\ref{nurbs_disc_general_dir_var}) into Eq.\,(\ref{pre_d_direc_deriv_del_pi}) gives
\begin{equation}
	\Delta \delta {\Pi^h_\mathrm{D}} = \delta {{\bf{q}}^{\rm{T}}}{{\bf{K}}_{\rm{D}}}\Delta {\bf{q}},
\end{equation}
with ${{\bf{K}}_{\rm{D}}} \coloneqq \mathop {\bf{A}}_{e = 1}^{{n_{{\rm{el}}}}} {\bf{K}}_{\rm{D}}^e$, where the element tangent stiffness matrix is defined as
\begin{equation}
	{\bf{K}}_{\rm{D}}^e \coloneqq \int_{\varXi_e} {{{\Bbb{N}}_e^{\rm{T}}}{{{\boldsymbol{k}}}_{\rm{D}}}{{\Bbb{N}}_e}\,{\tilde j}\,{\rm{d}}\xi^{1}}.
\end{equation}
\section{}
\label{app_algorithm_chart}
\begin{algorithm*}
	\SetAlgoLined
	\For{every Gauss points on the axis of slave body}{
		The convective coordinate $\xi^1_{(1)}$ of the Gauss point is given\;
		Evaluate the position vector $\boldsymbol{\varphi}^{(1)}\equiv\boldsymbol{\varphi}^{(1)}\big(\xi^1_{(1)}\big)$\;
		For the given point $\boldsymbol{\varphi}^{(1)}$, find the closest point on the axis of master body, see Algorithm \ref{proc_algo_ax_gp_ax_prj}\;
		\If{a projection point $\widetilde {\boldsymbol{\varphi}}$ is found within the range of cutoff radius $r_\mathrm{c}$}{
			Determine an initial guess of the convective circumferential coordinate ${\bar \xi^2_{(0)}}$ for the local contact search, see Algorithm \ref{proc_algo_ig_lcsearch}\;
			For an initial guess of the convective axial coordinate, use ${\bar \xi^1_{(0)}}\equiv{\widetilde{\xi^1}}$, see Algorithm\,\ref{proc_algo_ax_gp_ax_prj}\;
			\For{every Gauss point in the boundary of the cross-section at $\xi^1_{(1)}$}{
				Calculate the angle $\theta_\mathrm{G}$ of Eq.\,(\ref{ang_tht_g_cr})\;
				\If{$0\le\theta_\mathrm{G}\le\varepsilon_\theta$}{
					Start the local contact search, see Algorithm\,\ref{proc_algo_lcsearch}\;
				}
			}
		}
	}
	\caption{Overall procedure of the contact search}
	\label{proc_algo_overall}
\end{algorithm*}

\subsection{Global contact search}
The global contact search scheme in Section \ref{glob_csearch} finds contact candidate Gauss integration points on the lateral surface of the slave body. Algorithm\,\ref{proc_algo_overall} shows the overall procedure of the global contact search, and Algorithm\,\ref{proc_algo_ax_gp_ax_prj} presents the procedure of the Gauss point-to-axis closest point projection, presented in Section\,\ref{app_glob_contact_search}. Algorithm \ref{proc_algo_lcsearch} shows the overall process of the local contact search. Since the NURBS basis functions are evaluated in the range of parametric coordinate $\left[\xi^1_\mathrm{min},\xi^1_\mathrm{max}\right]\ni\xi^1$ only, we move to the next initial guess if the coordinate $\xi^1$ goes outside of this range, see lines \ref{xi1_range_treat_s}-\ref{xi1_range_treat_e} of Algorithm\,\ref{proc_algo_ax_gp_ax_prj}, and lines \ref{loc_xi1_range_treat_s}-\ref{loc_xi1_range_treat_e} of Algorithm\,\ref{proc_algo_lcsearch}. We still need to check if the solution of Eq.\,(\ref{cpp_uni_lat}) found by the iterative method is a local maximum solution, based on the condition 
\begin{subequations}
\begin{equation}
	\label{crit_ang_loc_max_cond}
	c_\mathrm{\varphi}>0,
\end{equation}
with
\begin{equation}
	\label{crit_ang_loc_max}
	c_\mathrm{\varphi}\coloneqq\left\{{\boldsymbol{\varphi}^{(1)}}-{\boldsymbol{\varphi}}^{(2)}\big({}^{(i)}{\xi}^1_{(2)}\big)\right\}\cdot\boldsymbol{\nu}_t^{(2)}\big({}^{(i)}\boldsymbol{\xi}^{(2)}\big).
\end{equation}
\end{subequations}
Thus, if the condition of Eq.\,(\ref{crit_ang_loc_max_cond}) is violated, we could try a new initial guess or simply omit the given Gauss integration point. In this paper, we utilize the latter approach (see lines \ref{crit_ang_loc_skip_s}-\ref{crit_ang_loc_skip_e} in Algorithm\,\ref{proc_algo_lcsearch}), since it does not significantly affect the accuracy of the solution, if sufficient number of Gauss integration points are used.
\begin{algorithm*}
	\SetAlgoLined
\KwResult{Closest point projection of a given Gauss integration point on the axis of the slave body to the axis of the master body.}
The position of a Gauss integration point on the axis of slave body $\boldsymbol{\varphi}^{(1)}$ is given\;
Select initial guesses by several points with uniform intervals in the entire domain $\big(\xi^1_\mathrm{min},\xi^1_\mathrm{max}\big)\ni{}^{(0)}{\xi^1_{(2)}}$ of the axis in the master body\;
\For{every initial guess ${}^{(0)}{\xi_{(2)}^1}$}{
	Start an iterative solution process from the initial guess ${}^{(0)}{\xi_{(2)}^1}$\; 
	$n_\mathrm{maxit}$ denotes the chosen limit of the number of iterations\;	
	Initialize the iteration count $i\leftarrow1$\;
	\While{$i\le{n_\mathrm{maxit}}$}{		
		Evaluate $f_\varphi$ of Eq.\,(\ref{ul_dist_prob_nec_cond})\;
		\If{$\left\lvert{f_\varphi}\right\rvert<\varepsilon^\varphi_\mathrm{cpp}$}{
		 	\eIf{$\left\|{\boldsymbol{\varphi}^{(2)}\big({}^{(i)}{\xi^1_{(2)}}\big)-{\boldsymbol{\varphi}}^{(1)}}\right\|\le{r_\mathrm{c}}$}{
		 		$\widetilde{\xi^1}\leftarrow{\prescript{(i)}{}\xi^1_{(2)}}$\;	
		 		$\widetilde{\boldsymbol{\varphi}}\coloneqq\boldsymbol{\varphi}^{(2)}(\widetilde{\xi^1})$\;
		 		Go to line \ref{L_final} (skip the remaining initial guesses)\;
		 	}{
		 		Break (move to the next initial guess);
	 		}
		}
		Calculate the increment $\Delta{\xi^1_{(2)}}$ using Eq.\,(\ref{lcontact_xi2_inc})\;
		Update the solution $\prescript{(i)}{}\xi^1_{(2)}\leftarrow\prescript{(i-1)}{}\xi^1_{(2)} + \Delta\xi^1_{(2)}$\;
		\If{${\prescript{(i)}{}\xi^1_{(2)}}\notin\left[\xi^1_\mathrm{min},\xi^1_\mathrm{max}\right]$}{	\label{xi1_range_treat_s}
			Break (move to the next initial guess);			
		}	\label{xi1_range_treat_e}
		$i\leftarrow{i+1}$\;
	}
}
The solution $\widetilde{\xi^1}$ is used as an initial guess of the convective axial coordinate in the local contact search, see Algorithm \ref{proc_algo_lcsearch}\; \label{L_final}
\caption{Gauss point-to-axis closest point projection}
\label{proc_algo_ax_gp_ax_prj}
\end{algorithm*}

\subsection{Local contact search}
\begin{algorithm*}
	\SetAlgoLined
	\KwResult{Convective coordinates ${\bar{\boldsymbol{\xi}}}=\left[\bar{\xi}^1,\bar{\xi}^2\right]^\mathrm{T}$ of the closest point on the lateral surface of master body for a given surface point $\boldsymbol{x}^{(1)}$ in the slave body.}
	The position of a slave point $\boldsymbol{x}^{(1)}$ on the lateral surface of slave body is given\;
	An initial guess of the convective axial coordinate ${\bar \xi^1_{(0)}}$ is given, see Algorithm \ref{proc_algo_ax_gp_ax_prj}\;
	Select an initial guess of the convective circumferential coordinate ${\bar \xi^2_{(0)}}$, see Algorithm \ref{proc_algo_ig_lcsearch}\;
	Start an iterative solution process to find the convective coordinates ${\bar{\boldsymbol{\xi}}}$ using the initial guess ${\boldsymbol{\xi}}^{(2)}_{(0)}\equiv{\bar{\boldsymbol{\xi}}}_{(0)}=\left[{\bar {\xi}^1_{(0)}},{\bar {\xi}^2_{(0)}}\right]^\mathrm{T}$\;
	$n_\mathrm{maxit}$ denotes the chosen limit of the number of iterations\;
	Initialize the iteration count $i\leftarrow1$\;
	\While{$i\le{n_\mathrm{maxit}}$}{
		Calculate ${{\boldsymbol{f}}\big( {{{{\boldsymbol{\xi }}}^{(2)}_{(i)}}} \big)}$ of Eq.\,(\ref{cpp_unit_lat_a_f})\;
		\If{$\left\| {{\boldsymbol{f}}\big( {{{\boldsymbol{\xi }}_{(i)}^{(2)}}} \big)} \right\| < {\varepsilon _{{\mathrm{cpp}}}}$}{
			Calculate $c_\mathrm{\varphi}$ of Eq.\,(\ref{crit_ang_loc_max})\;
			\eIf{$c_\mathrm{\varphi}<0$}{	\label{crit_ang_loc_skip_s}
				Break (skip the given Gauss point or try a new initial guess)\;	
			}{
				${\bar{\boldsymbol{\xi}}}\leftarrow{{\boldsymbol{\xi}}^{(2)}_{(i)}}$\;
				Break;
			}	\label{crit_ang_loc_skip_e}
		}
		Calculate the increment $\Delta{{\boldsymbol{\xi }}^{(2)}}$ using Eq.\,(\ref{cpp_uni_lat_linearized_sol_inv})\;
		Update the solution using Eq.\,(\ref{cpp_update_add})\;
		Apply the periodicity to $\xi^2_{(2)}$ using Eq.\,(\ref{cyclic_xi2})\;
		%
		\If{$\xi_{(2)}^1\notin\left[\xi^1_\mathrm{min},\xi^1_\mathrm{max}\right]$}{	\label{loc_xi1_range_treat_s}
			Break (skip the given Gauss point or try a new initial guess)\;			
		}	\label{loc_xi1_range_treat_e}
		
		$i\leftarrow{i+1}$\;
	}
	\caption{Local contact search}
	\label{proc_algo_lcsearch}
\end{algorithm*}
\subsubsection{Determination of an initial guess}
\label{lsearch_ig_det}
We determine the initial guess ${\bar \xi^2_{(0)}}$ as the position of the intersection point, shown in Fig.\,\ref{contact_lsearch_ig}, by solving Eq.\,(\ref{intersect_cond_ipd}) using a Newton-Raphson iteration. For a given $\big({}^{(i-1)}{\bar \xi^2_{(0)}},{}^{(i-1)}\alpha_\mathrm{ig}\big)$, we first calculate
\begin{equation}
	\label{inc_xi_alp_det_ig}
	\left\{ {\begin{array}{*{20}{c}}
			{\Delta {\bar \xi^2_{(0)}}}\\
			{\Delta {\alpha _{{\rm{ig}}}}}
	\end{array}} \right\} =  - \prescript{(i-1)}{}{\boldsymbol{e}}^{*\,-1}\left\{ {\begin{array}{*{20}{c}}
			{{}^{(i-1)}{e^1}}\\
			{{}^{(i-1)}{e^2}}
	\end{array}} \right\}
\end{equation}
where ${}^{(i - 1)}{e^\gamma } \coloneqq {e^\gamma }\big({}^{(i-1)}{\bar \xi^{2} _{(0)}},{}^{(i - 1)}\alpha_\mathrm{ig} \big)$, $\gamma\in\left\{1,2\right\}$, and ${}^{(i - 1)}{\boldsymbol{e}^* } \coloneqq {\boldsymbol{e}^*}\big(\prescript{(i-1)}{}{\bar \xi^2_{(0)}}\big)$ with
\begin{align}
	\boldsymbol{e}^*\big(\xi _{(2)}^{2\,}\big)&\coloneqq \left[ \setlength{\arraycolsep}{5pt}
	\renewcommand{\arraystretch}{1.2}{\begin{array}{*{20}{c}}
			{{{\partial {e^1}}}/{{\partial \xi _{(2)}^{2\,}}}}&{{{\partial {e^1}}}/{{\partial \alpha_\mathrm{ig} }}}\\
			{{{\partial {e^2}}}/{{\partial \xi _{(2)}^{2\,}}}}&{{{\partial {e^2}}}/{{\partial \alpha_\mathrm{ig} }}}
	\end{array}} \right] \nonumber\\
	&= \left[ {\begin{array}{*{20}{c}}
			{\zeta _{,2}^{1\,(2)}}&{ - {\boldsymbol{d}}_{(2)}^1 \cdot \widetilde {{{\boldsymbol{\varphi }}_{\rm{d}}}}}\\
			{\zeta _{,2}^{2\,(2)}}&{ - {\boldsymbol{d}}_{(2)}^2 \cdot \widetilde {{{\boldsymbol{\varphi }}_{\rm{d}}}}}
	\end{array}} \right],
\end{align}
and then update
\begin{equation}
	\left.\begin{array}{l}
		\begin{aligned}
		\prescript{(i)}{}{\bar \xi ^2_{(0)}} &= \prescript{(i-1)}{}{\bar \xi_{(0)}^2} + \Delta {\bar \xi _{(0)}^2},\\
		\prescript{(i)}{}{\alpha _{{\rm{ig}}}} &= {}^{(i - 1)}{\alpha _{{\rm{ig}}}} + \Delta {\alpha _{{\rm{ig}}}},
		\end{aligned}
	\end{array}\right\},\,i=1,2,...,
\end{equation}
until a convergence criterion $e < {\varepsilon _{{\rm{ig}}}}$ is satisfied, where 
\begin{equation}
	\label{def_e_det_ig_ints}
	e\coloneqq\sqrt{{\big({}^{(i - 1)}{e^1}\big)^2}+{\big({}^{(i - 1)}{e^2}\big)^2}},
\end{equation}
and $\varepsilon_\mathrm{ig}>0$ is a chosen tolerance. We choose $\prescript{(0)}{}\alpha_\mathrm{ig}=1$, and several values of $\prescript{(0)}{}{\bar \xi^2_{(0)}}\in\left[\xi^2_\mathrm{min},\xi^2_\mathrm{max}\right]$ with a uniform interval, see Algorithm \ref{proc_algo_ig_lcsearch} for the details.
\begin{algorithm*}
	\SetAlgoLined
	\KwResult{Find an intersection point between $\partial \widetilde{\mathcal{A}_t}$ and the projected vector $\widetilde{\boldsymbol{\varphi}_\mathrm{d}}$, which is used as an initial guess in the local contact search.}
	${\bar \xi^1_{(0)}} \leftarrow {\widetilde{\xi^1}}$, see Algorithm \ref{proc_algo_ax_gp_ax_prj}\;
	$\big({\bar \xi^2_{(0)}},\alpha_\mathrm{ig}\big)$ satisfying Eq.\,(\ref{intersect_cond_ipd}) is found by an iterative process\;
	Select several initial guesses ${}^{(0)}{\bar \xi^2_{(0)}}$ of the coordinate $\bar \xi^2_{(0)}$ with uniform interval in the entire domain $\left[\xi^2_\mathrm{min},\xi^2_\mathrm{max}\right]\ni{{}^{(0)}{\bar \xi^2_{(0)}}}$\;
	\For{every initial guess ${}^{(0)}{\bar \xi^2_{(0)}}$}{
		Start an iterative process from the initial guess ${}^{(0)}{\bar \xi^2_{(0)}}$\;
		${}^{(0)}\alpha_\mathrm{ig}\leftarrow1$\;
		$n_\mathrm{maxit}$ denotes the chosen limit of the number of iterations\;
		Initialize the iteration count $i\leftarrow1$\;
		\While{$i\le{n_\mathrm{maxit}}$}{
			Calculate $e$ using Eq.\,(\ref{def_e_det_ig_ints})\;
			\If{$e < {\varepsilon _{{\rm{ig}}}}$}{
				${\bar \xi^2_{(0)}}\leftarrow{}^{(i)}{\bar \xi^2_{(0)}}$\;
				Go to line \ref{finish_det_ig} (skip the remaining initial guesses)\;
			}
			Calculate the increment $\Delta{\bar \xi^2_{(0)}}$ and $\Delta{\alpha}_\mathrm{ig}$ using Eq.\,(\ref{inc_xi_alp_det_ig})\;
			Update the solution ${}^{(i)}{\bar \xi^2_{(0)}}\leftarrow{}^{(i-1)}{\bar \xi^2_{(0)}} + \Delta{\bar \xi^2_{(0)}}$ and ${}^{(i)}\alpha_\mathrm{ig}\leftarrow{}^{(i-1)}\alpha_\mathrm{ig}+\Delta\alpha_\mathrm{ig}$\;
			$i\leftarrow{i+1}$\;
		}
	}
	We use $\big(\bar \xi_{(0)}^1,{\bar \xi_{(0)}^2}\big)$ as an initial guess in the local contact search, see Algorithm \ref{proc_algo_lcsearch}.\label{finish_det_ig}
	\caption{Determination of an initial guess $\big({{\bar \xi^1_{(0)}}},{{\bar \xi^2_{(0)}}}\big)$ in the local contact search}
	\label{proc_algo_ig_lcsearch}
\end{algorithm*}
\section{}
In this appendix, we provide supplementary information in the numerical examples.
\subsection{Sliding contact between two initially straight beams}
\label{app_num_ex_supp_2bslide}
\subsubsection{Case 1}
Table \ref{app_2beam_rigid_flex_ngp_info} shows the selected number of sub-elements for the contact integral. Table \ref{app_2beam_rigid_flex_nload} shows the selected load increment size in each interval of the load parameter for each case of the penalty parameters.
\begin{figure}
	\centering
	\begin{subfigure}[b] {0.475\textwidth} \centering
		\includegraphics[width=\linewidth]{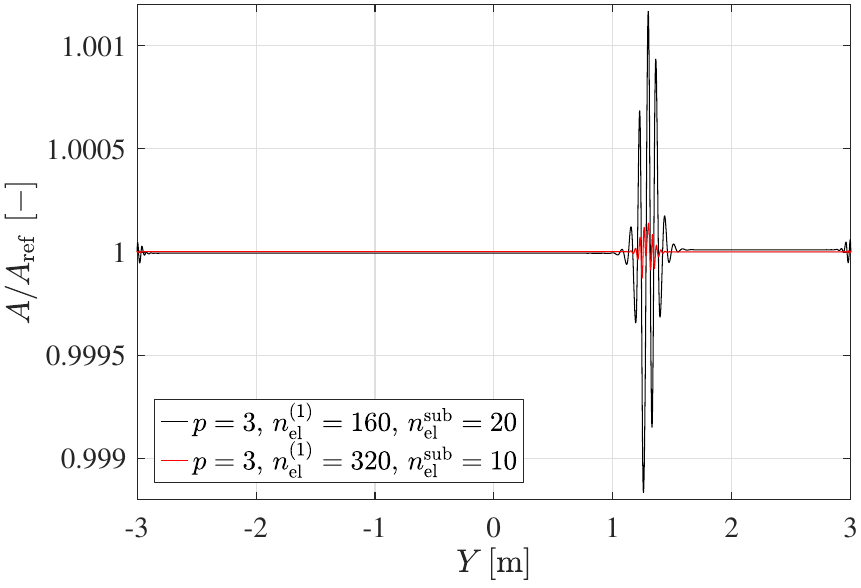}
	\end{subfigure}
	\caption{Contact between two initially straight beams (case 2): Distribution of the cross-sectional area along the axis.}
	\label{cont_slide_carea_org_data}	
\end{figure}
\begin{table}[]
	\scriptsize
	\centering     
	\caption{Sliding contact between two beams (case 1): The selected number of sub-elements in axial ($n_\mathrm{el}^\mathrm{sub}$) and circumferential ($m_\mathrm{el}^\mathrm{sub}$) directions for the contact integral.}
	\label{app_2beam_rigid_flex_ngp_info}	
	\begin{tabular}{ccccccc}
		\Xhline{3\arrayrulewidth}
		$\epsilon_\mathrm{N}{L_0}/E\,[-]$              & 5   & 10  & 20  & 50  & 100 & 200 \\ 
		\Xhline{3\arrayrulewidth}
		$n_\mathrm{el}^\mathrm{sub}$ & 20  & 20  & 20  & 20  & 20  & 40  \\
		$m_\mathrm{el}^\mathrm{sub}$& 100 & 100 & 100 & 100 & 100 & 200 \\ 
		\Xhline{3\arrayrulewidth}
	\end{tabular}
\end{table}
\begin{table*}[]
	\scriptsize
	\centering     
	\caption{Sliding contact between two beams (case 1): Load increment sizes used for each case of penalty parameters in the results of Fig.\,\ref{btb_slide_case1_conv_penet_cforce}.}
	\label{app_2beam_rigid_flex_nload}	
	\begin{tabular}{cccccccc}
		\Xhline{3\arrayrulewidth}
		\multicolumn{2}{c}{$\epsilon_\mathrm{N}{L_0}/{E}\,\,[-]$}& 5                      & 10                     & 20                     & 50     & 100   & 200    \\ 
		\Xhline{3\arrayrulewidth}
		\multirow{2}{*}{$\Delta \lambda_\mathrm{load}$}&${\lambda_\mathrm{load}}\in[{0},{0.5}]$                     & {0.01} & {0.01} & {0.01} & 0.01   & 0.01  & 0.0025 \\
		&${\lambda_\mathrm{load}}\in[{0.5},{1}]$                     &                  0.005      &           0.005             &          0.005              & 0.0025 & 0.001 & 0.001  \\ \hline
		\multicolumn{2}{c}{Total \#load   steps}  & 150                    & 150                    & 150                    & 250    & 550   & 700    \\
		\Xhline{3\arrayrulewidth}		
	\end{tabular}
\end{table*}
\subsubsection{Case 2}
Fig.\,\ref{cont_slide_carea_org_data} shows the original graph of the magnified one in Fig.\,\ref{cont_slide_carea_ratio_mag_view}. Table \ref{app_2beam_flex_flex_nload} shows the selected load increment size in each interval of the load parameter.
\begin{table*}[]
	\scriptsize
	\centering     
	\caption{Sliding contact between two beams (case 2): Load increment sizes used for each case of the results in Figs.\,\ref{cont_slide_cent_ydisp_master} and \ref{cont_slide_flex_carea_ratio}.}
	\label{app_2beam_flex_flex_nload}	
	\begin{tabular}{ccccc}
		\Xhline{3\arrayrulewidth}	
		\multicolumn{1}{c}{}           & \multicolumn{1}{c}{} & \multicolumn{1}{c}{\multirow{3}{*}{\begin{tabular}[c]{@{}l@{}}$p=3$, $n_\mathrm{el}^{(1)}=160$,\\ $n_\mathrm{el}^\mathrm{sub}=20$, $n_\mathrm{el}^{(2)}=160$\end{tabular}}} & \multicolumn{1}{c}{\multirow{3}{*}{\begin{tabular}[c]{@{}l@{}}$p=3$, $n_\mathrm{el}^{(1)}=320$,\\$n_\mathrm{el}^\mathrm{sub}=10$, $n_\mathrm{el}^{(2)}=160$\end{tabular}}} & \multicolumn{1}{c}{\multirow{3}{*}{\begin{tabular}[c]{@{}l@{}}$p=4$, $n_\mathrm{el}^{(1)}=320$,\\ $n_\mathrm{el}^\mathrm{sub}=20$, $n_\mathrm{el}^{(2)}=320$\end{tabular}}}\\
		       &     &   &    &                      \\ 
				&	&	&	&	\\
		\Xhline{3\arrayrulewidth}
		\multirow{2}{*}{$\Delta{\lambda}_\mathrm{load}$} & ${\lambda_\mathrm{load}}\in[{0},{0.5}]$                & 0.005                                                                                     & 0.002                                                      & {0.002}                               \\
		& ${\lambda_\mathrm{load}}\in[{0.5},{1}]$                & 0.0025                                                                                    & 0.001                                                                                     & 0.001  \\ \hline
		\multicolumn{2}{c}{Total \#load steps}                & 300                                                                                       & 750                                                                                       & 750\\ 
		\Xhline{3\arrayrulewidth}		
	\end{tabular}
\end{table*}
\subsection{Twisting of wire strands}
\subsubsection{A strand of two wires}
Table \ref{app_twist2b_nlstep_info_r025_r0125_r00625} shows the selected number of sub-elements for the contact integral. Tables \ref{app_twist2b_nlstep_info_r025}, \ref{app_twist2b_nlstep_info_r0125}, and \ref{app_twist2b_nlstep_info_r00625} show the chosen load increment sizes in each case of the initial radii $R=0.25\,\mathrm{m}$, $0.125\,\mathrm{m}$, and $0.0625\,\mathrm{m}$, respectively.
\begin{table}[]
	\scriptsize
	\centering     
	\caption{Twisting of wire strands (two wires, $R=0.25\,\mathrm{m},\,0.125\,\mathrm{m},\,0.0625\,\mathrm{m}$): The selected number of sub-elements in axial $(n^\mathrm{sub}_\mathrm{el})$ and circumferential directions $(m^\mathrm{sub}_\mathrm{el})$ for the contact integral.}
	\label{app_twist2b_nlstep_info_r025_r0125_r00625}	
	\begin{tabular}{clcccc}
		\Xhline{3\arrayrulewidth}		
		\multicolumn{2}{c}{$\epsilon_\mathrm{N}{L_0}/E\,\,[-]$}                                                          & $10^0$ & $10^1$ & $10^2$                     & $10^3$                     \\ 
		\Xhline{3\arrayrulewidth}
		& $R=0.25\,\mathrm{m}$   & 10     & 10     & 10                         & 20                         \\
		& $R=0.125\,\mathrm{m}$  & 10     & 10     & 10                         & 20                         \\
		\multirow{-3}{*}{$n^\mathrm{sub}_\mathrm{el}$}                                                     & $R=0.0625\,\mathrm{m}$ & 10     & 10     & {20}  & {30}  \\ \hline
		& $R=0.25\,\mathrm{m}$   & 200    & 200    & 200                        & 300                        \\
		& $R=0.125\,\mathrm{m}$  & 200    & 200    & 200                        & 300                        \\
		\multirow{-3}{*}{$m^\mathrm{sub}_\mathrm{el}$} & $R=0.0625\,\mathrm{m}$ & 200    & 200    & {300} & {400} \\ 
		\Xhline{3\arrayrulewidth}
	\end{tabular}
\end{table}
\begin{table*}[]
	\scriptsize
	\centering
	\caption{Twisting of wire strands (two wires, $R=0.25\,\mathrm{m}$): Load increment sizes used for each case of penalty parameters in the results of Fig.\,\ref{twist_2beam_caxis_diff_ref}.}
	\label{app_twist2b_nlstep_info_r025}
	\begin{tabular}{clcccc}
		\Xhline{3\arrayrulewidth}		
		\multicolumn{2}{c}{$\epsilon_\mathrm{N}{L_0}/E\,[-]$}  & $10^0$                      & $10^1$                     & $10^2$                    & $10^3$                     \\
		\Xhline{3\arrayrulewidth}		
		\multirow{4}{*}{$\Delta{\lambda}_\mathrm{load}$} & ${\lambda_\mathrm{load}}\in[{0},{0.3}]$   & \multirow{4}{*}{$0.005$} & \multirow{4}{*}{$0.005$} & {$0.0025$}                & 0.00125 \\
		& ${\lambda_\mathrm{load}}\in[{0.3},{0.5}]$ &                        &                        & 0.0025 &        0.00125                  \\
		& ${\lambda_\mathrm{load}}\in[{0.5},{0.95}]$ &                        &                        & 0.001                       & $0.001$                   \\
		& ${\lambda_\mathrm{load}}\in[{0.95},{1}]$   &                        &                        &  0.001                      & $0.0001$                  \\ \hline
		\multicolumn{2}{c}{Total \#load steps}   & 200                    & 200                    & 700                    & 1350                     \\
		\Xhline{3\arrayrulewidth}		
	\end{tabular}
\end{table*}
\begin{table*}[]
	\scriptsize
	\centering
	\caption{Twisting of wire strands (two wires, $R=0.125\,\mathrm{m}$): Load increment sizes used for each case of penalty parameters in the results of Fig.\,\ref{twist_2beam_caxis_diff_ref}.}
	\label{app_twist2b_nlstep_info_r0125}	
	\begin{tabular}{clcccc}
		\Xhline{3\arrayrulewidth}		
		\multicolumn{2}{c}{$\epsilon_\mathrm{N}{L_0}/E\,[-]$}  & $10^0$                       & $10^1$                      & $10^2$                    & $10^3$   \\ 
		\Xhline{3\arrayrulewidth}
		\multirow{4}{*}{$\Delta{\lambda}_\mathrm{load}$} & ${\lambda_\mathrm{load}}\in[{0},{0.5}]$   & 0.005                   & 0.005                   & 0.0025                 & 0.0025 \\
		& ${\lambda_\mathrm{load}}\in[{0.5},{0.8}]$ & 0.0025  & 0.0025 & 0.001 & 0.001  \\
		& ${\lambda_\mathrm{load}}\in[{0.8},{0.9}]$   &  0.0025  & 0.0025                        &  0.001                      & 0.0005 \\ 
		& ${\lambda_\mathrm{load}}\in[{0.9},{1}]$   &  0.0025  & 0.0025                        &  0.001                      & 0.00005 \\ \hline		
		\multicolumn{2}{c}{Total \#load steps}   & 300                     & 300                     & 700                    & 2700   \\
		\Xhline{3\arrayrulewidth}		
	\end{tabular}
\end{table*}
\begin{table*}[]
	\scriptsize
	\centering
	\caption{Twisting of wire strands (two wires, $R=0.0625\,\mathrm{m}$): Load increment sizes in each case of penalty parameters in the results of Fig.\,\ref{twist_2beam_caxis_diff_ref}.}
	\label{app_twist2b_nlstep_info_r00625}
	\begin{tabular}{cccccc}
		\Xhline{3\arrayrulewidth}		
		\multicolumn{2}{c}{$\epsilon_\mathrm{N}{L_0}/E\,[-]$}   & $10^0$                      & $10^1$                     & $10^2$                    & $10^3$   \\ 
		\Xhline{3\arrayrulewidth}		
		\multirow{5}{*}{$\Delta{\lambda}_\mathrm{load}$} & ${\lambda_\mathrm{load}}\in[{0},{0.3}]$    & 0.005  & 0.005  & 0.0025 & 0.0025  \\
		& ${\lambda_\mathrm{load}}\in[{0.3},{0.5}]$  & 0.005  & 0.005  & 0.0025 & 0.0025  \\
		& ${\lambda_\mathrm{load}}\in[{0.5},{0.9}]$  & 0.0025 & 0.0025 & 0.0005 & 0.0005  \\
		& ${\lambda_\mathrm{load}}\in[{0.9},{0.93}]$ & 0.0025 & 0.0025 & 0.0005 & 0.0001  \\
		& ${\lambda_\mathrm{load}}\in[{0.93},{1}]$   & 0.0025 & 0.0025 & 0.0005 & 0.00005 \\ \hline
		\multicolumn{2}{c}{Total \#load steps}    & 300    & 300    & 1200   & 2700    \\ 
		\Xhline{3\arrayrulewidth}		
	\end{tabular}
\end{table*}
\subsubsection{A strand of seven wires}
Table \ref{app_twist7b_pair_info_m6s1a} shows the information of the chosen slave-master contact pairs for the case of M6S1A in the example of seven wires. Table \ref{app_twist7b_nlstep_info_2case} shows the selected load increment sizes in each case of M1S6 and M6S1A.
\begin{table}[]
	\small
	\centering
	\caption{Twisting of wire strands (seven wires): The information of slave-master contact pair in the case of M6S1A. See Fig.\,\ref{twist_7b_init_yz} for the numbering of bodies.}
	\label{app_twist7b_pair_info_m6s1a}
	\bgroup
	\def\arraystretch{0.95}
\begin{tabular}{cccc}
	\Xhline{3\arrayrulewidth}
	\multicolumn{1}{l}{}                                                                             & \begin{tabular}[c]{@{}c@{}}Contact\\ pair\#\end{tabular} & \begin{tabular}[c]{@{}c@{}}Slave\\ body\end{tabular} & \begin{tabular}[c]{@{}c@{}}Master\\ body\end{tabular} \\ 
	\Xhline{3\arrayrulewidth}
	\multirow{6}{*}{\begin{tabular}[c]{@{}c@{}}Inner wire\\ -outer wire\\      contact\end{tabular}} & 1                                                        & 7                                                    & 1                                                     \\
	& 2                                                        & 7                                                    & 2                                                     \\
	& 3                                                        & 7                                                    & 3                                                     \\
	& 4                                                        & 7                                                    & 4                                                     \\
	& 5                                                        & 7                                                    & 5                                                     \\
	& 6                                                        & 7                                                    & 6                                                     \\ \hline
	\multirow{6}{*}{\begin{tabular}[c]{@{}c@{}}Outer wire\\ -outer wire\\      contact\end{tabular}} & 7                                                        & 1                                                    & 2                                                     \\
	& 8                                                        & 2                                                    & 3                                                     \\
	& 9                                                        & 3                                                    & 4                                                     \\
	& 10                                                       & 4                                                    & 5                                                     \\
	& 11                                                       & 5                                                    & 6                                                     \\
	& 12                                                       & 6                                                    & 1                                                     \\ 
	\Xhline{3\arrayrulewidth}
\end{tabular}
\egroup
\end{table}
\begin{table}[]
	\scriptsize
	\centering
	\caption{Twisting of wire strands (seven wires): Load increment sizes.}
	\label{app_twist7b_nlstep_info_2case}		
	\begin{tabular}{clcc}
		\Xhline{3\arrayrulewidth}
		\multicolumn{2}{c}{}                   & M1S6   & M6S1A  \\ 
		\Xhline{3\arrayrulewidth}
		\multirow{4}{*}{$\Delta \lambda_\mathrm{load}$} & ${\lambda_\mathrm{load}}\in{[0,0.4]}$ & 0.005  & 0.004  \\
		& ${\lambda_\mathrm{load}}\in{[0.4,0.5]}$ & 0.005 & 0.0004 \\
		& ${\lambda_\mathrm{load}}\in{[0.5,0.8]}$ & 0.001 & 0.0004 \\
		& ${\lambda_\mathrm{load}}\in{[0.8,1]}$ & 0.0005 & 0.0004 \\ \hline
		\multicolumn{2}{c}{Total \#load steps} & 800   & 1600   \\ 
		\Xhline{3\arrayrulewidth}
	\end{tabular}
\end{table}

\end{appendices}

\bibliography{sn-bibliography}


\end{document}